\documentclass[11pt]{amsart}
\usepackage[utf8]{inputenc}
\usepackage{amsmath, amssymb, amsthm,mathtools}
\usepackage{hyperref}
\usepackage{geometry}
\usepackage{enumerate}
\usepackage{tikz-cd}
\usepackage{fullpage}
\usepackage{dsfont}
\usepackage{babel}
\usepackage{nicefrac}

\newtheorem{theorem}{Theorem}[section]
\newtheorem*{theorem*}{Theorem}
\newtheorem{lemma}[theorem]{Lemma}
\newtheorem{proposition}[theorem]{Proposition}
\newtheorem{corollary}[theorem]{Corollary}
\newtheorem{conjecture}[theorem]{Conjecture}

\theoremstyle{definition}
\newtheorem{definition}[theorem]{Definition}
\newtheorem*{definition*}{Definition}
\newtheorem{remark}[theorem]{Remark}
\newtheorem{example}[theorem]{Example}

\numberwithin{equation}{section}

\newcommand{\A}{\mathbb{A}}
\newcommand{\Al}{\mathbb{A}^1}

\newcommand{\Z}{\mathbb{Z}}
\newcommand{\Q}{\mathbb{Q}}
\newcommand{\C}{\mathbb{C}}
\newcommand{\R}{\mathbb{R}}
\newcommand{\N}{\mathbb{N}}
\newcommand{\T}{\mathbb{T}}
\newcommand{\KK}{\mathbb{K}}
\newcommand{\KKR}{\mathbb{K}_{\R}}
\newcommand{\Spec}{\operatorname{Spec}}
\newcommand{\Hom}{\mathrm{Hom}}

\newcommand{\GW}{\operatorname{GW}}
\newcommand{\W}{\operatorname{W}}
\newcommand{\ind}{\operatorname{ind}}
\newcommand{\Tr}{\operatorname{Tr}}
\newcommand{\gw}[1]{\left\langle #1 \right\rangle}
\newcommand{\val}{\operatorname{val}}
\newcommand{\G}{\mathbb{G}}
\newcommand{\NP}{\operatorname{NP}}
\newcommand{\DS}{\operatorname{DS}}
\newcommand{\mce}{\mathcal{E}}
\newcommand{\mcx}{Y}

\newcommand{\mcl}{\mathcal{L}}
\newcommand{\PP}{\mathbb{P}}

\newcommand{\trop}{\operatorname{trop}}
\newcommand{\mult}{\operatorname{mult}}

\newcommand{\Div}{\operatorname{Div}}
\newcommand{\Divo}{\operatorname{Div}_0}
\DeclareMathOperator{\divi}{div}

\newcommand{\mcOX}{\mathcal{O}_{\XSig}}

\newcommand{\bfx}{\mathbf{x}}

\newcommand{\wivasig}{\widetilde{\varphi}_\sigma}
\newcommand{\mct}{\mathcal{T}}
\DeclareMathOperator{\Cl}{Cl}
\DeclareMathOperator{\Gl}{Gl}

\newcommand{\XSig}{X_\Sigma}
\newcommand{\bso}{b_\sigma^{\omega}}
\newcommand{\bst}{b_\tau^{\omega}}
\newcommand{\diff}[2][]{\frac{\partial #1}{\partial #2}}
\newcommand{\dsso}[1][\sigma]{d_{\sigma_0}^#1}
\DeclareMathOperator{\diag}{diag}
\DeclareMathOperator{\MV}{MV}
\DeclareMathOperator{\vol}{vol}
\DeclareMathOperator{\conv}{conv}

\DeclareMathOperator{\charact}{char}
\DeclareMathOperator{\Area}{Area}

\DeclareMathOperator{\Jac}{Jac}
\DeclareMathOperator{\spann}{span}

\newcommand{\Homalg}{\Hom_{\KK\text{-}\text{alg}}}

\newcommand{\Homgp}{\Hom_{\text{gp}}}

\newcommand{\kk}{\mathds{k}}
\newcommand{\multA}{\mult^{\Al}}
\newcommand{\GWk}{\GW(\kk)}

\newcommand{\kkq}[1][\kk]{\nicefrac{#1^\times}{(#1^\times)^2}}
\newcommand{\kkqo}[1][\kk]{\nicefrac{#1}{(#1^\times)^2}}

\newcommand{\psisd}[1][D]{\psi_\sigma^{#1}}

\title{Combinatorial relative orientability for toric and tropical geometry}

\author{Tiggemann, Nathan}

\address{Department of Mathematics, TU Darmstadt}

\email{tiggemann@mathematik.tu-darmstadt.de}

\date{\today}

\begin{document}
	\maketitle
\begin{abstract}
	Many problems in enumerative geometry are solved by computing the degree of the top Chern class of a suitable vector bundle. Results from $\Al$-homotopy theory allow to answer enumerative questions over arbitrary fields, replacing the top Chern class by the Euler class in Chow-Witt groups. For its degree to be well-defined, one needs the vector bundle to be relatively oriented.\\
	Another way to solve enumerative problems is to translate them to tropical geometry, where they become combinatorial problems. In this paper, we give a combinatorial description of relative orientability for vector bundles over smooth proper toric varieties. We show that it behaves well under tropicalization, establishing a tropical version of relative orientability, which will allow to the use of tropical geometry to answer enumerative questions over arbitrary fields. As a demonstration, we prove a quadratically enriched version of the Bernstein-Khovanskii-Kushnirenko theorem via tropical geometry. Lastly, we define the notion of $\kk$-phase structures on tropical hypersurfaces, generalizing the real case.
\end{abstract}
\tableofcontents
\section{Introduction}
	In $\Al$\emph{-enumerative geometry}, one poses questions in enumerative geometry over an arbitrary field $\kk$. Here, the answers take values in the Grothendieck-Witt ring $\GW(\kk)$ of $\kk$. In particular, one can also solve real counting problems.\\
	One way to tackle enumerative questions is via \emph{tropical geometry}. It allows counting problems from enumerative geometry to be translated to a tropical counterpart where one can use combinatorial methods to solve the problem. This was successfully generalized to $\Al$-counts e.g.\@ by Jaramillo Puentes and Pauli in \cite{JPPbezout} and Markwig, Payne, and Shaw in \cite{MarkwigPayneShaw}.\\
	In classical enumerative geometry, many problems can be solved by computing the degree of the top Chern class of an appropriate vector bundle $V \to X$. In $\A^1$-enumerative geometry, the top Chern class is replaced by the Euler class in twisted Chow-Witt groups. For its degree $n^{\Al}$ to be well defined, the vector bundle $V \to X$ must be relatively oriented in the following sense.
	\begin{definition*}
		A rank $n$ vector bundle $V \to X$ over a smooth $n$-dimensional $\kk$-scheme is \emph{relatively orientable} if there exists a line bundle $\mathcal{L} \to X$ and an isomorphism  
		\[
		\rho\colon \det V \otimes \omega_{X/\kk} \xrightarrow{\;\cong\;} \mathcal{L}^{\otimes 2},
		\]  
		where $\omega_{X/\kk}$ denotes the canonical line bundle. In this case, we say that $V$ is \emph{relatively oriented by $(\mathcal{L},\rho)$}.
	\end{definition*}
	In their paper, Jaramillo Puentes and Pauli used tropical geometry to address the following problem in $\Al$-enumerative geometry, which is known as the Bernstein-Khovanskii-Kushnirenko (BKK) theorem in the classical case.\\
	Let $X = X_\Sigma$ be a smooth, proper toric variety over $\kk$ of dimension $n$. For $i = 1, \dots, n$, let $V_i \subset X$ be hypersurfaces, i.e., the zero loci of global sections of line bundles $\mathcal{L}_i$ on $X$. How many points are there in $V_1\cap \dots \cap V_n$?\\
	It turns out that, assuming the relative orientability of $V = \mcl_1\oplus \dots \oplus \mcl_n\to X$, one can count the intersections $x\in V_1\cap \dots\cap V_n$ with a quadratically enriched multiplicity $\multA(x)$ valued in $\GWk$. The resulting sum is a well-defined element in $\GWk$ independent of the specific sections of the line bundles, as it equals the degree of the Chow-Witt-valued Euler class $e(V)$ \cite{KassWickelgrenCubicSurface}\cite{LevineAspects}\cite{BachmannWickelgren}.\\
	Here is how one can use tropical geometry to compute the degree of $e(V)$ valued in $\GWk$.
	The intersection points of the $V_i$ tropicalize to intersection points of their tropicalizations $\trop(V_i)$. \cite{JPPbezout} find a purely combinatorial rule assigning a $\GWk$-valued multiplicity $\operatorname{mult}^{\A^1}(p)$ to each point $p \in \trop(V_1) \cap \dots \cap \trop(V_n)$ which equals the sum of quadratically enriched intersection multiplicities $\mult^{\A^1}(x)$ of points $x\in V_1\cap\dots\cap V_n$ with $\trop x=p$.\\
	To define the tropical enriched intersection multiplicity, the authors decorate tropical hypersurfaces by placing certain coefficients in its complement. We will present an alternative way to remember quadratic data, $\kk$-phase structures. These directly generalize real phase structures to arbitrary fields for tropical hypersurfaces. The definition of $\kk$-phases is a first step towards a possible generalization of tropical enriched intersection multiplicity to not-$0$-dimensional intersections.\\
	In \cite{JPPbezout} it is shown that the sum of the enriched tropical intersection multiplicities is independent of the chosen hypersurfaces, provided that the Minkowski sum of the Newton polytopes of their defining polynomials has no boundary lattice point with all coordinates odd.\\
	This inspires the authors to define this property as the combinatorial version of relative orientability. Phrased for one hypersurface only:
	\begin{definition*}
		A tropical hypersurface is \emph{combinatorially orientable} if the Newton polytope of its defining polynomial has no boundary points in $\Z^n$ with all entries odd.
	\end{definition*}
	In \cite[Conjecture 5.14.]{JPPbezout} it is conjectured that combinatorial orientability of the tropicalization of the hypersurfaces is equivalent to the relative orientability of the sum of the line bundles.\\
	However, this definition has two flaws: First, it depends on the exact polynomial defining the hypersurface. For example both tropical polynomials $x\oplus y\oplus 0$ and $x^2\oplus xy \oplus x$ define the same plane tropical curve, but only one of them has a Newton polytope with no odd boundary points. Second, for general finite rank lattices, there is no well-defined notion of odd points, so this definition can not directly be generalized.\\
	In this paper, we give an alternative definition of combinatorial orientability that does not depend on the polynomial and generalizes to arbitrary finite rank lattices. To motivate the definition, we will translate the relative orientability of a torus invariant line bundle $\mcOX(D)$, that is, $\mcOX(D)\otimes \omega_{X_\Sigma/\kk}$ being isomorphic to a tensor square of a line bundle, on a smooth proper toric variety to a combinatorial one on the equations defining the associated polytope $P_D$. Notice that considering line bundles is sufficient for the general case, as $\det V \otimes \omega_{X/\kk}$ is a line bundle anyways. For that, let $\pi_2\colon M\to M/2M$ be the reduction modulo $2$. Further, let $\Div_T(\XSig)$ be the group of torus invariant divisors. 
	\begin{theorem*}[\ref{thm:DrelOrientInLattice}]
		Let $\XSig$ be a smooth proper toric variety and $D=\sum_{\rho}a_\rho D_\rho\in\Div_T(\XSig)$. Then $\mcOX(D)$ is relatively orientable, i.e. $\mcOX(D)\otimes \omega_{X_\Sigma/\kk}$ isomorphic to a square of a line bundle, if and only if $\pi_2({L_D\cap M})\not = M/2M$, where $L_D = \{m\in M_\R\mid \langle m,v_\rho \rangle = -a_\rho \text{ for some }\rho \in \Sigma(1)\}$.
	\end{theorem*}
	This inspires the following definition.	
	\begin{definition*}[\ref{def:combOrientP}]
		We say that a polytope $P \subset M_\R$ is \emph{combinatorially orientable} if $P$ is full-dimensional and $\pi_2(L_P\cap M)\not = M/2M$, where $L_P$ is the union of the supporting hyperplanes of the minimal description of $P$ as an intersection of half-spaces.
	\end{definition*}
	Further, we will answer the conjecture mentioned above, building a bridge between classical relative orientability and the combinatorial one.
	\begin{theorem*}[\ref{thm:JPPforNef}]
		Let $\XSig$ be a smooth proper toric variety of dimension $n$ and let $\mcOX(E_1)$, $\ldots$, $\mcOX(E_n)$ be toric line bundles on $\XSig$ and all $E_i$ be nef. Further, let $A_i \coloneqq P_{E_i}\cap M\subset M$.\\
		If the vector bundle $\mcOX(E_1)\oplus \dots \oplus\mcOX(E_n)\to \XSig$ is relatively orientable, and $P_{E_1}+\dots+P_{E_n}$ of full dimension, then $\conv(A_1+\dots+A_n)$, the convex hull of the Minkowski sum of $A_1,\dots ,A_n$, is combinatorially orientable.
	\end{theorem*}
	We will give counterexamples for not-nef divisors. Starting with a combinatorially orientable polytope like in the case of the BKK theorem, we will prove that one is always able to find a smooth proper toric variety for which the polytope comes from a relatively orientable divisor.
	\begin{theorem*}[\ref{thm:combOrientGivesOrientDivisor}]
		Let $P$ be a combinatorially orientable lattice polytope. Then, for every smooth refinement $\Sigma$ of the normal fan of $P$, there exists a relatively orientable $D\in\Div_T(\XSig)$ such that $P = P_D$.
	\end{theorem*}
	Lastly, connecting the new definition to the one in \cite{JPPbezout} using odd boundary points, we get the following result in the charts $U_\sigma$, indexed by the maximal cones. It gives a sufficient criterion for the inverse of \cite{JPPbezout}'s conjecture to hold. $L_{D,\sigma}$ and $\NP_{\sigma}(D)$ are in-chart versions of the general constructions.
	\begin{theorem*}[\ref{thm:orientAndOddBoundary}]
		Let $\XSig$ be a smooth proper toric variety and $D\in\Div_T(\XSig)$ a torus invariant divisor. Then the following statements are true.
		\begin{itemize}
			\item $D$ is relatively orientable if and only if  $L_{D,\sigma}$ contains no odd points for some/all $\sigma\in\Sigma(n)$.
			\item If $P_D$ is a lattice polytope, $D$ being relatively orientable implies that $\NP_{\sigma}(D)$ has no odd boundary points for some/all $\sigma\in\Sigma(n)$.
			\item If further all edges of $P_D$ have lattice length $\geq n-1$ and $D$ is ample, then $\NP_{\sigma}(D)$ having no odd boundary points for some  $\sigma\in\Sigma(n)$ implies that $D$ is relatively orientable.
		\end{itemize}
		In particular, if $D$ is the sum of $n$ ample divisors, $D$ is relatively orientable if and only if $\NP_{\sigma}(D)$ has no odd boundary points for some/all $\sigma\in\Sigma(n)$		
	\end{theorem*}
	For the tropical case, this implies the following.
	\begin{theorem*}[\ref{thm:tropicalizeOrient}]
		Let $\XSig$ be a smooth proper toric variety and $D\in \Div_T(\XSig)$ be relatively orientable, i.e. $\mcOX(D)\otimes \omega_{X_\Sigma/\kk}$ a tensor square, with $P_D$ a full-dimensional lattice polytope.	
		Then, given a general section $f$ of $\mcOX(D)$, the Newton polytope of its tropicalization is combinatorially orientable. Further, in each chart $U_\sigma$, its Newton polytope has no odd boundary points.\\
		In particular, \cite{JPPbezout}'s Conjecture 5.14 is true if $D$ is the sum of $n$ ample divisors. 
	\end{theorem*}
	To give an application, we will give a new combinatorial tropical proof of the $\Al$-BKK theorem using our results. We do not claim originality of this result, it follows e.g.\@ from \cite[Proposition 19]{SrinivasanWickelgren}, which is based on work in \cite{LevineAspects}, together with the classical result.
	\begin{theorem*}[\ref{thm:A1BKK}, $\Al$-BKK theorem]
		Let $\Sigma$ be a smooth complete fan of dimension $n$ and let
		$E_1,\dots,E_n\in\Div_T(X_\Sigma)$ be nef such that
		$V = \mathcal{O}(E_1)\oplus\dots\oplus\mathcal{O}(E_n)$ is relatively
		orientable. Then for every field $\kk$
		$$n^{\Al}(V_\kk) = \frac{\MV(P_{E_1},\dots,P_{E_n})}{2}\,h\in\GW(\kk).$$
		In particular, if $\kk$ is infinite, then for a general section $s$ of $V_\kk$
		all zeros are simple and lie in the torus, and
		$\sum_{z} \ind_z s = \frac{\MV(P_{E_1},\dots,P_{E_n})}{2}\,h$.
	\end{theorem*}
\subsection{Overview}
\label{subsection:Overview}
	Section \ref{section:TorVarsAndFans} gives a brief introduction to toric varieties and the Grothendieck-Witt ring.
	In the next section, Section \ref{section:CombRelOrient}, we discuss how far it is possible to translate relative orientability into a combinatorial property, give a new definition of combinatorial orientability and discuss its relation to \cite{JPPbezout}'s combinatorial orientability. 
	In Section \ref{section:ExplicitRelOrientOnTorVar}, we construct explicit relative orientations for torus invariant vector bundles and prove a formula for the $\Al$-multiplicities of zeros of sections. Section \ref{section:CombRelOrientInTropGeo} gives a quick introduction to (quadratically enriched) tropical geometry, the original inspiration for this paper. We discuss the implications of the results in the previous sections for tropical geometry and we use them to (re-)prove the $\Al$-BKK theorem. The last Section \ref{section:kphase} establishes an alternative way to enrich tropical hypersurfaces via $\kk$-phase structures, a generalization of real phase structures.
\subsection{Acknowledgments}
	The author acknowledges support by Deu\-tsche For\-schungs\-ge\-mein\-schaft (DFG, German Research Foundation) through the Collaborative Research Centre TRR 326 \textit{Geometry and Arithmetic of Uniformized Structures}, project number 444845124.\\
	Further, the author wants to thank Sabrina Pauli for suggesting this project, proof reading this paper and answering many questions.
\subsection{Declaration of AI usage}
	For this paper, AI (Claude Opus 5.5) was used to check correctness, to find some examples, and to simplify some proofs and statements. Every line in this text was written by hand.
\section{Preliminaries}
\label{section:TorVarsAndFans}
\subsection{Toric varieties}
\label{subsection:TorVars}
	In this section, we fix the notation we use in the context of toric varieties. For an introduction to toric geometry, see e.g.\@ \cite{cox2011toric}.\\
	From now on, $\XSig$ will always be a smooth complete toric variety over a field $\kk$, with a split torus $T$ acting on it, defined by a convex rational polyhedral fan $\Sigma$ in $N_\R\eqqcolon N \otimes \R$ where $N$ is a rank $n$ lattice with dual lattice $M\coloneqq \Hom(N,\Z)$. 
	\begin{example}
		$\PP^n = \XSig$ for $\Sigma$ the smooth complete fan with rays defined by $v_0 = -\sum_{i=1}^{n}e_i$ and $v_i=e_i$ for $i=1,\dots, n$ and maximal cones $\sigma_i$ defined by all rays except the $i$-th one.
	\end{example}
	The set of all $r$-dimensional cones in $\Sigma$ will be denoted by $\Sigma(r)$. For $\sigma\in\Sigma$ we will write $U_\sigma$ for the corresponding affine variety $\Spec \kk[\sigma^\vee \cap M]$ where $\sigma^\vee$ denotes the dual cone. We order the rays of $\Sigma$ and write $\Sigma(1) = \{\rho_1,\dots,\rho_s\}$, further, $\sum_{\rho}$ means the sum over all $\rho\in\Sigma(1)$. If we choose a subset $\{\rho_{i_1},\dots\rho_{i_m}\}$ of the rays, we always assume that $i_1<\dots<i_m$. Given a ray $\rho = \rho_i$, its unique minimal generator will be denoted $v_\rho = v_i$ and the corresponding torus invariant prime divisor $D_\rho = D_i$. The sets of Weil divisors, torus-invariant Weil divisors, and principal divisors of $\XSig$ are $\Div(\XSig),\Div_T(\XSig)$ and $\Divo(\XSig)$, respectively. Given a torus invariant divisor $D=\sum_{\rho}a_\rho D_\rho$ we write $\{b_\sigma^D\}_{\sigma\in\Sigma(n)}$ for its Cartier data, i.e.\@ for the collection of elements in $M$ for which $\langle b_\sigma^D, v_\rho\rangle = -a_\rho$ for every $\rho\in\sigma(1)$. Explicitly, $b_\sigma^D = -\sum_{m\in B_\sigma} a_\rho m$.\\
	Since $\Sigma$ is smooth and complete, the minimal generators $\{v_{i_1},\dots,v_{i_n}\} = \{v_\rho\}_{\rho\in\sigma(1)}$ of the rays of a maximal cone $\sigma\in\Sigma(n)$ form a basis of $N$. We denote its dual basis by $B_\sigma = \{m_{i_1}^\sigma,\dots,m_{i_n}^\sigma\} = \{m^\sigma_\rho\}_{\rho\in\sigma(1)}$. The canonical line bundle on $\XSig$ is $\omega_{X_\Sigma/\kk}\cong\mcOX(-\sum_{\rho}D_\rho)$. We will write $\{\bso\}_{\sigma\in\Sigma(n)}$ for the Cartier data of $-\sum_{\rho}D_\rho$. We have $\bso = \sum_{m\in B_\sigma} m$, so $\bso$ is simply the sum over the dual basis defined by sigma.\\
	Every $m\in M$ defines a rational function $\chi^m$ on $\XSig$ with $\divi(\chi^m) = \sum_\rho \langle m,v_{\rho} \rangle D_{\rho}$ and two divisors $D = \sum_{\rho}a_\rho D_\rho, D' = \sum_{\rho}b_\rho D_\rho$ are rationally equivalent, if and only if $\sum_{\rho}a_\rho D_\rho = \sum_{\rho}(b_\rho + \langle m,v_\rho \rangle) D_\rho$ for some $m\in M$. In that case, we write $D\sim D'$. We denote the class group by $\Cl(X_\Sigma) = \Div(\XSig)/\Divo(\XSig) = \Div_T(\XSig)/\sim$.
	\begin{definition}
		\label{def:PD}
		Given a torus invariant divisor $D = \sum_\rho a_\rho D_\rho$, we define $P_D \subseteq M_\R$ to be the rational polyhedral polytope
		$$P_D=\{m\in M_\R \mid \langle m, v_{\rho}\rangle \geq - a_\rho \text{ for all } \rho\in\Sigma(1)\}.$$
	\end{definition}
	If $D' = D + \divi \chi^m$ for some $m\in M$, then $P_{D'} = P_D - m$. Moreover, $P_{D} + P_{D'}\subseteq P_{D+D'}$ and equality holds for $D,D'$ nef, which is equivalent to basepoint free in our case. 
	Recall that the global sections $\Gamma(X_\Sigma,\mcOX(D))$ of $\mcOX(D)$ are then given by $\bigoplus_{m\in P_D\cap M}\kk\cdot \chi^m.$\\
	Whenever we write polytope in this paper, we will always mean rational polyhedral polytope. A polytope is called a \emph{lattice polytope} if all its vertices lie in $M$.
\subsection{The quotient construction of toric varieties}
\label{subsection:QuotConstr}
	For a (smooth proper) toric variety $\XSig$, we define $S\coloneqq\kk[x_\rho \mid \rho\in\Sigma(1)]$ to be its \emph{Cox ring}. Further, set $Z(\Sigma)\coloneqq\mathcal{V}(\langle \prod_{\rho\not \in \sigma(1)}x_\rho\mid \sigma\in\Sigma(n)\rangle)$, $G_\Sigma\coloneqq\Spec \kk[\Cl(\XSig)]$, and we may identify $T$ with $\G_m^{\Sigma(1)}/G_\Sigma$ and interpret $G_\Sigma\subseteq \G_m^{\Sigma(1)}$. We have the following theorem.
	\begin{theorem}[{\cite[Theorem 2.1]{Elizondo_arithm_tor_var}}]
		\label{thm:CoxCoordsQuotientDescription}
		Let $\XSig$ be a smooth proper toric variety. Then $\XSig$ is the geometric quotient $(\A^{\Sigma(1)}\backslash Z(\Sigma))/G_\Sigma$, i.e.\@ for the quotient map $\pi$ and $x,y\in \A^{\Sigma(1)}\backslash Z(\Sigma)$ we have $\pi(x) = \pi(y)$ if and only if $x$ and $y$ lie in the same $G_\Sigma$-orbit.
	\end{theorem}
	Using this theorem, we write points in $\XSig$ as $[\bfx]$, the equivalence class of a tuple $\bfx = (x_{\rho})_{\rho\in\Sigma(1)}$ under the action of $G_\Sigma$ on $\A^{\Sigma(1)}\backslash Z(\Sigma)$.
\subsection{Local coordinates on toric varieties}
	\label{subsection:LocalCoordinates}
	The Cox ring $S$ has a grading by the class group $\Cl(X_\Sigma)$ with degree map $\prod_\rho x_\rho^{a_\rho}\mapsto [\sum_{\rho} a_\rho D_\rho]$. As usual for graded rings, we will write $S_\alpha$ for the set of homogeneous degree $\alpha$ elements of $S$. It is $\Gamma(\XSig,\mcOX(D)) \cong S_{[D]}$ via the $D$-homogenization $\chi^m\mapsto \prod_\rho x_\rho^{\langle m,v_\rho\rangle + a_\rho}$.\\
	There is a standard affine cover $\{U_\sigma\}_{\sigma\in \Sigma(n)}$ of $X_\Sigma$, where $U_\sigma = \Spec \kk[\sigma^\vee \cap M]$. In the quotient construction, $U_\sigma$ corresponds to the quotient of $\{\bfx\in\A^{\Sigma(1)}\backslash Z(\Sigma) \mid x_\rho\not=0\text{ for all }\rho\not\in\sigma(1)\}$; we will identify these two sets.\\
	For $\sigma(1) = \{\rho_{i_1},\dots,\rho_{i_n}\}$ the local charts $\varphi_\sigma\colon U_\sigma\to \A^n$ send $[\bfx]\in U_\sigma$ to $(\bfx^{m_{i_1}^\sigma},\dots,\bfx^{m_{i_n}^\sigma})\in\A^n$ where $B_\sigma = (m_{i_1}^\sigma,\dots,m_{i_n}^\sigma)$ is the basis of $\sigma^\vee\cap M$ dual to $v_{{i_1}},\dots,v_{{i_n}}$ and $\bfx^m \coloneqq \prod_{\rho\in\Sigma(1)}x_\rho^{\langle m, v_\rho \rangle}$.
\subsection{The Grothendieck-Witt ring}
	For more details, we refer to \cite{pauli2026pcmilecturenotesmotivic} and  \cite{JPPbezout}, the content presented here is taken from there.
	\begin{definition}
		The \emph{Grothendieck-Witt ring} $\GW(R)$ of a commutative ring with unity $R$ is the group completion of the semi-ring of isometry classes of non-degenerate symmetric bilinear forms over $R$ with respect to taking the direct sum where multiplication is given by the taking the tensor product.
	\end{definition}
	We will mainly be interested in the cases where $R$ is a field or the integers $\Z$. For enumerative problems over a general field $\kk$, the counting is done in its Grothendieck-Witt ring.
	\begin{proposition}
		If $\charact\kk\not=2$, $\GWk$ is generated by symbols $\langle a \rangle$ for $a\in\kk^\times$ subject to the relations 
		\begin{itemize}
			\item $\gw{a} = \gw{ab^2}$ for $a,b\in\kk^\times$,
			\item $\gw{a} + \gw{b} = \gw{a+b} + \gw{ab(a+b)}$ for $a,b,a+b\in\kk^\times$,
			\item $\gw{a}\gw{b} = \gw{ab}$ for $a,b\in\kk^\times$. 
		\end{itemize}
	\end{proposition}
	In particular, we can see $\kkq$ as a subset of $\GWk$ via $[a] \mapsto \gw{a}$.
	\begin{example}
		Here are some examples.
		\begin{itemize}
			\item For algebraically closed fields, $\kkq = \{1\}$ and $\GWk \cong \Z$.
			\item For the real numbers, $\kkq[\R] =  \{\pm1\}$ and $\GW(\R) \cong \gw{1}\Z \oplus \gw{-1}\Z$.
			\item For the integers, $\GW(\Z) = \gw{1}\Z \oplus \gw{-1}\Z$ \cite[Theorem II.4.3]{MilnorHusemoller1973}.
			\item For $\kk = \Q_p$ with $p\not = 2$, it is, by \cite[Lemma 12.1.7.]{Voight2021}, $\kkq = \{1,e,p,ep\}$ where $e$ is a quadratic non-residue modulo $p$. If $p\equiv3\mod{4}$, we can choose $e = -1$.
		\end{itemize}
	\end{example}
	\begin{definition}
		The class in $\GW(R)$ of the form on a free rank-$2$ module with Gram matrix $$\begin{pmatrix}
			0&1\\1&0
		\end{pmatrix}$$ is called the \emph{hyperbolic form}.
	\end{definition}
	\begin{remark}
		\label{rem:a-aish}
		It is $\gw{a} + \gw{-a} =  h$ for all $a\in R^\times$ .
	\end{remark}
	This allows the following definition.
	\begin{definition}
		The \emph{Witt ring} of $\kk$ of characteristic not two is $\W(\kk)\coloneqq\frac{\GWk}{\Z\cdot h}$.
	\end{definition}
	\begin{definition}
		We say that the \emph{rank} of a symmetric bilinear form $\beta\colon V\times V\to R$ is the rank of the $R$-module $V$.
	\end{definition}
	For $\charact \kk \not = 2$, an element in $\GWk$ is completely determined by its rank and its image in $\W(\kk)$.
	\begin{lemma}[\cite{MarkwigPayneShaw}]
		The map $\GW(\KK)\to\GWk; \gw{\sum_{k=k_0}c_kt^{k/n}}\mapsto \gw{c_{k_0}}$ is an isomorphism.
	\end{lemma}
	Via this map, we will identify $\GWk$ with $\GW(\KK)$ and $\kkq$ with $\kkq[\KK]$.\\
	Suppose $L/\kk$ is a finite separable field extension and let $b\colon V\times V\to L$ be a non-degenerate symmetric bilinear form over $L$. Composing with the field trace $\Tr_{L/\kk}$ and viewing $V$ as a $\kk$-vector space yields a non-degenerate symmetric bilinear form over $\kk$ 
	$$\Tr_{L/\kk}(b)\colon V\times V \overset{b}{\to} L \overset{\Tr_{L/\kk}}{\to} \kk.$$
	In fact, this defines a homomorphism of Grothendieck-Witt rings
	$$\Tr_{L/\kk}\colon \GW(L)\to\GWk.$$
\section{Combinatorial relative orientability}
	\label{section:CombRelOrient}
	We briefly explain the conjecture made in \cite{JPPbezout}. Recall from the introduction that a line bundle $\mcl$ is called \emph{relatively orientable} if $\mcl\otimes\omega_{X_\Sigma/\kk}$ is isomorphic to a square of a line bundle. In Jaramillo Puentes and Pauli's tropical proof of the enriched BKK theorem, the following plays the role of relative orientability on the tropical side.
	\begin{definition}[{\cite[Definition 5.12]{JPPbezout}}]
		\label{def:JPPcombrelorient}
		A tuple $(A_1,\dots,A_n)$ of $n$ finite subsets of $\Z^n$ is called \emph{combinatorially oriented} if the Minkowski sum of the convex hulls of $A_1,\dots,A_n$ has no points on the boundary with all coordinates odd.
	\end{definition}
	Further, they state the following conjecture.
	\begin{conjecture}[{\cite[Conjecture 5.14.]{JPPbezout}}]
		Let $\XSig$ be a smooth toric variety of dimension $n$ and let $\mcOX(E_1),\dots,\mcOX(E_n)$ be toric line bundles on $\XSig$. Further, let $A_i\subset \Z^n$ be the Newton polytopes of general sections of the vector bundles in chosen coordinates.\\
		Then $(A_1,\dots,A_n)$ is combinatorially oriented if and only if the vector bundle $\mcOX(E_1)\oplus \dots \oplus\mcOX(E_n)\to \XSig$ is relatively orientable. 
	\end{conjecture}
	Our goal is to answer this conjecture. While doing so, we will also introduce a new definition of combinatorial orientability that is, unlike this one using odd points, independent of bases and translations. In particular, it can be generalized to arbitrary lattices, does not require choosing coordinates and can nicely be translated to tropical geometry as a property of balanced polyhedral complexes.\\
	In the following sections, we are only interested in the relative orientability of vector bundles that are direct sums of line bundles. Hence, we are essentially only interested in the question whether a line bundle $\mcOX(D) \otimes \omega_{X_\Sigma/\kk}$ is isomorphic to a square, as $\det(\mcOX(E_1)\oplus\dots\oplus \mcOX(E_n)) \cong \mcOX(E_1+\dots+E_n)$.
	\begin{definition}
		We call a divisor $D$ \emph{relatively orientable} if $\mcOX(D)$ is relatively orientable.
	\end{definition}
	\begin{remark}
		\label{rem:relOrientOddCoeff}
		Note that the line bundle $\mcOX(D)$ of some $D\in \Div_T(\XSig)$ is isomorphic to a square if and only if $D\sim D'$ for some $D' = \sum_{\rho}a_\rho D_\rho\in \Div_T(\XSig)$ where all $a_\rho$ are even. Hence, a line bundle $\mcOX(D)$ is relatively orientable if and only if there is a divisor $D'=\sum_{\rho}a_\rho D_\rho$ with $D\sim D'$ and all $a_\rho$ odd, since one has to twist by the canonical bundle $\mcOX(-\sum_\rho D_\rho)$.
	\end{remark}
\subsection{Ample and nef divisors and their $P_D$}
	We recall translations of properties of divisors into the polytope $P_D\subseteq M_\R$ of Definition \ref{def:PD}. See \cite{cox2011toric} for more details and definitions. This theorem helps us to identify nef divisors.
	\begin{theorem}[{\cite[Theorem 6.1.7., Theorem 6.1.14., and Theorem 6.3.12.]{cox2011toric}}]
		\label{thm:BasepointFreeCharacterizations}
		Let $\XSig$ be a smooth proper toric variety. Then the following are equivalent for $D\in\Div_T(\XSig)$.
		\begin{enumerate}
			\item $D$ is basepoint free.
			\item $D$ is nef (numerically effective, i.e.\@ $D\cdot C\geq 0$ for every irreducible complete curve $C\subseteq \XSig$).
			\item $b_\sigma^D\in P_D$ for all $\sigma\in \Sigma(n)$.
			\item $P_D$ is the convex hull of $\{b_\sigma^D\mid \sigma\in\Sigma(n)\}$.
			\item $\{b_\sigma^D\mid \sigma\in\Sigma(n)\}$ is the set of vertices of $P_D$.
		\end{enumerate}
		In particular, $P_D$ is a lattice polytope if $D$ is nef.
	\end{theorem}
	\begin{example}
		For $\XSig$ a smooth proper toric variety, every ample $D\in\Div_T(\XSig)$ is nef.
	\end{example}
	\begin{example}
		For $\PP^n$ and $(\PP^1)^n$, every torus invariant divisor $D$ with nonempty $P_D$ is nef.
	\end{example}
	From \cite{cox2011toric}[Section 6.1,6.2] we get the following.
	\begin{proposition}
		Let $\XSig$ be a smooth proper toric variety. Then $\Sigma$ equals the normal fan of $P_D$ if and only if $D$ is ample. 
	\end{proposition}
	\begin{proof}
		All numbers in this proof reference to \cite{cox2011toric}. By Theorem 6.1.14, $D$ is ample if and only of $\varphi_D$ is strictly convex. The latter implies that $\Sigma$ is the normal fan of $P_D$, see equation (6.1.8) in the proof of the theorem. If $\Sigma$ is the normal fan of $P_D$, then $D = D_{P_D}$. The ampleness of $D$ follows now from Proposition 6.1.10.
	\end{proof}
\subsection{Relative orientability in $P_D$}
	\label{subsection:RelOrientInPD}
	In this subsection, we examine if it is possible to give a description of relative orientability of line bundles $\mcOX(D)$ in terms of their associated polytope $P_D$.\\
	The first thing to notice is that, as we are working over an arbitrary finite rank lattice $M$, \emph{odd points} can only be defined relative to a basis of $M$.	Recall from Subsection \ref{subsection:TorVars} that for $D \sim D'$, the polytopes $P_D$ and $P_{D'}$ are translates of each other. However, in that case $\mcOX(D) \cong \mcOX(D')$ and consequently, if one is relatively orientable then so is the other.\\
	Hence, to give a meaningful interpretation of relative orientability in terms of $P_D$, we need to find a property that is independent of translation and doesn't require a choice of basis. As a first step, we introduce the following notation.
	\begin{definition}
		For $D = \sum_\rho a_\rho D_\rho\in \Div_T(\XSig)$ we write 
		$$L_D \coloneqq \{m\in M_\R\mid \langle m,v_\rho \rangle = -a_\rho \text{ for some }\rho \in \Sigma(1)\} = \bigcup_{{\rho\in\Sigma(1)}} H^D_\rho$$
		for the union of the supporting hyperplanes $H^D_\rho \coloneqq \{m\in M_\R\mid \langle m,v_\rho \rangle = -a_\rho\}$ in the definition of $P_D$.
	\end{definition}	
	Clearly, just as for $P_D$, we have $L_{D'} = L_D - m$ if $D' = D+ \divi \chi^m$ for some $m\in M$. Notice however, that $L_{D+D'} \not =  L_{D} + L_{D'} = M_\R$. Furthermore, $\partial P_D\subset L_D$.
	\begin{example}
		\label{ex:defFirstHirz}
		In Figure \ref{fig:HirzebruchSections}, we consider divisors on the \emph{first Hirzebruch surface}, that is, the toric variety defined by the smooth complete fan with rays generated by $v_1 = (0,1), v_2 = (1,0), v_3 = (0,-1), v_4 = (-1,-1)$.
	\end{example}
	The following proposition allows us to translate relative orientability into a more combinatorial property. Let $\pi_2\colon M \to M/2M$ denote the reduction modulo $2$.
	\begin{proposition}
		\label{prop:theEquivalence}
		Let $\XSig$ be a smooth proper toric variety and $D = \sum_\rho a_\rho D_\rho\in \Div_T(\XSig)$. Then $D$ is relatively orientable, if and only if $\pi_2(L_D\cap M) \not = M/2M$.
	\end{proposition}
	As we will need them later, we source parts of the proof out into two lemmas.
	\begin{lemma}
		\label{lem:whenOddSolutions}
		Let $v\in N$ be primitive and $u\in M$. Then $$\{\langle m, v\rangle \mid m\in u+2M\} = \langle u,v\rangle + 2\Z.$$
		In particular, for $z\in\Z$, the hyperplane $\{m\in M\mid \langle m,v\rangle = z\}$ meets the coset $u+2M$ if and only if $z\equiv \langle u,v\rangle \mod 2$.
	\end{lemma}
	\begin{proof}
		Every $m\in u+2M$ is of the form $m = u+2k$ with $k\in M$, so $\langle m,v\rangle = \langle u,v\rangle + 2\langle k,v\rangle$. As $v$ is primitive, it can be extended to a basis of $N$; the corresponding dual basis element $w\in M$ satisfies $\langle w,v\rangle = 1$. Hence $\langle\cdot,v\rangle\colon M\to\Z$ is surjective and $\langle k,v\rangle$ runs through all of $\Z$.
	\end{proof}
	\begin{lemma}
		\label{lem:notInPiMakesOddCoeff}
		Let $\XSig$ be a smooth proper toric variety and $D = \sum_\rho a_\rho D_\rho\in \Div_T(\XSig)$ and $u\in M$. Then $u+2M\not\in\pi_2(L_D\cap M)$ if and only if all coefficients of $D+\divi(\chi^u)$ are odd.
	\end{lemma}
	\begin{proof}
		By Lemma \ref{lem:whenOddSolutions}, the hyperplane $H^D_\rho$ does not meet $u+2M$ if and only if $-a_\rho\not\equiv\langle u,v_\rho\rangle\mod 2$, i.e.\@ if and only if $a_\rho+\langle u,v_\rho\rangle$ is odd. Since $D+\divi(\chi^u) = \sum_\rho \bigl(a_\rho+\langle u,v_\rho\rangle\bigr)D_\rho,$ the claim follows by considering all $\rho\in\Sigma(1)$ simultaneously.
	\end{proof}
	\begin{proof}[Proof of the proposition]
		It is $\pi_2(L_D\cap M) \not = M/2M$ if and only if there exists a $u\in M$ such that $u+2M \not\in\pi_2(L_D\cap M)$. This condition is, by Lemma \ref{lem:notInPiMakesOddCoeff}, equivalent to $D+\divi \chi^u \sim D$ having only odd coefficients, which is, finally, equivalent to $D$ being relatively orientable by Remark \ref{rem:relOrientOddCoeff}.
	\end{proof}
	In fact, given that the divisor is relatively orientable, we know which element of $M/2M$ is not hit and how to find a $D'\sim D$ with only odd coefficients.
	\begin{proposition}
		\label{prop:theMissingPoint}
		Let $\XSig$ be a smooth proper toric variety and $D\in\Div_T(\XSig)$ with Cartier data $\{b_\sigma^D\}_{\sigma\in\Sigma(n)}$ and let $u_\sigma = b_\sigma^D + \bso$.\\
		Then $D$ is relatively orientable, if and only if $u_\sigma + 2M$ does not depend on $\sigma$.\\
		Further, if $D$ is relatively orientable, $$\pi_2(L_D\cap M) = (M/2M)\backslash \{u_\sigma + 2M\}$$
		and $D + \divi \chi^{u_\sigma}$ has only odd coefficients.
	\end{proposition}
	\begin{proof}
		Let $D$ be relatively orientable. For $\rho\in\sigma(1)$, we have $a_\rho + \langle u_\sigma,v_\rho\rangle = a_\rho - a_\rho + 1 = 1$ by definition of the Cartier data, so $u_\sigma+2M\not\in\pi_2(L_D\cap M)$ by Lemma \ref{lem:notInPiMakesOddCoeff}. Now let $w+2M\not\in\pi_2(L_D\cap M)$. Then $\langle u_\sigma - w,v_\rho\rangle$ is even for all $\rho\in\sigma(1)$, by, again, Lemma \ref{lem:notInPiMakesOddCoeff}. As $\{v_\rho\}_{\rho\in\sigma(1)}$ is a basis of $N$, this implies $u_\sigma - w\in 2M$. This shows the further-part and the first direction.\\		
		For the other direction, let $u = u_\sigma + 2k\in u_\sigma + 2M$. The divisor $D + \divi \chi^u$ has coefficient $a_\rho + \langle (b_\sigma^D + \bso) + 2k,v_\rho\rangle = a_\rho - a_\rho + 1 + 2 \langle k,v_\rho\rangle \equiv 1 \mod{2}$ for $\rho\in\sigma(1)$. As each ray ray lies in some maximal cone, the relative orientability follows from Remark \ref{rem:relOrientOddCoeff}.
	\end{proof}
	\begin{corollary}
		Let $\XSig$ be a smooth proper toric variety and $D = \sum_\rho a_\rho D_\rho\in \Div_T(\XSig)$. Then $D$ is relatively orientable, if and only if $|\pi_2(L_D\cap M)| = 2^n -1$.
	\end{corollary}
	As the boundary of $P_D$ is a subset of $L_D$, we get the following corollary, which gives a nice necessary condition for relative orientability that only depends on the actual polytope $P_D$.
	\begin{corollary}
		\label{cor:RelOrientImpliesCombRelOrient}
		Let $\XSig$ be a smooth proper toric variety and $D\in\Div_T(\XSig)$. If $D$ is relatively orientable, then $\pi(\partial P_D\cap M) \not= M/2M$.
	\end{corollary}
	Here are two examples that show that this condition is not sufficient.
	\begin{example}
		\label{ex:FirstHirzNotRelOrientRedundantLine}
		Consider the first Hirzebruch surface and the divisor $D = D_3+D_4$. It is
		$$P_D\cap \Z^2 = \partial P_D\cap \Z^2 = \{(0,0), (1,0),(0,1)\},$$
		and obviously, $\pi_2(\partial P_D\cap \Z^2)\not=\Z^2 / 2\Z^2$. However, the equation $\langle m, v_3\rangle = -m_2 = -1$ has an odd solution. Hence, $\pi_2(L_D\cap \Z^2) = \Z^2 / 2\Z^2$, so $D$ is not relatively orientable by Proposition \ref{prop:theEquivalence}.
		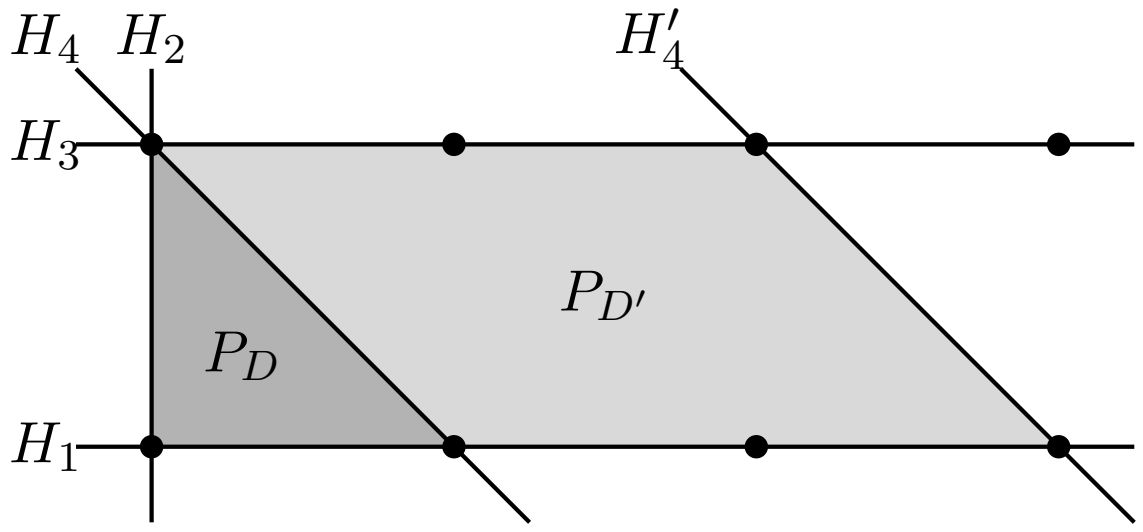
\begin{figure}
			\centering
			\begin{tikzpicture}
				\begin{scope}[scale=2, every node/.style={scale=1}]
					\fill[gray!30] (0,0) -- (0,1) -- (2,1) -- (3,0) -- cycle;
					\fill[gray!60] (0,0) -- (0,1) -- (1,0) -- cycle;
					
					\node at (0.3,0.3) {$P_{D}$};
					\node at (1.5,0.5) {$P_{D'}$};
					
					\node at (-0.35,1) {$H_3$};
					\node at (-0.35,1.35) {$H_4$};
					\node at (0,1.35) {$H_2$};
					\node at (-0.35,0) {$H_1$};
					\node at (1.65,1.35) {$H_4'$};
					
					\filldraw[color = black, fill=black] (0,0) circle (1pt);
					\filldraw[color = black, fill=black] (1,0) circle (1pt);
					\filldraw[color = black, fill=black] (0,1) circle (1pt);
					\filldraw[color = black, fill=black] (1,1) circle (1pt);
					\filldraw[color = black, fill=black] (2,1) circle (1pt);
					\filldraw[color = black, fill=black] (2,0) circle (1pt);
					\filldraw[color = black, fill=black] (3,1) circle (1pt);
					\filldraw[color = black, fill=black] (3,0) circle (1pt);
					
					\draw[thick,black] (-0.25,1) -- (3.25,1);					
					\draw[thick,black] (0,-0.25) -- (0,1.25);
					\draw[thick,black] (-0.25,1.25) -- (1.25,-0.25);
					\draw[thick,black] (-0.25,0) -- (3.25,0);
					\draw[thick, black] (1.75,1.25) -- (3.25,-0.25);
				\end{scope}
			\end{tikzpicture}
			\caption{A drawing of $P_{D}$ and $P_{D'}$ with $D = D_3 + D_4$ and $D' = D_3+3D_4$ for the first Hirzebruch surface. Line $l_i$ comes from the equation of ray $i$. It is $L_D = \bigcup_{i=1}^4H^D_i$ and $L_{D'} = H_4^{D'}\cup\bigcup_{i=1}^3H_i$. Further, $L_{P_D} = H^D_1\cup H^D_2\cup H^D_4$ and $L_{P_{D'}} = L_{D'}$ (See Definition \ref{def:LP}).}
			\label{fig:HirzebruchSections}
		\end{figure}
	\end{example}
	The reason in this example for the converse of Corollary \ref{cor:RelOrientImpliesCombRelOrient} to fail, is that one of the inequalities describing $P_{D}$ is redundant, see Figure \ref{fig:HirzebruchSections}. But even if all inequalities are non-redundant, we can run into a similar problem.
	\begin{example}
		\label{ex:NoOddNotRelOrient}
		Consider the projective space $\PP^3$ and the nef divisor $D = D_0$. We get 
		\begin{align*}
			P_D \cap \Z^3= \partial P_D\cap \Z^3 &= \{m\in\Z^3\mid m_1 \geq 0,\ m_2\geq 0,\ m_3\geq 0,\ -m_1-m_2-m_3 \geq -1\} \\&= \{(0,0,0),(1,0,0),(0,1,0),(0,0,1)\}.
		\end{align*}
		As $\partial P_D \cap \Z^3$ contains only four points, $\pi_2(\partial P_D \cap \Z^3)\not = \Z^3 / 2\Z^3$. One can easily check that $\pi_2(L_D\cap M) = \Z^3 / 2\Z^3$, so $D$ is not relatively orientable.
	\end{example}
	By Proposition \ref{prop:theEquivalence}, the converse of Corollary \ref{cor:RelOrientImpliesCombRelOrient} is true if $\pi_{2}(L_D\cap M) = \pi_{2}(\partial P_D\cap M)$. This is e.g.\@ the case for $X_\Sigma = \PP^2$ and $X_\Sigma = (\PP^1)^n$ for divisors with full-dimensional $P_D$. In Lemma \ref{lem:NecCondForLP=DelP} we will see a sufficient condition for this to hold.
	\begin{remark}
		Call $m\in M$ \emph{odd} with respect to a $\Z$-basis $B$ of $M$, if all its coordinates in $B$ are odd. This is equivalent to $\pi_2(m) = \sum_{b\in B}b + 2M$.\\ Then the set $L_D + u$ contains no odd points w.r.t.\@ $B$ if and only if $u + \sum_{b\in B} b + 2M \not\in \pi_2(L_D\cap M)$.
	\end{remark}
	Summing up all statements of this subsection, we get the following theorem.
	\begin{theorem}
		\label{thm:DrelOrientInLattice}
		Let $\XSig$ be a smooth proper toric variety and $D\in\Div_T(\XSig)$ with Cartier data $\{b_\sigma^D\}_{\sigma\in\Sigma(n)}$. Then the following are equivalent.
		\begin{itemize}
			\item $\mcOX(D)$ is relatively orientable, i.e. $\mcOX(D)\otimes \omega_{X_\Sigma/\kk}$ isomorphic to a square.
			\item $D$ is relatively orientable.
			\item $D\sim D'=\sum_{\rho}a_\rho D_\rho$ with all $a_\rho$  odd.
			\item It is $\pi_2(L_D\cap M)\not=M/2M$.
			\item It is $\pi_2(L_D\cap M)=(M/2M)\backslash \{(b_\sigma^D + \bso) + 2M\}$ for some/all $\sigma\in\Sigma(n)$.
			\item It is $|\pi_2(L_D\cap M)| = 2^n-1$.
			\item $(b_\sigma^D + \bso) + 2M$ does not depend on $\sigma$.
		\end{itemize}
		Further, they imply $\pi_2(\partial P_D\cap M) \not = M/2M$.
	\end{theorem}
\subsection{Combinatorial orientability}
	\label{subsection:combinatorialOrientability}
	The next definition is basically the $L_D$-construction from above, except that we start with only a polytope instead of a divisor.
	\begin{definition}
		\label{def:LP}
		For a full-dimensional polytope $P\subset M_\R$, we define $L_P$ to be the union of the supporting hyperplanes of the minimal description of $P$ as an intersection of half-spaces.
	\end{definition}
	As this minimal description of $P$ is unique for a full-dimensional polytope (see e.g.\@ \cite{grunbaum2003convex}), this is well-defined. 
	\begin{remark}
		Clearly, $\partial P \subseteq L_P$. See Figure \ref{fig:HirzebruchSections} and Example \ref{ex:defFirstHirz} for an example. Further, $L_{P_D} \subseteq L_D$ (see again Figure \ref{fig:HirzebruchSections}) and equality holds for $\Sigma$ the normal fan of $P_D$. Recall that this is for example the case if $D$ is ample \cite{cox2011toric}[Equation (6.1.8)].
	\end{remark}
	We can now give our new definition of combinatorial orientability, inspired by Proposition \ref{prop:theEquivalence}.
	\begin{definition}
		\label{def:combOrientP}
		We say that a polytope $P \subset M_\R$ is \emph{combinatorially orientable} if $P$ is full-dimensional and $\pi_2(L_P\cap M)\not = M/2M$.
	\end{definition}
	\begin{remark}
		Adapting the proof of Proposition \ref{prop:theMissingPoint}, one can see that for a full-dimensional $P$ it holds that $\pi_2(L_P\cap M)\not = M/2M$ is equivalent to $|\pi_2(L_P\cap M)|=2^n-1$ if for one vertex the normal vectors of the adjacent facets form a basis modulo $2$.
	\end{remark}
	Sometimes, considering the boundary of $P$ suffices.
	\begin{lemma}
	\label{lem:NecCondForLP=DelP}
		Let $\dim M_\R = n\geq 2$ and be $P\subseteq M_\R$ a smooth full-dimensional lattice polytope with all edges having lattice length $\geq n-1$. \\Then $\pi_2(L_P\cap M) = \pi_2(\partial P\cap M)$.
		In particular, this is always the case for $n=2$.
	\end{lemma}
	\begin{proof}
		Consider a facet $F$ of $P$ containing a vertex $v$. Let $H_F$ be the affine hypersurface containing $F$. As $P$ is smooth and full-dimensional, the directions $w_1,\dots,w_n$ of the edges adjacent to $v$ form a $\Z$-basis of $M$. Let w.l.o.g.\@ $w_1,\dots,w_{n-1}$ be the directions of the edges contained in $F$. By convexity of $F$ and the lengths of the edges being greater or equal to $n-1$, it is
		$$\conv\{v,v+(n-1)w_1,\dots,v+(n-1)w_{n-1}\}\subseteq F.$$
		Hence, for $\varepsilon\in\{0,1\}^{n-1}$,
		$$M\ni v+ \sum_{i} \varepsilon_i w_i = v+ \sum_{i} \frac{\varepsilon_i}{n-1} (n-1)w_i \in F.$$
		It follows that $\pi_2(F\cap M) = \pi_2(H_F\cap M)$. Taking the union over all facets proves the claim.
	\end{proof}
	\begin{example}
		This inequality is sharp: For the divisor $(n-2)D_0$ on $\PP^n$ where $n\geq 3$ and $P = P_{(n-2)D_0}$ one checks easily that $\pi_2(L_P\cap M) = \Z^n/2\Z^n$ while $(1,\dots,1)\not\in\pi_2(\partial P\cap M)$. 
	\end{example}
	From  Proposition \ref{prop:theEquivalence} we immediately get the following theorem.
	\begin{theorem}
		\label{thm:RelOrientImpliesCombOrient}
		Let $\XSig$ be a smooth proper toric variety and $D\in\Div_T(\XSig)$.\\
		If $D$ is relatively orientable with full-dimensional $P_D$, then the polytope $P_D$ is combinatorially orientable.\\
		Further, if $D$ is ample, then $D$ is relatively orientable if and only if the polytope $P_D$ is combinatorially orientable.
	\end{theorem}
	\begin{example}
		\label{ex:blowUpP3}
		Consider the blow up of $\PP^3$ at $[1:0:0:0]$. This is, the toric variety defined by the same complete smooth fan as $\PP^n$, except that one adds a ray $\rho_4$ in the $e_1+e_2+e_3$ direction and subdivides the maximal cone $\sigma_0$ into cones $\sigma'_i$ with rays $\{\rho_1,\dots,\rho_4\}\backslash\{\rho_i\}$ for $i=1,2,3$.\\
		One computes that $P_{aD_0-bD_4} = \{m\in\R^3_{\geq0}\mid b\leq m_1+m_2+m_3\leq a\}$. Consider the divisors $E_1 = E_2 = D_0$ and $E_3 = 2D_0-D_4$. Then $E_1$ and $E_2$ are nef and have full dimensional polytopes corresponding to them, $E_3$ even is ample. Their sum, $E \coloneqq 4D_0 - D_4$ is ample, too.\\
		However, $P_E$ has edges of length one and in fact $\pi_2(\partial  P_E\cap \Z^3)\not= \pi_2(L_D\cap \Z^3) =  \pi_2(L_{P_D}\cap \Z^3) = \Z^3/2\Z^3$, so $D$ is not relatively orientable.
	\end{example}
	If we start with a polytope and are allowed to choose a fan, the following holds.
	\begin{theorem}
		\label{thm:combOrientGivesOrientDivisor}
		Let $P$ be a combinatorially orientable lattice polytope. Then, for every smooth refinement $\Sigma$ of the normal fan of $P$, there exists a relatively orientable $D\in\Div_T(\XSig)$ such that $P = P_D$.
	\end{theorem}
	\begin{proof}
		Let $w\in M$ such that $\pi_2(w)\not\in (M/2M)\backslash\pi_2(L_P\cap M)$. Let $\Sigma$ be a smooth refinement of the normal fan of $P$. Number its rays such that the first $r$ rays are the ones from the normal fan. For $i=r+1,\dots,s$ choose $a_{\rho_i}\in\Z$ such that $\pi_{2}(w)\not\in\pi_{2}(\{m\in M\mid \langle m ,v_{\rho_i}\rangle = -a_{\rho_i}\})$ and $P\subseteq\{m\in M\mid \langle m ,v_{\rho_i}\rangle \geq -a_{\rho_i}\}$. By Lemma \ref{lem:whenOddSolutions}, this is the case for every $a_\rho\in (1+\langle w, v_{\rho_i}\rangle) + 2\Z$. Then $P = P_D$ for $D\coloneqq \sum_{\rho}a_\rho D_\rho$ and $D$ is relatively orientable by Proposition \ref{prop:theEquivalence}.
	\end{proof}
	\begin{remark}
		\label{rem:findTheSigma}
		The usual approach for the classical BKK theorem that makes a statement over solutions in $(\C^\times)^n$, see e.g.\@ \cite{fulton1993introductionTorVars}, is this other way around: Starting with a lattice polytope $P$ (or several polytopes) one constructs a fan $\Sigma$ such that there is a divisor $D$ with $P_D = P$. Unsurprisingly, the question whether this divisor is relatively orientable depends on the explicit toric variety. The polytope in Example \ref{ex:FirstHirzNotRelOrientRedundantLine}, seen in Figure \ref{fig:HirzebruchSections}, can either be interpreted as $P_D$ for the first Hirzebruch surface with $D = D_3+D_4$ or for $\PP^2$ with $D = D_3$. In the first case, $D$ is not relatively orientable, in the second one it is.\\
		Further, the choice of $D$ is not unique either: It is $P_D = P_{D'}$ for the first Hirzebruch surface with $D = D_3+D_4$ and $D' = 2D_3 + D_4$. In this case $D$ is not relatively orientable while $D'$ is.
	\end{remark}
\subsection{Relative orientability in Newton polytopes}
	\label{subsection:RelOrientInNP}
	Recall from Section \ref{section:TorVarsAndFans} that the global sections of $\mcOX(D)$ are given by $\bigoplus_{m\in P_D\cap M} \kk\cdot \chi^m$ for $D\in \Div_T(\XSig)$. This motivates the following definition.
	\begin{definition}
		\label{def:newtonpolytope}
		Let $f = \sum_{m\in P_D\cap M} a_m \chi^m$ be a global section of $\mcOX(D)$. We define the \emph{Newton polytope} of $f$, $\NP(f)\subseteq M_\R$, to be the convex hull of $\{m\in P_D\cap M\mid a_m\not = 0\}$ in $M_\R$.\\
		Further, we define the \emph{Newton polytope of $D$}, $\NP(D)$ to be the convex hull of $P_D\cap M$  in $M_\R$.
	\end{definition}
	Clearly, $\NP(f) \subseteq \NP(D) \subseteq P_D$. If $f$ is a general section, then $\NP(f) = \NP(D)$. It is $\NP(D) = P_D$ if and only if $P_D$ is a lattice polytope. We get the following by applying Theorem \ref{thm:RelOrientImpliesCombOrient}.
	\begin{corollary}
		If $D\in\Div_T(\XSig)$ is relatively orientable and $P_D$ a full-dimensional lattice polytope, then $\NP(D)$ is combinatorially orientable.
	\end{corollary}
	The following example shows that we should restrict ourselves mostly to divisors with $P_D$ a lattice polytope (e.g. nef divisors) if we use $\NP(D)$. This is not a strong restriction.
	\begin{example}
		\label{ex:NPNotCombOrientButDOrient}
		Let $\XSig$ be the third Hirzebruch surface, that is, the toric variety defined by the smooth complete fan with rays generated by $v_1 = (0,1), v_2 = (1,0), v_3 = (0,-1), v_4 = (-1,-3)$, and $D = \sum_{i=1}^{4}D_i$. In that case, $D$ is relatively orientable, $P_D$ full-dimensional but not a lattice polytope. In this example, $\NP(D)$ is not combinatorially orientable.
	\end{example}
	\begin{example}
		\label{ex:notLatticeNotLatticeLattice}
		Consider the third Hirzebruch surface and the divisor $D = \sum_{i=1}^{4}3D_i$. Then $P_D$ is the triangle with corners $(-3,-3),(12,-3),(-3,2)$, a lattice polygon. However, as $b_{\sigma_2'} = (-3,3)\not\in P_D$, $D$ is not nef.\\
		Further, let $D' = \sum_{i=1}^{4}2D_i$ and $D'' = \sum_{i=1}^{4}D_i$, so $D = D' + D''$. Then neither $P_{D'}$ nor $P_{D''}$ are lattice polytopes. In particular, $\NP(D') + \NP(D'') \subsetneqq \NP(D) = P_D = P_{D'} + P_{D''}$.
	\end{example}
	We get the following version of Jaramillo Puentes and Pauli's conjecture.
	\begin{theorem}
		\label{thm:JPPforNef}
		Let $\XSig$ be a smooth proper toric variety of dimension $n$ and let $\mcOX(E_1)$, $\ldots$, $\mcOX(E_n)$ be toric line bundles on $\XSig$ and all $E_i$ be nef. Further, let $A_i \coloneqq P_{E_i}\cap M\subset M$.\\
		If the vector bundle $\mcOX(E_1)\oplus \dots \oplus\mcOX(E_n)\to \XSig$ is relatively orientable, and $P_{E_1}+\dots+P_{E_n}$ of full dimension, then $\conv(A_1+\dots+A_)n$ is combinatorially orientable.
	\end{theorem}
	\begin{proof}
		The vector bundle $\mcOX(E_1)\oplus \dots \oplus\mcOX(E_n)$ being relatively orientable is equivalent to $\mcOX(\sum_{i=1}^n E_i)$ being relatively orientable. As $P_{\sum_{i=1}^n E_i} = \sum_{i=1}^n P_{E_i}$, we get the statement by applying Theorem \ref{thm:RelOrientImpliesCombOrient}, as the convex hull of $A_i$ is again $P_{E_i}$, since the $P_{E_i}$ are lattice polytopes.
	\end{proof}
	In the theorem above, it is crucial that $P_{\sum_{i=1}^n E_i} = \sum_{i=1}^n P_{E_i}$, as the next example demonstrates.
	\begin{example}
		Consider the first Hirzebruch surface. Let $E_1 = 2D_3+D_4$ and $E_2 = 2D_3+4D_4$. Then $P_{E_1}+P_{E_2} = \conv\{(0,0),(0,3),(2,3),(5,0)\}$ is not combinatorially orientable. However, $$E_1+E_2 = 4D_3+5D_4 \sim 4D_3+5D_4+\divi(\chi^{(1,1)}) = D_1+D_2+3D_3+3D_4$$
		has only odd coefficients, so $E_1+E_2$ is relatively orientable by Remark \ref{rem:relOrientOddCoeff}. Hence, by Theorem \ref{thm:RelOrientImpliesCombOrient}, the polytope $P_{E_1+E_2}\supsetneqq P_{E_1}+P_{E_2}$ is combinatorially orientable.\\
		Let now $E_1 = 2D_3+D_4$ and $E_2 = D_3+4D_4$. Then $P_{E_1}+P_{E_2} = \conv\{(0,0),(0,2),(3,2),(5,0)\}$ is combinatorially orientable, but $P_{E_1+E_2} = \conv\{(0,0),(0,3),(2,3),(5,0)\}\supsetneqq P_{E_1}+P_{E_2}$ is not. In particular, $E_1+E_2$ is not relatively orientable by Theorem \ref{thm:RelOrientImpliesCombOrient}.\\
		In both cases, the inequality $-m_2\geq -2$ is redundant for $P_{E_1}$, but contributes to $P_{E_1+E_2}$.
	\end{example}
\subsection{Odd boundary points}
	In this subsection, we want to relate our new definition of combinatorial orientability to the original definition given in \cite{JPPbezout}. Let $D = \sum_{\rho}a_\rho D_\rho$ be a divisor. As mentioned earlier, a local trivialization of $\mcOX(D)$ on $U_\sigma$ is given by the section $\chi^{b_\sigma^D}$. In other words, we have an isomorphism $\mcOX(D)|_{U_\sigma} \xrightarrow{\;\cong\;} \mcOX|_{U_\sigma}$ given by $s\mapsto s / \chi^{b_\sigma^D}$. Under this local trivialization, we get that the sections of $\mcOX|_{U_\sigma}$ are given by $\bigoplus_{m\in P_D\cap M} \kk\cdot\chi^{m-b_\sigma^D}$.\\
	Writing these sections in the local coordinates $(y_1,\dots,y_n) = (\bfx^{m_{i_1}^\sigma},\dots,\bfx^{m_{i_n}^\sigma})$ on $U_\sigma$, we end up with
	$$\bigoplus_{m\in P_D\cap M} \kk\cdot y_1^{\langle m-b_\sigma^D, v_{i_1}\rangle}\cdots y_n^{\langle m-b_\sigma^D, v_{i_n}\rangle}.$$
	In particular, we can define a Newton polytope that lies in $\R^n$ for each chart. While the coordinates we define in Section \ref{section:ExplicitRelOrientOnTorVar} differ from these by a sign in the first coordinate, this does not affect the Newton polytope, and hence we are going to omit it.
	\begin{definition}
		Given a divisor $D\in \Div_T(\XSig)$ and $\sigma\in\Sigma(n)$ we write $\psisd\colon M_\R \to \R^n$ for the affine map that first shifts by $b_\sigma^D$ and then gives the coordinates with respect to the basis $\{m_i^\sigma\}_{i\in\sigma(1)}$.
	\end{definition}
	This map sends $M$ to $\Z^n$, as $b_\sigma^D\in M$ and $\{m_i^\sigma\}_{i\in\sigma(1)}$ is a $\Z$-basis. In formulas, $\psisd(m) = (\langle m-b_\sigma^D, v_{i_1}\rangle,\dots,\langle m-b_\sigma^D, v_{i_n}\rangle)$ where $\sigma(1) = \{\rho_{i_1},\dots,\rho_{i_n}\}$ with $i_1<\dots<i_n$. As this map is affine, it commutes with taking the convex hull of a set.
	\begin{definition}
		Given a global section $f = \sum_{m\in P_D\cap M} a_m \chi^m$ of $\mcOX(D)$, we define its \emph{Newton polytope with respect to $\sigma\in\Sigma(n)$}, $\NP_{\sigma}(f) \subseteq\R^n$ to be
		$$\psisd(\NP(f)) = \conv\{(\langle m-b_\sigma^D, v_{i_1}\rangle,\dots,\langle m-b_\sigma^D, v_{i_n}\rangle)\mid m\in P_D\cap M,\ a_m\not=0\},$$ where $\sigma(1) = \{\rho_{i_1},\dots,\rho_{i_n}\}$ with $i_1<\dots<i_n$.\\
		Further, we define the \emph{Newton polytope of $D$ with respect to $\sigma$}, $\NP_\sigma(D)$, to be $\psisd(\NP(D))$ and $L_{D,\sigma} \coloneqq \psisd(L_D)$.
	\end{definition}
	\begin{example}
		Let $\Sigma$ be the following (smooth proper) fan.
		
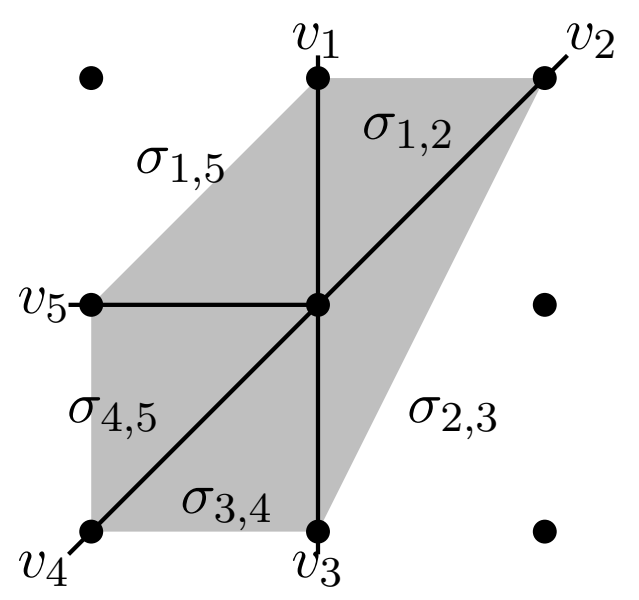
\begin{center}
			\begin{tikzpicture}[scale=1.5, every node/.style={scale=1}]
				\fill[lightgray] (0,1) -- (1,1) -- (0,-1) -- (-1,-1) -- (-1,0) -- cycle;
				
				\draw[thick,black] (0.0,0.0) -- (0,1.1);
				\draw[thick,black] (0.0,0.0) -- (1.1,1.1);
				\draw[thick,black] (0.0,0.0) -- (0,-1.1);
				\draw[thick,black] (0.0,0.0) -- (-1.1,-1.1);
				\draw[thick,black] (0.0,0.0) -- (-1.1,0);
				
				\node[anchor=south] at (0,1) {$v_1$};
				\node[anchor=south west] at (1,1) {$v_2$};
				\node[anchor=north] at (0,-1) {$v_3$};
				\node[anchor=north east] at (-1,-1) {$v_4$};
				\node[anchor=east] at (-1,0) {$v_5$};
				\node at (0.4,0.75) {$\sigma_{1,2}$};
				\node at (0.6,-0.5) {$\sigma_{2,3}$};
				\node at (-0.4,-0.9) {$\sigma_{3,4}$};
				\node at (-0.9,-0.5) {$\sigma_{4,5}$};
				\node at (-0.6,0.6) {$\sigma_{1,5}$};
				
				\foreach \x in {-1,...,1}{
					\foreach \y in {-1,...,1}{
						\fill (\x,\y) circle (1.5pt);
					}
				}				
			\end{tikzpicture}
		\end{center}
		This fan has primitive ray generators $$v_1 = (0,1), v_2 = (1,1), v_3=(0,-1), v_4 = (-1,-1), v_5 = (-1,0).$$
		Consider the divisor $D = -2D_2+2D_3+4D_4+3D_5$; $P_D$ can be seen in Figure \ref{fig:NPsigmas}. $D$ is relatively orientable, as $(0,1)\not\in\pi_2(L_D\cap M)$.
		The five other pictures in Figure \ref{fig:NPsigmas}, are its Newton polytope with respect to the different two dimensional cones of $\Sigma$.
		\begin{figure}
			\begin{center}
				\begin{tikzpicture}[scale=1.2, every node/.style={scale=1}]
					\pgfmathsetmacro{\ltl}{0.25}
					
					\begin{scope}[shift={(-3,0)}]
						\node at (-0.5,1) {$P_D$:};
						\fill[lightgray] (0,2) -- (2,2) -- (3,1) -- (3,0) -- (2,0) -- cycle;
						
						\draw[thick, black] (-0.25,2.25) -- (2.25,-0.25);
						\draw[thick, black] (-0.25,0) -- (3.25,0);
						\draw[thick, black] (3,-0.25) -- (3,2.25);
						\draw[thick, black] (-0.25,2) -- (3.25,2);
						\draw[thick, black] (1.75,2.25) -- (3.25,0.75);
						\foreach \x in {0,...,3}{
							\foreach \y in {0,...,2}{
								\fill (\x,\y) circle (1.5pt);
							}
						}				
						\fill (0,0) circle (2.5pt);
					\end{scope}
					\begin{scope}[shift={(2.5,0)}]
						\fill[lightgray] (0,0) -- (0,1) -- (1,2) -- (2,2) -- (2,0) -- cycle;
						\draw[thick, black] ({0-\ltl},0) -- ({2+\ltl},0);
						\draw[thick, black] ({0-\ltl},2) -- ({2+\ltl},2);
						\draw[thick, black] (0,{0-\ltl}) -- (0,{2+\ltl});
						\draw[thick, black] (2,{0-\ltl}) -- (2,{2+\ltl});
						\draw[thick, black] ({0-\ltl},{1-\ltl}) -- ({1+\ltl},{2+\ltl});
						
						\node at (-0.6,1) {$\NP_{\sigma_{1,2}}$:};					
						\foreach \x in {0,...,2}{
							\foreach \y in {0,...,2}{
								\fill (\x,\y) circle (1.5pt);
							}
						}				
						\fill (0,0) circle (2.5pt);
					\end{scope}
					\begin{scope}[shift={(7,0)}]
						\node at (-0.6,1) {$\NP_{\sigma_{2,3}}$:};					
						\fill[lightgray] (0,0) -- (0,2) -- (1,2) -- (2,1) -- (2,0) -- cycle;
						\draw[thick, black] ({0-\ltl},0) -- ({2+\ltl},0);
						\draw[thick, black] ({0-\ltl},2) -- ({2+\ltl},2);
						\draw[thick, black] (0,{0-\ltl}) -- (0,{2+\ltl});
						\draw[thick, black] (2,{0-\ltl}) -- (2,{2+\ltl});
						\draw[thick, black] ({1-\ltl},{2+\ltl}) -- ({2+\ltl},{1-\ltl});
						
						\foreach \x in {0,...,2}{
							\foreach \y in {0,...,2}{
								\fill (\x,\y) circle (1.5pt);
							}
						}	
						\fill (0,0) circle (2.5pt);			
					\end{scope}
					\begin{scope}[shift={(-3,-4.5)}]
						\node at (-0.6,1.5) {$\NP_{\sigma_{3,4}}$:};					
						\fill[lightgray] (0,0) -- (0,2) -- (2,2) -- (2,1) -- (1,0) -- cycle;
						\draw[thick, black] ({0-\ltl},0) -- ({2+\ltl},0);
						\draw[thick, black] ({0-\ltl},2) -- ({2+\ltl},2);
						\draw[thick, black] (0,{0-\ltl}) -- (0,{3+\ltl});
						\draw[thick, black] (2,{0-\ltl}) -- (2,{3+\ltl});
						\draw[thick, black] ({1-\ltl},{0-\ltl}) -- ({2+\ltl},{1+\ltl});
						
						\foreach \x in {0,...,2}{
							\foreach \y in {0,...,3}{
								\fill (\x,\y) circle (1.5pt);
							}
						}				
						\fill (0,0) circle (2.5pt);
					\end{scope}
					\begin{scope}[shift={(2.5,-4.5)}]
						\node at (-0.6,1.5) {$\NP_{\sigma_{4,5}}$:};					
						\fill[lightgray] (0,0) -- (0,1) -- (2,3) -- (2,1) -- (1,0) -- cycle;
						\draw[thick, black] ({0-\ltl},0) -- ({2+\ltl},0);
						\draw[thick, black] (0,{0-\ltl}) -- (0,{3+\ltl});
						\draw[thick, black] (2,{0-\ltl}) -- (2,{3+\ltl});
						\draw[thick, black] ({1-\ltl},{0-\ltl}) -- ({2+\ltl},{1+\ltl});
						\draw[thick, black] ({0-\ltl},{1-\ltl}) -- ({2+\ltl},{3+\ltl});
						
						\foreach \x in {0,...,2}{
							\foreach \y in {0,...,3}{
								\fill (\x,\y) circle (1.5pt);
							}
						}				
						\fill (0,0) circle (2.5pt);
					\end{scope}
					\begin{scope}[shift={(7,-4.5)}]
						\node at (-0.6,1.5) {$\NP_{\sigma_{1,5}}$:};					
						\fill[lightgray] (0,0) -- (0,1) -- (2,3) -- (2,1) -- (1,0) -- cycle;
						\draw[thick, black] ({0-\ltl},0) -- ({2+\ltl},0);
						\draw[thick, black] (0,{0-\ltl}) -- (0,{3+\ltl});
						\draw[thick, black] (2,{0-\ltl}) -- (2,{3+\ltl});
						\draw[thick, black] ({1-\ltl},{0-\ltl}) -- ({2+\ltl},{1+\ltl});
						\draw[thick, black] ({0-\ltl},{1-\ltl}) -- ({2+\ltl},{3+\ltl});
						
						\foreach \x in {0,...,2}{
							\foreach \y in {0,...,3}{
								\fill (\x,\y) circle (1.5pt);
							}
						}				
						\fill (0,0) circle (2.5pt);
					\end{scope}
				\end{tikzpicture}
			\end{center}
			\caption{The Newton polytopes with respect to different $\sigma$. The fat point indicates the origin.}
			\label{fig:NPsigmas}
		\end{figure}
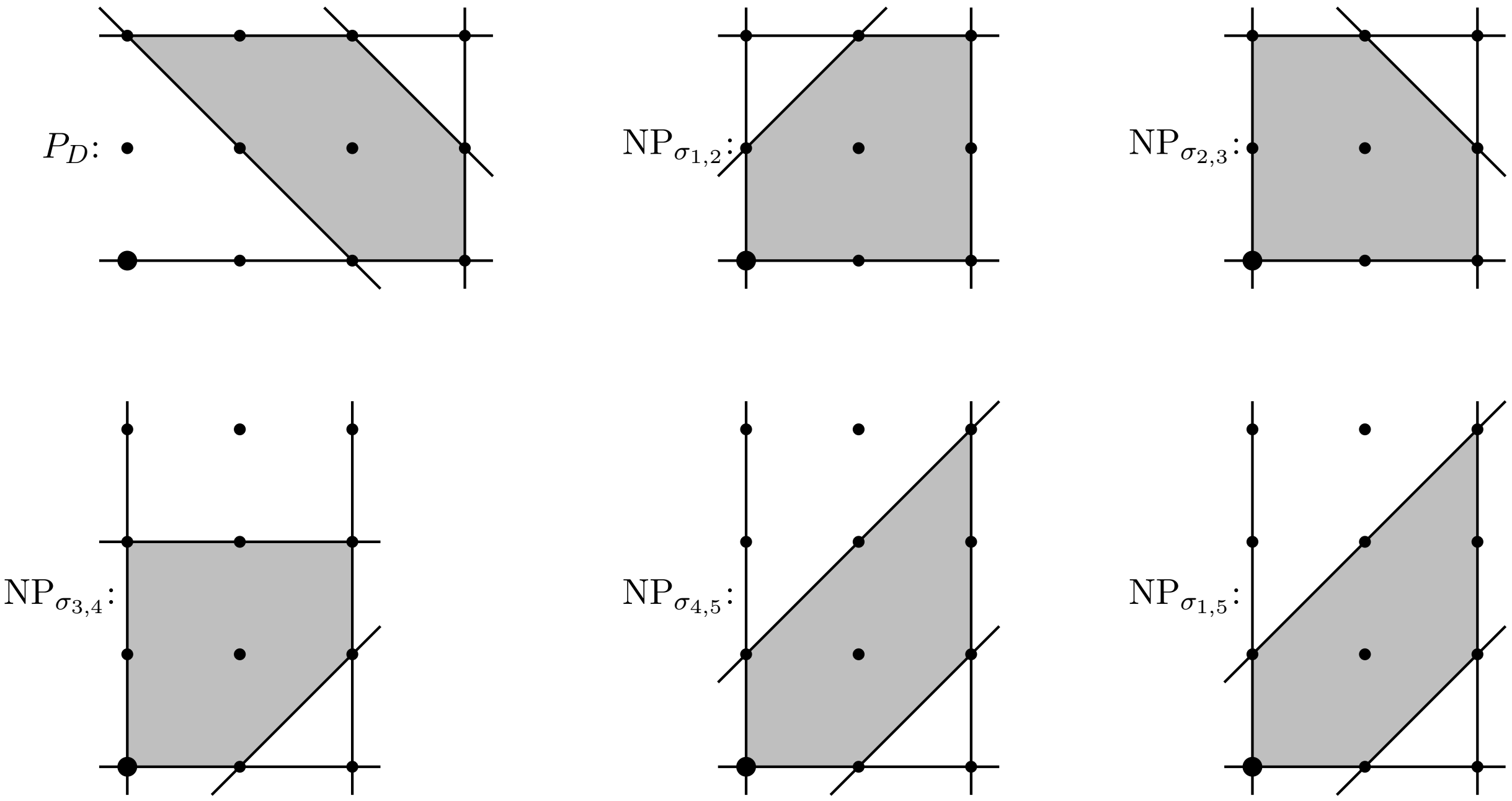
	\end{example}
	One advantage of having sets in $\R^n$ is that there is a well-defined notion of odd lattice points -- tuples in $\Z^n$ with only odd entries.\\
	Moreover, as we will show in the next lemma, the Newton polytope with respect to $\sigma$ is the same for rationally equivalent divisors.
	\begin{lemma}
		Let $D\in\Div_T(\XSig)$, $\sigma,\sigma'\in\Sigma(n)$ and $D' = D + \divi(\chi^\mu)$, i.e.\@ $D\sim D'$. Then the following statements are true.
		\begin{itemize}
			\item $\NP_\sigma(D) = \NP_\sigma(D')$.
			\item $L_{D,\sigma} = L_{D',\sigma}$.
			\item $\NP_\sigma(D)\subset \R_{\geq0}^n$.
			\item $\vol(\NP_\sigma(D)) = \vol(\NP_{\sigma'}(D))$.
		\end{itemize}
	\end{lemma}
	\begin{proof}
		Recall that $P_{D'} = P_D - \mu$. Further, we know that $b_\sigma^{D'} = b_\sigma^{D} + b_\sigma^{\divi \chi^\mu} = b_\sigma^D - \mu$. We compute
		\begin{align*}
			&\{(\langle m-b_\sigma^D, v_{i_1}\rangle,\dots,\langle m-b_\sigma^D, v_{i_n}\rangle)\mid m\in P_D\}\\
			&= \{(\langle m-\mu -(b_\sigma^D-\mu), v_{i_1}\rangle,\dots,\langle m-\mu -(b_\sigma^D-\mu), v_{i_n}\rangle)\mid m\in P_D  = P_{D'}-\mu\}\\
			&= \{(\langle m -b_\sigma^{D'}, v_{i_1}\rangle,\dots,\langle m-b_\sigma^{D'}, v_{i_n}\rangle)\mid m\in P_{D'}\}.
		\end{align*}
		So the convex hull is taken over the same set for $\NP_{\sigma}(D)$ and $\NP_\sigma(D')$, proving the first point. The second point follows similarly.\\
		It is $\langle m-b_\sigma^D, v_i\rangle = \langle m, v_i\rangle - \langle b_\sigma^D, v_i\rangle \geq -a_i + a_i = 0$ for $i\in\sigma(1)$ as $m\in P_D$ and by definition of $b_\sigma^D$, so the third point follows.\\
		The last point is due to the fact, that $\NP_{\sigma'}(D)$ gets mapped to $\NP_{\sigma}(D)$ by an invertible $\Z$-linear map and a shift.
	\end{proof}
	The fact that the Newton polytope with respect to $\sigma$ lies in the positive orthant means that the global sections of $\mcOX(D)$ are polynomials in our chosen coordinates. The definition of combinatorial orientability given in \cite{JPPbezout} is based on this coordinate-wise interpretation of global sections. Our next goal is to show that in the relatively orientable case we actually don't have odd boundary points on $\NP_{\sigma}(D)$. The next Lemma is the main missing ingredient.
	\begin{lemma}
		\label{lem:bsdbsoToOdd}
		Let $D\in\Div_T(\XSig)$ and $\sigma\in\Sigma(n)$. Then 
		$$\psisd((b_\sigma^D + \bso)+2M) = (1,\dots,1)+2\Z^n.$$
	\end{lemma}
	\begin{proof}
		Let $k\in M$. We compute
		\begin{align*}
			\psisd((b_\sigma^D + \bso)+2k) &= (\langle (b_\sigma^D + \bso)+2k-b_\sigma^D, v_{i_j}\rangle)_{j=1,...,n} \\
			&=(\langle  \bso,v_{i_j}\rangle+2\langle k, v_{i_j}\rangle)_{j=1,...,n}\\
			&=(1+2\langle k, v_{i_j}\rangle)_{j=1,...,n}
		\end{align*}
		As $v_{i_1},\dots,v_{i_n}$ is a basis of $N$, the claim follows.
	\end{proof}
	\begin{theorem}
		\label{thm:orientAndOddBoundary}
		Let $\XSig$ be a smooth proper toric variety and $D\in\Div_T(\XSig)$. Then the following statements are true.
		\begin{itemize}
			\item $D$ is relatively orientable if and only if  $L_{D,\sigma}$ contains no odd points for some/all $\sigma\in\Sigma(n)$.
			\item If $P_D$ is a lattice polytope, $D$ being relatively orientable implies that $\NP_{\sigma}(D)$ has no odd boundary points for some/all $\sigma\in\Sigma(n)$.
			\item If further all edges of $P_D$ have lattice length $\geq n-1$ and $D$ is ample, then $\NP_{\sigma}(D)$ having no odd boundary points for some  $\sigma\in\Sigma(n)$ implies that $D$ is relatively orientable.
		\end{itemize}
		In particular, if $D$ is the sum of $n$ ample divisors, $D$ is relatively orientable if and only if $\NP_{\sigma}(D)$ has no odd boundary points for some/all $\sigma\in\Sigma(n)$		
	\end{theorem}
	\begin{proof}
		By Propositions \ref{prop:theEquivalence} and \ref{prop:theMissingPoint}, $D\in\Div_T(\XSig)$ is relatively orientable if and only if $(b_\sigma^D + \bso)+2M\not\in\pi_2(L_D\cap M)$ for some/all $\sigma\in\Sigma(n)$. The first point follows now by applying Lemma \ref{lem:bsdbsoToOdd}.\\
		The second point follows from the  first as $P_D$ is a lattice polytope  and $\partial\NP_{\sigma}(D)\subseteq L_{D,\sigma}$.\\
		The last point is Lemma \ref{lem:NecCondForLP=DelP} combined with the first point and the fact that for ample divisors $L_D$ equals  $L_{P_D}$ and $P_D$ is a smooth full-dimensional lattice polytope.
	\end{proof}
	\begin{example}
		In Example \ref{ex:blowUpP3}, $\NP_{\sigma_1}(E)$ has no odd boundary point, while $\NP_{\sigma_1'}(D)$ has one, as computations in Macaulay2 show. You can find the computations \href{https://ntiggemann.github.io/Examples.m2}{here}.
	\end{example}
	\begin{example}
		Let $\XSig$ be the third Hirzebruch surface and $D = D_1+D_2+D_3+D_4$. In that case, $D$ is relatively orientable, $P_D$ not a lattice polytope and $\NP_{\sigma_1'}(D)$ has an odd boundary point, as can be seen in Figure \ref{fig:NPsigma1strich}.
		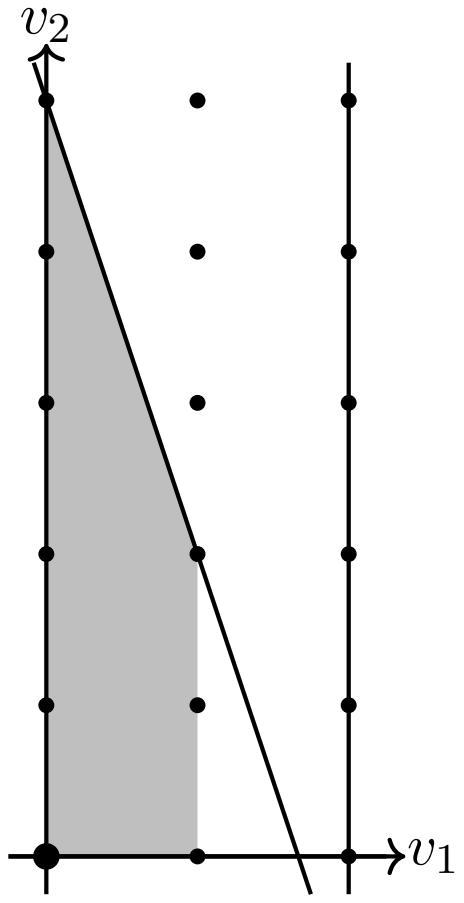
\begin{figure}
			\begin{center}
				\begin{tikzpicture}
					\pgfmathsetmacro{\ltl}{0.25}
					\fill[lightgray] (0,0) -- (0,5) -- (1,2) -- (1,0) -- cycle;
					\draw[->,thick] (0,{0-\ltl}) -- (0,{5+1.5*\ltl}); 	
					
					\draw[thick, black] ({0-\ltl},0) -- ({2+\ltl},0);
					
					\draw[->,thick] ({0-\ltl},0) -- ({2+1.5*\ltl},0);
					\draw[thick, black] ({0-\ltl/3},{5+\ltl}) -- ({1+\ltl/3+2/3},{0-\ltl});
					\draw[thick, black] (2,{0-\ltl}) -- (2,{5+\ltl});
					
					\node[anchor=south] at (0,5+\ltl) {$v_2$}; 					
					\node[anchor=west] at (2+\ltl,0) {$v_1$};
					
					\foreach \x in {0,...,2}{
						\foreach \y in {0,...,5}{
							\fill (\x,\y) circle (1.5pt);
						}
					}				
					\fill (0,0) circle (2.5pt);
				\end{tikzpicture}
				\caption{An example of a relatively orientable divisor $D$ for which $\NP_{\sigma}(D)$ does have an odd boundary point. $\NP_{\sigma_1'}(D)$ is indicated by the gray area, the union of the lines is $L_{D,\sigma_1'}$.}
				\label{fig:NPsigma1strich}
			\end{center}
		\end{figure}
	\end{example}
	\begin{example}
		Consider the smooth proper toric variety defined by the fan with rays $(0,1)$, $(1,1)$, $(1,0)$, $(0,-1)$, $(-1,0)$ and the torus invariant, not relatively orientable, divisor $D = -D_2+3D_3+D_4+D_5$. $P_D = \conv\{(1,0),(0,1),(1,1)\}$ is a lattice polytope with $\pi_2(L_{P_D}\cap \Z^2)\not = \Z^2/2\Z^2$. However, $\NP_{\sigma_{1,5}}(D)$ has an odd boundary point and $\pi(\NP_{\sigma_{1,5}}(D))\not = \Z^2/2\Z^2$.
	\end{example}
	\begin{remark}
		In the Macaulay2 \cite{M2} package \emph{LatticePolytopes} \cite{LatticePolytopesSource} (which uses the packages \emph{Polyhedra} \cite{PolyhedraSource} and \emph{NormalToricVarieties} \cite{NormalToricVarietiesSource}) the $41$ smooth $2$-polytopes with up to $12$ lattice points and the $103$ smooth $3$-polytopes with up to $16$ lattice points are listed. There are $6$ combinatorially orientable polytopes among the $41$ smooth $2$-polytopes. Out of the $103$ smooth $3$-polytopes, $30$ are combinatorially orientable. For $71$ of the $3$-polytopes it is $\pi_2(\partial P \cap \Z^3)\not= \Z^3/2\Z^3$. Finally, $69$ satisfy that applying ''$\psisd$`` (where $b_\sigma^D$ is the vertex and the basis the one of the directions of the edges at this vertex) for each corner never yields a polytope with an odd boundary point. The computations were run on Macaulay2 version 1.26.05. The source code for these computations can be found \href{https://ntiggemann.github.io/Examples.m2}{here}.\\
		Notice that by our definition of combinatorial orientability we ''lost`` some of the polytopes for which all $\NP_{\sigma}(D)$ have no odd boundary point. However, the approach to take the latter property as a definition does not make sense for non-smooth polytopes, as their edges don't have to give $\Z$-bases. Only considering the actual boundary points when applying $\pi_2$ is too weak as a condition as the next examples shows.
	\end{remark}
	\begin{example}
		\label{ex:BadGuy}
		Consider the polytope $$P \coloneqq \conv\{(0,0,0),(0,1,0),(1,1,0),(0,0,1),(0,2,1),(2,2,1)\}$$ (this is the smooth $3$-polytope with index $71$ in the package mentioned in the previous remark). It is not combinatorially orientable but $\pi_2(\partial P \cap \Z^3)\not= \Z^3/2\Z^3$. In Figure \ref{fig:BadGuy} we can see a plot of $P$. The red dot is an odd point with respect to the basis of the adjacent edges of the vertex in the origin. If one chooses e.g.\@ the basis of the adjacent edges of the vertex $(0,0,1)$, there are no odd boundary points.
		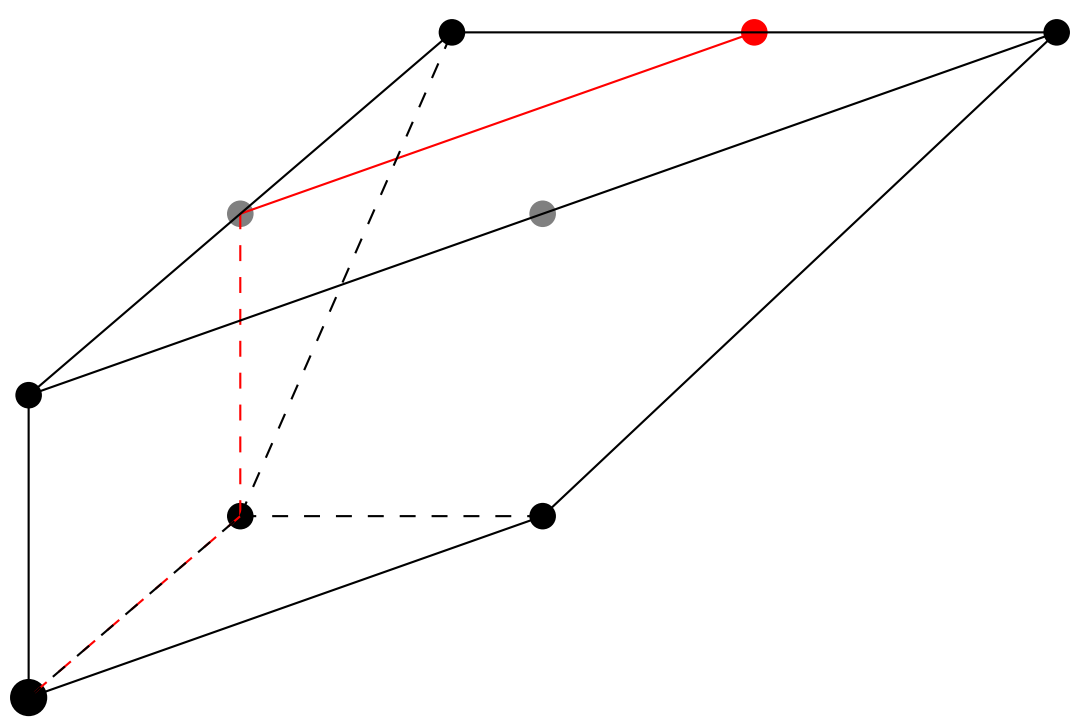
\begin{figure}
			\begin{center}
				\begin{tikzpicture}[scale=2]
					\pgfmathsetmacro{\ptSize}{0.5}
					
					\coordinate (v1) at (0,0);
					\coordinate (v2) at (0.7,0.6);
					\coordinate (v3) at (1.7,0.6);
					\coordinate (v4) at (0,1);
					\coordinate (v5) at (1.4,2.2);
					\coordinate (v6) at (3.4,2.2);

					\fill[gray] (0.7,1.6) circle (2.5*\ptSize pt);
					\fill[gray] (1.7,1.6) circle (2.5*\ptSize pt);

					\fill (v1) circle (3.5*\ptSize pt);
					\fill (v2) circle (2.5*\ptSize pt);
					\fill (v3) circle (2.5*\ptSize pt);
					\fill (v4) circle (2.5*\ptSize pt);
					\fill (v5) circle (2.5*\ptSize pt);
					\fill (v6) circle (2.5*\ptSize pt);
					
					\coordinate (badPoint) at (2.4,2.2);
					\fill[red] (badPoint) circle (2.5*\ptSize pt);
					\draw[red] (badPoint) -- (0.7,1.6);
					\draw[red,dashed] (v2) -- (0.7,1.6);
					\draw[red,dashed] (v2) -- (v1);
					
					\draw[black] (v1) -- (v4);
					\draw[black,dashed] (v1) -- (v2);
					\draw[black,dashed] (v2) -- (v5);
					\draw[black] (v5) -- (v4);
					\draw[black] (v1) -- (v3);
					\draw[black] (v3) -- (v6);
					\draw[black] (v6) -- (v4);
					\draw[black,dashed] (v2) -- (v3);
					\draw[black] (v5) -- (v6);
				\end{tikzpicture}
				\caption{The polytope mentioned in  Example \ref{ex:BadGuy}. The red point has coordinates $(1,1,1)$ with respect to the basis defined by the edges at the origin, the fat point. Clearly, $(1,0,0)\not\in\pi(\partial P \cap M)$ and $\pi(L_P\cap M) = (\Z/2\Z)^3$.}
				\label{fig:BadGuy}
			\end{center}
		\end{figure}
	\end{example}
\section{Nisnevich coordinates on toric varieties}
\label{section:ExplicitRelOrientOnTorVar}
	For the $\Al$-BKK theorem, we want to compute the Poincaré Hopf Euler number $n^{\Al}(V)$ of the relatively oriented vector bundle $V = \mcOX(E_1)\oplus\cdots\oplus\mcOX(E_n)\to \XSig$ which is defined to be the sum of local indices $\ind_x\sigma$ at the zeros of $\sigma$, see \cite{BachmannWickelgren}. In this section, we prove a formula for the local indices; we give local Nisnevich coordinates compatible with the relative orientation and compute the $\Al$-degree in those coordinates. While this formula was effectively already used in \cite{JPPbezout}, this is, to our knowledge, the first time an explicit proof of it is given.\\
	We refer the reader to \cite{pauli2026pcmilecturenotesmotivic} for an introduction to the $\Al$-degree and the Poincaré Hopf Euler number $n^{\Al}(V)$ of a vector bundle $V$.\\
	We closely follow \cite{McKean} and adapt the strategies used there for $\PP^n$ to our case of smooth proper $\XSig$. We give Nisnevich coordinates and compatible local trivializations for sums of torus invariant line bundles. Recall that we fix a numbering of the rays of $\Sigma$, $\Sigma(1) = \{\rho_1,\dots,\rho_s\}$.\\
	Recall that, since $\Sigma$ is smooth and complete, the minimal generators $\{v_{i_1},\dots,v_{i_n}\} = \{v_\rho\}_{\rho\in\sigma(1)}$ of the rays of a maximal cone $\sigma\in\Sigma(n)$ form a basis of $N$. We denote its dual basis by $B_\sigma = \{m_{i_1}^\sigma,\dots,m_{i_n}^\sigma\} = \{m^\sigma_\rho\}_{\rho\in\sigma(1)}$. Further, recall that the canonical line bundle on $\XSig$ is $\omega_{X_\Sigma/\kk}\cong\mcOX(-\sum_{\rho}D_\rho)$ and that we write $\{\bso\}_{\sigma\in\Sigma(n)}$ for the Cartier data of $-\sum_{\rho}D_\rho$. Explicitly, $\bso = \sum_{m\in B_\sigma} m$. In this section, we will work with $\Hom(\det \mct\XSig,\det V) = \det V \otimes \omega_{X_\Sigma/\kk}$ in the definition of relative orientability given in the introduction.
	\begin{definition}
		Let $p\colon V\to \XSig$ be an algebraic vector bundle of rank $n$. 
		\begin{itemize}
			\item Let $x\in \XSig$ be a closed point and let $U$ be a Zariski open neighborhood of $x$. An étale map $\psi\colon U \to \A_\kk^n$ which induces and isomorphism of residue fields at $x$ is called \emph{Nisnevich coordinates}.
			\item Assume $p$ is relatively oriented with relative orientation $\rho$. Let $x\in \XSig$ be a closed point and let $\psi$ be Nisnevich coordinates of $x$ such that there is a local trivialization of $V|_U$. We say that Nisnevich coordinates around $x$ and a trivialization $V|_U\cong U\times \A_\kk^n$ are \emph{compatible} with $\rho$ if the image of the distinguished section of $\Hom(\det \mct\XSig|_U,\det V|_U)$ taking the distinguished basis of $\det \mct\XSig|_U$ to $\det V|_U$ under $\rho$ can be written as $l\otimes l$, where $l$ is a section in $\mcl$.
		\end{itemize}
	\end{definition}
	\begin{definition}
		We define $\sigma_0\in\Sigma(n)$ to be the cone with $\sigma_0(1) = \{\rho_1,\dots,\rho_n\}$ and write $d^\sigma_\tau\in\{\pm 1\}$ for the determinant of the base change from $B_\sigma$ to $B_\tau$ for $\tau,\sigma\in\Sigma(n)$.\\
		We define twisted coordinates $\{(U_\sigma,\wivasig)\}_{\sigma\in\Sigma(n)}$ via
		$$\wivasig\colon U_\sigma \to \A^n_\kk;\ \wivasig([\bfx]) = (d^\sigma_{\sigma_0}\bfx^{m_{i_1}^\sigma},\dots,\bfx^{m_{i_n}^\sigma}).$$
	\end{definition}
	Note that conveniently $d^\tau_\sigma = d^\sigma_\tau = d^\pi_\tau d^\sigma_\pi = d_\pi^\tau d^\sigma_\pi$, for $\sigma,\tau,\pi\in\Sigma(n)$.
	We define $\bigwedge_{i\in \sigma(1)}$ to be the wedge product in the order of the numbering of the rays (which we have fixed in the beginning).
	\begin{proposition}
		The twisted coordinates are Nisnevich coordinates. Moreover, $\wivasig$ determines the distinguished basis element $d^\sigma_{\sigma_0}\partial_\sigma\coloneqq d^\sigma_{\sigma_0}\bigwedge_{i\in\sigma}\frac{\partial}{\partial \bfx^{m_i^\sigma}}$ of $\det \mct X_\Sigma|_{U_\sigma}$ with transition functions $\det g_{\sigma,\tau} = \bfx^{\bso-\bst}$.
	\end{proposition}
	\begin{proof}
		By construction, $\wivasig$ is an isomorphism, so it is étale and induces an isomorphism $\kappa(p)\cong \kappa(\widetilde{\varphi}_\sigma(p))$ for every $p\in U_\sigma$, and thus it defines Nisnevich coordinates.\\
		Recall that $\mct \A^n$ has the standard trivializations $\{\diff{t_1},\dots,\diff{t_n}\}$.
		Since $\wivasig$ induces an isomorphism $\mct \XSig|_{U_\sigma} \cong \mct \A^n$, we may pull back the standard trivializations $\mct \A^n \to \A^n$ by $\wivasig$ to obtain the twisted trivializations $\{\dsso \partial_{i_1/\sigma},\dots,\partial_{i_n/\sigma}\},\ i_1<\dots<i_n\in\sigma(1),$ where $\partial_{i/\sigma} = \diff{\bfx^{m^\sigma_i}}$. It follows that $\det \mct \XSig|_{U_\sigma}$ is trivialized by $\dsso \bigwedge_{i\in\sigma(1)} \partial_{i/\sigma}$.\\		
		Finally, we consider the transition functions $\det g_{\sigma,\tau}\colon \det\mct\XSig|_{U_\tau} \to \det\mct\XSig|_{U_\sigma}$. Let from now on $z_k = \bfx^{m^\sigma_{i_k}}$ and $y_k = \bfx^{m^\tau_{i_k}}$.
		Let $\lambda_{j,k}$ be such that $m^\tau_{i_k} = \sum_{j=1}^{n}\lambda_{j,k}m^\sigma_{i_j}$ and $A\in\Gl_n(\Z)$ the matrix with entries $\lambda_{j,k}$, that is, the base-change matrix from $B_\tau$ to $B_\sigma$. In particular, $\det A = d^\tau_\sigma = d_\tau^\sigma$. We get that $y_k = z_1^{\lambda_{1,k}}\cdot \dots\cdot z_n^{\lambda_{n,k}}$ and compute
		\begin{align*}
			\partial_\sigma = \bigwedge_{i\in \sigma(1)} \diff{\bfx^{m^\sigma_i}} &= \bigwedge_{i=1}^n \diff{z_i}\\
			&= \bigwedge_{i=1}^n \sum_{j=1}^{n} \diff[y_j]{z_i}\diff{y_j}\\
			&= \bigwedge_{i=1}^n \sum_{j=1}^{n} \left(\frac{\lambda_{i,j}y_j}{z_i}\right)\diff{y_j}\\
			&= \det\left(\frac{\lambda_{i,j}y_j}{z_i}\right) \cdot\bigwedge_{i=1}^n \diff{y_i}\\
			&= \det(\diag(y_1,\dots,y_n)) \cdot \det(A) \cdot \det(\diag (z_1^{-1},\dots,z_n^{-1}))\cdot\partial_\tau\\
			&= d^\sigma_\tau y_1\cdot \dots \cdot y_n \cdot z_1^{-1}\cdot \dots \cdot z^{-1}_n\cdot\partial_\tau\\
			&= d^\sigma_\tau \bfx^{\sum_{j=1}^{n} m_{i_j}^\tau - \sum_{j=1}^{n} m_{i_j}^\sigma}\cdot\partial_\tau\\
			&= d^\sigma_\tau \bfx^{\bst-\bso}\cdot\partial_\tau
		\end{align*}
		Multiplying this equation with $\dsso$ we get 
		$$\dsso \partial_\sigma = \bfx^{\bst-\bso}\dsso[\tau]\partial_\tau.$$
	\end{proof}
	Let $D\in\Div_T(\XSig)$ with Cartier data $\{b_\sigma^D\}_{\sigma\in\Sigma(n)}$. Then we have $D|_{U_\sigma} = \divi (\bfx^{-b_\sigma^D})|_{U_\sigma}$ for all $\sigma\in\Sigma(n)$ and in conclusion, we know that $\bfx^{b_\sigma^D}$ defines a local trivializing frame of $\mcOX(D)$ on $U_\sigma$. The transition functions from $U_\sigma$ to $U_\tau$ on $U_\sigma\cap U_\tau$ are $\bfx^{b_\tau^D - b_\sigma^D}$. This gives us the next proposition.\\
	From now on, let $E_1,\dots,E_n\in\Div_T(\XSig)$ with Cartier data $\{b_\sigma^1\}_{\sigma\in\Sigma(n)},\dots,\{b_\sigma^n\}_{\sigma\in\Sigma(n)}$ and $D = \sum_{i=1}^{n}E_i$ with Cartier data $\{b_\sigma^D\}_{\sigma\in\Sigma(n)}$.
	\begin{proposition}
		The local trivialization $\bfx^{b_\sigma^1} \oplus \dots \oplus \bfx^{b_\sigma^n}$ of $\bigoplus_{i=1}^n\mcOX(E_i)|_{U_\sigma}$ determines the distinguished basis element $\bfx^{\sum_{i=1}^{n}b_\sigma^i} = \bfx^{b_\sigma^D}$ of $\det \bigoplus_{i=1}^n\mcOX(E_i)|_{U_\sigma} = \mcOX(D)|_{U_\sigma}$ with transition functions $\det h_{\sigma, \tau} = \bfx^{\sum_{i=1}^{n}b_\sigma^i - b_\tau^i} = \bfx^{b_\sigma^D - b_\tau^D}$.
	\end{proposition}
	Let $E_1,\dots,E_n\in\Div_T(\XSig)$ be such that there is a $D'\in\Div_T(\XSig)$ with $D=\sum_{i=1}^{n}E_i \sim D' =  \sum_{\rho} a_\rho D_\rho$ and all $a_\rho$ odd. Let $E \coloneqq \sum_{\rho} (a_\rho-1)D_\rho$ and $E/2 = \sum_{\rho} \frac{a_\rho-1}{2}D_\rho$. Then $\mcOX(D) \otimes \omega_{X_\Sigma/\kk} \cong \mcOX(E) \cong \mcOX(E/2)^{\otimes 2}$ is a square, as we need it to be for the relative orientation. As the Cartier data $\{b_\sigma^E\}_{\sigma\in\Sigma(n)}$ of the divisor $E$ solve $\langle b_\sigma^E, v_\rho\rangle = a_\rho -1\in2\Z$ for $\rho \in\sigma(1)$ and the $v_\rho$ are primitive, we have $\frac{1}{2}b_\sigma^E\in M$ for each $\sigma\in\Sigma(n)$.\\
	We will give a relative orientation of $\bigoplus_{i=1}^n\mcOX(E_i)$ and show that the local trivializations of $\bigoplus_{i=1}^n\mcOX(E_i)$ discussed above are compatible with this relative orientation and the Nisnevich coordinates from our twisted cover.
	For our relative orientation, we give an isomorphism 
	$$\psi\colon \mcOX(E/2)^{\otimes 2} \to \Hom(\det \mct X_\Sigma, \mcOX(D))\quad (\cong \omega_{X_\Sigma/\kk}\otimes \mcOX(D))$$
	by defining $\psi|_{U_\sigma}$ for each $\sigma$. Since $\mcOX(E/2)^{\otimes 2}|_{U_\sigma}$ is generated by $\bfx^{\frac{1}{2}b_\sigma^E}\otimes \bfx^{\frac{1}{2}b_\sigma^E}$, it suffices to define $\alpha_\sigma \coloneqq \psi|_{U_\sigma}(\bfx^{\frac{1}{2}b_\sigma^E}\otimes \bfx^{\frac{1}{2}b_\sigma^E})$, which is a homomorphism from $\det \mct X_\Sigma$ to $\mcOX(D)$. Both are one dimensional and above we have given local trivializations, so we may define $\alpha_\sigma$ to be the homomorphism taking $\dsso\partial_\sigma$ to $\bfx^{b_\sigma^D}$.\\
	To show that $\psi$ is well-defined, we need to show that on $U_\sigma\cap U_\tau$ the maps $\psi|_{U_\sigma}$ and $\psi|_{U_\tau}$ differ by the transition function $\bfx^{(b_\sigma^E - b_\tau^E)/2}\otimes\bfx^{(b_\sigma^E- b_\tau^E)/2}\colon \mcOX(E/2)^{\otimes 2}|_{U_\tau}\to \mcOX(E/2)^{\otimes 2}|_{U_\sigma}$. In other words, we need to show that $\alpha_\sigma = \bfx^{b_\sigma^E- b_\tau^E}\alpha_\tau$ on $U_\sigma\cap U_\tau$. To this end, let $\det g_{\sigma,\tau}$ and $\det h_{\sigma,\tau}$ be the transition functions given above. We compute
	\begin{align*}
		\alpha_\sigma\circ \det g_{\sigma,\tau}(\dsso[\tau]\partial_\tau) &= \alpha_\sigma(\bfx^{\bst - b_\sigma^{\omega}}\dsso \partial_\sigma)\\
		&= \bfx^{b_\sigma^{\omega}-\bst} \cdot \bfx^{b_\sigma^D}\\
		&= \bfx^{b_\sigma^{\omega}-\bst}\cdot  \det h_{\sigma,\tau}\cdot  \bfx^{b_\tau^D}\\
		&= \bfx^{b_\sigma^{\omega}-\bst} \cdot \bfx^{b_\sigma^D-b_\tau^D}\cdot  \bfx^{b_\tau^D}\\
		&= \bfx^{b_\sigma^{\omega}-\bst} \cdot \bfx^{b_\sigma^D-b_\tau^D}\cdot  \alpha_\tau(\dsso[\tau]\partial_\tau)\\
		&= \bfx^{b_\sigma^E - b_\tau^E}\cdot  \alpha_\tau(\dsso[\tau]\partial_\tau)
	\end{align*}
	\begin{lemma}
		The local trivializations $\bfx^{b_\sigma^1} \oplus \dots \oplus \bfx^{b_\sigma^n}$ of $\bigoplus_{i=1}^n\mcOX(E_i)|_{U_\sigma}$ are compatible with the Nisnevich coordinates $\{(U_\sigma,\wivasig)\}$ and the relative orientation $(\mcOX(E/2),\psi)$ of $\bigoplus_{i=1}^n\mcOX(E_i) \to \XSig$.
	\end{lemma}
	\begin{proof}
		By construction, $\alpha_\sigma$ is the element of $\Hom(\det \mct \XSig|_{U_\sigma}, \mcOX(D)|_{U_\sigma})$ taking the distinguished basis element $\dsso\partial_\sigma$ of $\det \mct \XSig|_{U_\sigma}$ to the distinguished basis element $\bfx^{b_\sigma^D}$ of $\mcOX(D)$. The relative orientation $(\mcOX(E/2),\psi)$ was built such that $\psi|_{U_\sigma}(\bfx^{\frac{1}{2}b_\sigma^E}\otimes \bfx^{\frac{1}{2}b_\sigma^E})\eqqcolon \alpha_\sigma$, so $\alpha_\sigma$ is a tensor square in $\mcOX(E/2)^{\otimes 2}|_{U_\sigma}$.
	\end{proof}
	\begin{lemma}
	\label{lem:degAsTrace}
		Let $s_1\in \mcOX(D_1), \dots, s_n\in \mcOX(D_n)$, $D_i\in\Div_T(\XSig)$, define hypersurfaces of $\XSig$ that intersect transversely at a point $p$. Let $f_l\in S_{[D_l]}$ be the Cox-homogeneous polynomial of $s_l$ and $F_l^\sigma$ its dehomogenization. Write $f^\sigma_l \coloneqq F_l^\sigma(d_{\sigma_0}^\sigma \bfx^{m_{i_1}},\dots,\bfx^{m_{i_n}})$. Let $$J_\sigma = \det \left(\frac{\partial f^\sigma_l}{\partial \bfx^{m_{i_j}}}\right)_{l,j=1,\dots,n}$$
		and let $\Tr_{\kappa(p)/\kk}\colon \GW(\kappa(p))\to \GW(\kk)$ be the trace map on bilinear forms obtained by post composing with the field trace $\kappa(p)\to \kk$. Then
		$$\ind_p(s_1,\dots,s_n) = \Tr_{\kappa(p)/\kk}\langle J_\sigma(\wivasig(p))\rangle.$$
	\end{lemma}
	\begin{proof}
		The local trivialization $(f^\sigma_1,\dots,f^\sigma_n)$ of the section $(s_1,\dots,s_n)$ is compatible with our chosen relative orientation and Nisnevich coordinates $\{(U_\sigma,\wivasig)\}$ as we have proven above. Since $p$ is a simple zero of the section, \cite{KassWickelgrenCubicSurface} (p. 17) gives us the claimed result.
	\end{proof}
\section{Quadratically enriched tropical geometry}
\label{section:CombRelOrientInTropGeo}
	In this section, we will see that the definition of combinatorial orientability translates nicely to the tropical case. In particular, we will see that this definition is well-defined even for balanced polyhedral complexes.
\subsection{Tropical numbers and tropical geometry}
\begin{figure}[htb]
		\centering
		\begin{tikzpicture}
			\draw [black] (1.5,1.5) -- (0,0);
			\draw [black] (-2,0) -- (0,0);
			\draw [black] (0,-2) -- (0,0);
			
			\filldraw[color = black, fill=black] (0.75,0.75) circle (1.5pt);
			\node[black] at (0.3,0.85) {$(0,2)$};
			
			\node[black] at (-2.05,0.2) {$\sigma_1$};
			\node[black] at (0.2,-2) {$\sigma_2$};
			\node[black] at (1.55,1.55) {$\sigma_3$};
		\end{tikzpicture}
		\caption{A tropical line in $\R^2$ with rays $\sigma_1,\sigma_2,\sigma_3$.}
		\label{Fig:TropLine}
	\end{figure}
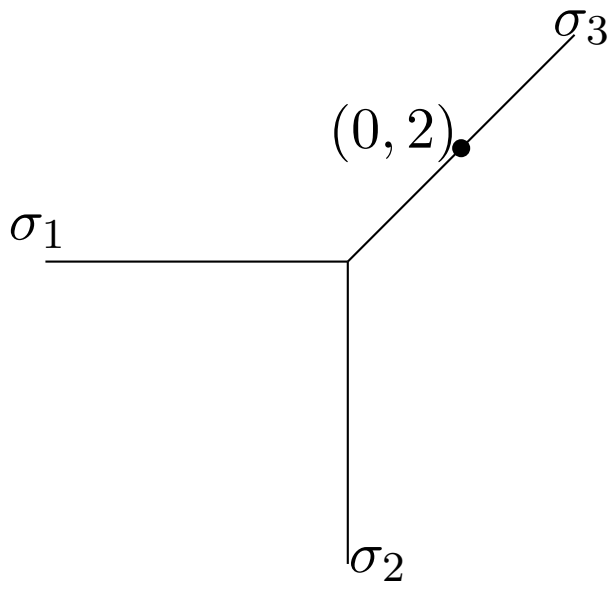

	We give a brief introduction to tropical geometry. For more details, we refer the reader to \cite{Payne_2008}, \cite{MikhalkinRau_TropicalGeometry} or \cite{SturmfelsMaclaganIntroTropGeo}. We define the \emph{tropical numbers} to be $\T = (\T,\oplus,\odot) = (\R\cup\{-\infty\}, \max,+)$. The tropical numbers form a semi-field, i.e.\@  they satisfy the axioms of a field except that there are no additive inverses.
	\begin{definition}
		For a field $\kk$ we define the \emph{field of Puiseux series} $\KK \coloneqq  \kk\{\{t\}\} = \bigcup_{n=1}^\infty \kk((t^{\frac{1}{n}}))$.\\
		It is a valued field with valuation map 
		$$\val\colon \KK^\times \to \R; \sum_{k=k_0}^{\infty}c_kt^{k/n}\mapsto k_0/n,$$
		we define $\nu\coloneqq-\val$.
	\end{definition}
	In this subsection, we assume $\KK$ to be algebraically closed.
	\begin{remark}
		Let $\kk$ be algebraically closed. If $\charact \kk = 0$, then $\KK$ is algebraically closed. If $\charact \kk = p > 0$, it is not. However, for $f\in\KK[x]$ of degree smaller than $p$, $\KK$ contains all zeros of $f$. I.e., as long as the degree of the polynomials we consider is not too big, $\KK$ is essentially algebraically closed for positive characteristic, too \cite[Theorem 5.14]{Abhyankar}.
	\end{remark}
	\begin{remark}
		The tropical numbers can be defined to be $(\R\cup\{\infty\},\min, +)$, too. This is equivalent to our definition via sign change.\\
		The constructions and statements in this subsection are true for a general algebraically closed valued field. However, we have to restrict ourselves to the Puiseux series from the next subsection onward anyways. 
	\end{remark}
	Let $N\cong\Z^n$ be a lattice with dual lattice $M = \Hom(N,\Z)$, and $\G_m^n = \Spec \KK[M]$ be a torus. Its $\KK$-rational points $\Homalg(\KK[M],\KK)$ are in a one-to-one correspondence to group homomorphisms $\Homgp(M,\KK^\times)$ via $\phi\mapsto (m\mapsto\phi(m))$. We define the \emph{tropicalization map} to be 
	$$\trop\colon \G_m^n(\KK) = \Homgp(M,\KK^\times)\to \Homgp(M,\R) = N_\R;\ \phi\mapsto \nu\circ\phi,$$
	and call $N_\R \cong \R^n$ the \emph{tropical torus}.\\
	Given a character $\chi^m$ and $\phi\in\Homgp(M,\KK^\times) = \G^n_m(\KK)$, and writing $x = \trop\phi = \nu\circ\phi\in N_\R = \Hom(M,\R)$, one computes that its tropicalization is 
	$\trop\chi^m(\phi) = \nu(\phi(m))=\langle m,x\rangle$. Given a Laurent polynomial $f = \sum_{i=0}^{k}a_i\chi^{m_i}$, the properties of the valuation map yield $\trop f = \max\{\nu (a_i) + \langle m_i,x\rangle\mid {i=0,\dots, k}\}$. We call such an expression a tropical (Laurent) polynomial and introduce the notation $\trop f = \bigoplus_{i=0}^k \nu (a_i) \odot \chi^{m_i}$. If a tropical polynomial $F$ is of the form $F = \trop f$, we call it the tropicalization of $f$.\\
	The \emph{tropical zero set $V^{\trop}(F)$ of a tropical polynomial $F$} is defined to be the set of all $x\in N_\R$ for which the maximum is attained at least twice. To each facet of a tropical zero set, we assign a weight: If the maximal monomials defining the facet $\sigma$ are $a_{m_1}\odot \chi^{m_1}$ and $a_{m_2}\odot \chi^{m_2}$, we assign the weight $[M\cap \R(m_1-m_2): \Z(m_1-m_2)]$ to it.
	\begin{remark}
		For all tropical polynomials $F$, we assume that there is a bijection between its monomial and the connected components of $N_\R\backslash V^{\trop}(F)$. I.e. we leave out monomials that do not contribute to the hypersurface.
	\end{remark}
	Tropical zero sets are balanced polyhedral complexes, and every balanced polyhedral complex in $N_\R$ of codimension one with positive weights is given by a tropical polynomial which is unique up to tropical multiplication with a tropical monomial, see e.g.\@ \cite[Theorem 2.4.10]{MikhalkinRau_TropicalGeometry}.\\
	Given a subvariety $Y$ of $\G_m^n$, its \emph{tropicalization} is defined to be $\trop Y = \overline{\trop Y(\KK)}\subseteq N_\R$, the topological closure of $\trop Y(\KK)$. If $Y$ is given by a Laurent polynomial $f\in\KK[M]$, it is $\trop Y = V^{\trop} (\trop f)$.
	\begin{example}
		Let $\KK = \C\{\{t\}\} = \bigcup_{n=1}^\infty \C((t^{\frac{1}{n}}))$ be the field of Puiseux series over $\C$. 
		Further, let $N = \Z^2$, we choose the standard basis $e_1 = (1,0), e_2 = (0,1)$. Then $M = \Hom(\Z^2,\Z)\cong \Z^2$ via $\varphi\mapsto (\varphi(e_1),\varphi(e_2))$ with dual basis $f_1 = (1,0), f_2 = (0,1)$. The pairing of $M$ and $N$ takes the form $\langle (m_1,m_2),(n_1,n_2)\rangle = m_1n_1+m_2n_2$ of the standard scalar product.\\
		With our choice of basis, $\G_m^2(\KK) = \Homgp(M,\KK^\times) \cong (\KK^\times)^2\eqqcolon T$ via $\varphi\mapsto (\varphi(e_1),\varphi(e_2))$. Tropicalization now simply maps a point $(f,g)\in T\cong\G_m^2(\KK)$ to $(-\val f,-\val g)\in\R^2\cong N_\R$.\\
		As an example of a subvariety, consider the line $\{(x,y)\in T\mid (t^{-2}+t)x + y + t^{1/2} = 0\} = V((t^{-2}+t)x + y + t^{1/2})$ in $T$. It contains e.g.\@ the point $(-1,t^{-2}- t^{1/2}+t)$, which tropicalizes to $(0,2)\in\R^2$.\\
		Simultaneously, $(t^{-2}+t)x + y + t^{1/2} = (t^{-2}+t)\chi^{(1,0)} + \chi^{(0,1)} + t^{1/2}\chi^{(0,0)}$ has tropicalization
		\begin{align*}
			&\nu (t^{-2}+t)\odot\chi^{(1,0)} \oplus \nu (1)\odot\chi^{(0,1)} \oplus \nu(t^{1/2})\odot\chi^{(0,0)} \\
			&= \max\{2+\langle (u,v),(1,0)\rangle, 0+\langle (u,v),(0,1)\rangle, -\frac{1}{2}+\langle (u,v),(0,0)\rangle\}\\ &= 2\odot u\oplus v \oplus -\frac{1}{2}.
		\end{align*} 
		The tropical line defined by this tropical polynomial can be seen in Figure \ref{Fig:TropLine}. Note that its vertex is at $(-\frac{5}{2},-\frac{1}{2})$, as the maximum there is attained thrice.\\
		In fact, for the point $(0,2)$, the maximum in the polynomial is attained twice, by $2\cdot u$ and $v$, so the point lies in the ray going to the upper right.
	\end{example}
	We define the \emph{Newton polytope of a tropical polynomial} $F = \bigoplus_{i=0}^k c_i \odot \chi^{m_i}$ to be $\NP(F)\coloneqq\conv\{m_1,\dots,m_k\}$. The \emph{dual subdivision $\DS(F)$ of $F$} is defined to be the refinement of $\NP(F)$ given by the projection to $M_\R$ of the boundary of the upper faces of the polyhedron $\conv(\{(m_i,c_i)\in M_\R\times\R\})$.\\
	There exists a one-to-one correspondence between vertices of $V^{\trop}(F)$ and connected components of $\NP(F)\backslash\DS(F)$, $l$-dimensional faces of  $V^{\trop}(F)$ and $n-l$-dimensional faces of $\DS(F)$, and connected components of $N_\R\backslash V^{\trop}(F)$ and vertices of $\DS(F)$. Inclusions are inverted in this correspondence.\\
	We say that $n$ tropical hypersurfaces $V_1,\dots,V_n$ in $N_\R$ intersect \emph{tropically transversely} at $p\in V_1\cap \dots\cap V_n$ if the point $p$ is an isolated point belonging to the interior  of a top dimensional face of $V_i$ for every $i=1.\dots,n$. Every transversal intersection point corresponds to a parallelepiped in the dual subdivision of $V_1\cup\dots\cup V_n$.
	\begin{definition}
		For a tropically transversal intersection point $p$ of $n$ tropical hypersurfaces $V_1,\dots,V_n$, its \emph{tropical intersection multiplicity} $\mult_p(V_1,\dots,V_n)$ is defined to be the sign-less volume of the parallelepiped dual to $p$.
	\end{definition}
	In the case that $V_i = V^{\trop}(\trop f_i)$, this multiplicity is exactly the number of points in $V(f_1)\cap\dots\cap V(f_n)$ that tropicalize to the tropical point (see e.g.\@ \cite{OssermanPayne}).\\
	In order to state the tropical BKK theorem, we first need the following definition.
	\begin{definition}
		Given polytopes $\Delta_1,\dots,\Delta_n$, we define the \emph{mixed volume}
		$$\MV(\Delta_1,\dots,\Delta_n)\coloneqq \text{coefficient of }\lambda_1\cdots\lambda_n \text{ in }\vol (\lambda_1\Delta_1+ \dots +\lambda_n\Delta_n).$$
	\end{definition}
	\begin{example}
		In the two-dimensional case this boils down to
		$$\MV(\Delta_1,\Delta_2) = \Area(\Delta_1+\Delta_2)-\Area(\Delta_1)-\Area(\Delta_2).$$
	\end{example}
	\begin{theorem}[Tropical BKK theorem, see e.g.\@ \cite{KazKhoEst21}]
		\label{thm:tropBKK}
		Let $F_1,\dots,F_n$ be tropical Laurent polynomials in $n$ variables with Newton polytopes $\Delta_1,\dots,\Delta_n$ and $V_i = V^{\trop}(F_i)$ for $i=1,\dots, n$ such that the intersections are tropical transversely. Then
		$$\sum_{p\in V_1\cap\dots\cap V_n}\mult_p(V_1,\dots,V_n) = \MV(\Delta_1,\dots,\Delta_n).$$
	\end{theorem}
	\begin{example}
		\label{ex:tropBKK}
		Consider the two-dimensional examples in Figure \ref{fig:tropBKK}. The curves of the same  color have the same Newton polytope but do not have the same dual subdivision. The mixed volume $\MV(\Delta_1,\Delta_2)$ equals the white area, i.e. the sum of the intersection multiplicities. In this example, this is $6$.
		\begin{figure}
		\begin{center}
			\pgfmathsetmacro{\zoom}{0.5}
			\pgfmathsetmacro{\eps}{0.28}
			\begin{tikzpicture}
				\begin{scope}[shift={(0,0)},scale={\zoom}, every node/.style={scale=1}]
					\node[blue] at (-5,-1.65) {$\Gamma_1$};
					\node[red] at (-4.5,-5) {$\Gamma_2$};
					
					\draw[blue] (1.0,2.0) -- (-1.0,-2.0);
					\draw[blue] (1.0,2.0) -- (8.0,6.666666666666666);
					\node[blue] at (-4.5,2.3) {2};
					\draw[blue] (1.0,2.0) -- (-6.0,2.0);
					\node[blue] at (-1.25,-5.0) {2};
					\draw[blue] (-1.0,-2.0) -- (-1.0,-6.0);
					\draw[blue] (-1.0,-2.0) -- (-6.0,-2.0);
					
					\draw[red] (0.0,-1.0) -- (-3.0,-1.0);
					\draw[red] (0.0,-1.0) -- (8.0,7.0);
					\draw[red] (0.0,-1.0) -- (0.0,-6.0);
					\draw[red] (-3.0,-1.0) -- (-5.0,-3.0);
					\draw[red] (-3.0,-1.0) -- (-3.0,7.0);
					\draw[red] (-5.0,-3.0) -- (-5.0,-6.0);
					\draw[red] (-5.0,-3.0) -- (-6.0,-3.0);

					\node[black] at ({-3+\eps},{1.95-\eps}) {$2$};
					\node[black] at ({7-\eps},{6+\eps}) {$1$};
					\node[black] at ({-0.5-\eps},{-1+\eps}) {$2$};
					\node[black] at ({-4+\eps},{-2.05-\eps}) {$1$};
				\end{scope}
				\begin{scope}[shift={(4.5,-2.5)},scale={2.7 * \zoom}, every node/.style={scale=1}]
					\fill[blue, opacity=0.25] (1,2) -- (1,4) -- (3,1) -- cycle;
					\fill[blue, opacity=0.25] (1,0) -- (3,0) -- (1,1) -- cycle;
					\fill[red, opacity=0.25] (0,1) -- (1,0) -- (0,0) -- cycle;
					\fill[red, opacity=0.25] (1,1) -- (0,2) -- (1,2) -- cycle;
					\fill[red, opacity=0.25] (3,0) -- (4,0) -- (3,1) -- cycle;
					
					\draw[red] (0,0) -- (1,0);
					\draw[blue] (1,0) -- (3,0);
					\draw[red] (4,0) -- (3,0);
					\draw[red] (0,0) -- (0,1);
					\draw[red] (1,0) -- (0,1);
					\draw[blue] (0,2) -- (0,1);
					\draw[red] (0,2) -- (1,1);
					\draw[blue] (1,0) -- (1,1);
					\draw[blue] (3,0) -- (1,1);
					\draw[red] (1,2) -- (1,1);
					\draw[red] (1,2) -- (0,2);
					\draw[red] (3,1) -- (4,0);
					\draw[red] (3,1) -- (3,0);
					\draw[blue] (3,1) -- (1,2);
					\draw[red] (0,2) -- (1,2);
					\draw[blue] (0,2) -- (0,4);
					\draw[blue] (1,2) -- (1,4);
					\draw[red] (0,4) -- (1,4);
					\draw[blue] (1,4) -- (3,1);
					\draw[red] (1,4) -- (2,3);
					\draw[blue] (2,3) -- (4,0);

					\filldraw[color = black, fill=black] (0,0) circle (2.1pt);
					\filldraw[color = black, fill=black] (0,1) circle (1.5pt);
					\filldraw[color = black, fill=black] (0,2) circle (1.5pt);
					\filldraw[color = black, fill=black] (0,3) circle (1.5pt);
					\filldraw[color = black, fill=black] (0,4) circle (1.5pt);
					\filldraw[color = black, fill=black] (1,0) circle (1.5pt);
					\filldraw[color = black, fill=black] (1,1) circle (1.5pt);
					\filldraw[color = black, fill=black] (1,2) circle (1.5pt);
					\filldraw[color = black, fill=black] (1,3) circle (1.5pt);
					\filldraw[color = black, fill=black] (1,4) circle (1.5pt);
					\filldraw[color = black, fill=black] (2,0) circle (1.5pt);
					\filldraw[color = black, fill=black] (2,1) circle (1.5pt);
					\filldraw[color = black, fill=black] (2,2) circle (1.5pt);
					\filldraw[color = black, fill=black] (2,3) circle (1.5pt);
					\filldraw[color = black, fill=black] (3,0) circle (1.5pt);
					\filldraw[color = black, fill=black] (3,1) circle (1.5pt);
					\filldraw[color = black, fill=black] (3,2) circle (1.5pt);
					\filldraw[color = black, fill=black] (4,0) circle (1.5pt);
					\filldraw[color = black, fill=black] (4,1) circle (1.5pt);
					
					\node[black] at (0.5,3) {$2$};
					\node[black] at (0.5,1) {$1$};
					\node[black] at (1.5,1.2) {$2$};
					\node[black] at (2.5,2) {$1$};
				\end{scope}
			\end{tikzpicture}
		\end{center}
		\begin{center}
			\pgfmathsetmacro{\zoom}{0.5}
			\pgfmathsetmacro{\eps}{0.28}
			\begin{tikzpicture}
				\begin{scope}[shift={(0,0)},scale={\zoom}, every node/.style={scale=1}]
					\node[blue] at (-5,-2.65) {$\Gamma_1$};
					\node[red] at (-4.5,-5) {$\Gamma_2$};
					
					\draw[blue] (2.0,0.0) -- (8.0,6.0);
					\draw[blue] (2.0,0.0) -- (0.0,0.0);
					\draw[blue] (2.0,0.0) -- (2.0,-6.0);
					\draw[blue] (8.0,6.0) -- (0.0,2.0);
					\draw[blue] (8.0,6.0) -- (9.0,6.666666666666666);
					\draw[blue] (0.0,0.0) -- (0.0,2.0);
					\draw[blue] (0.0,0.0) -- (-3.0,-3.0);
					\node[blue] at (-3.0,2.3) {2};
					\draw[blue] (0.0,2.0) -- (-6.0,2.0);
					\draw[blue] (-3.0,-3.0) -- (-3.0,-6.0);
					\draw[blue] (-3.0,-3.0) -- (-6.0,-3.0);
					
					\draw[red] (-2.0,-1.0) -- (-4.0,-1.0);
					\draw[red] (-2.0,-1.0) -- (6.0,7.0);
					\draw[red] (-2.0,-1.0) -- (-2.0,-6.0);
					\draw[red] (-4.0,-1.0) -- (-5.0,0.0);
					\draw[red] (-4.0,-1.0) -- (-4.0,-6.0);
					\draw[red] (-5.0,0.0) -- (-5.0,7.0);
					\draw[red] (-5.0,0.0) -- (-6.0,0.0);
					
					\node[black] at ({-5-\eps},{2+\eps}) {$2$};
					\node[black] at ({2-\eps},{3+\eps}) {$1$};
					\node[black] at ({0-\eps},{1+\eps}) {$1$};
					\node[black] at ({-2+\eps},{-2-\eps}) {$1$};
					\node[black] at ({-4+\eps},{-3.1-\eps}) {$1$};
				\end{scope}
				\begin{scope}[shift={(4.5,-2.5)},scale={2.7 * \zoom}, every node/.style={scale=1}]
					\fill[blue, opacity=0.25] (1,2) -- (2,2) -- (1,4) -- cycle;
					\fill[blue, opacity=0.25] (1,0) -- (2,0) -- (1,1) -- cycle;
					\fill[blue, opacity=0.25] (2,1) -- (3,1) -- (2,3) -- (4,0) --(3,0) -- cycle;
					\fill[red, opacity=0.25] (0,1) -- (0,2) -- (1,2) --(2,1) -- cycle;
					
					\draw[red] (0,0) -- (1,0);
					\draw[blue] (2,0) -- (1,0);
					\draw[red] (2,0) -- (3,0);
					\draw[blue] (3,0) -- (4,0);
					\draw[blue] (0,0) -- (0,1);
					\draw[red] (0,1) -- (0,2);
					\draw[blue] (0,4) -- (0,2);
					\draw[blue] (1,1) -- (1,0);
					\draw[red] (0,1) -- (1,2);
					\draw[red] (1,1) -- (1,2);
					\draw[red] (0,1) -- (1,1);
					\draw[red] (0,2) -- (1,2);
					\draw[red] (0,4) -- (1,4);
					\draw[blue] (1,2) -- (1,4);
					\draw[red] (1,4) -- (2,3);
					\draw[blue] (2,3) -- (4,0);
					\draw[blue] (4,0) -- (3,1);
					\draw[blue] (2,2) -- (1,4);
					\draw[blue] (2,2) -- (1,2);
					\draw[blue] (3,1) -- (2,3);
					\draw[red] (3,1) -- (2,2);
					\draw[blue] (1,1) -- (2,0);
					\draw[red] (1,1) -- (2,1);
					\draw[red] (1,2) -- (2,1);
					\draw[blue] (2,1) -- (3,0);
					\draw[blue] (2,1) -- (3,1);
					\draw[blue] (3,0) -- (3,1);
					
					\filldraw[color = black, fill=black] (0,0) circle (2.1pt);
					\filldraw[color = black, fill=black] (0,1) circle (1.5pt);
					\filldraw[color = black, fill=black] (0,2) circle (1.5pt);
					\filldraw[color = black, fill=black] (0,3) circle (1.5pt);
					\filldraw[color = black, fill=black] (0,4) circle (1.5pt);
					\filldraw[color = black, fill=black] (1,0) circle (1.5pt);
					\filldraw[color = black, fill=black] (1,1) circle (1.5pt);
					\filldraw[color = black, fill=black] (1,2) circle (1.5pt);
					\filldraw[color = black, fill=black] (1,3) circle (1.5pt);
					\filldraw[color = black, fill=black] (1,4) circle (1.5pt);
					\filldraw[color = black, fill=black] (2,0) circle (1.5pt);
					\filldraw[color = black, fill=black] (2,1) circle (1.5pt);
					\filldraw[color = black, fill=black] (2,2) circle (1.5pt);
					\filldraw[color = black, fill=black] (2,3) circle (1.5pt);
					\filldraw[color = black, fill=black] (3,0) circle (1.5pt);
					\filldraw[color = black, fill=black] (3,1) circle (1.5pt);
					\filldraw[color = black, fill=black] (3,2) circle (1.5pt);
					\filldraw[color = black, fill=black] (4,0) circle (1.5pt);
					\filldraw[color = black, fill=black] (4,1) circle (1.5pt);
					
					\node[black] at (0.5,3) {$2$};
					\node[black] at (0.5,0.5) {$1$};
					\node[black] at (2,0.5) {$1$};
					\node[black] at (2,1.5) {$1$};
					\node[black] at (2,2.5) {$1$};
				\end{scope}
			\end{tikzpicture}
		\end{center}
		\caption{The curves from Example \ref{ex:tropBKK}. The black numbers denote the intersection multiplicity / area.}
		\label{fig:tropBKK}
		\end{figure}
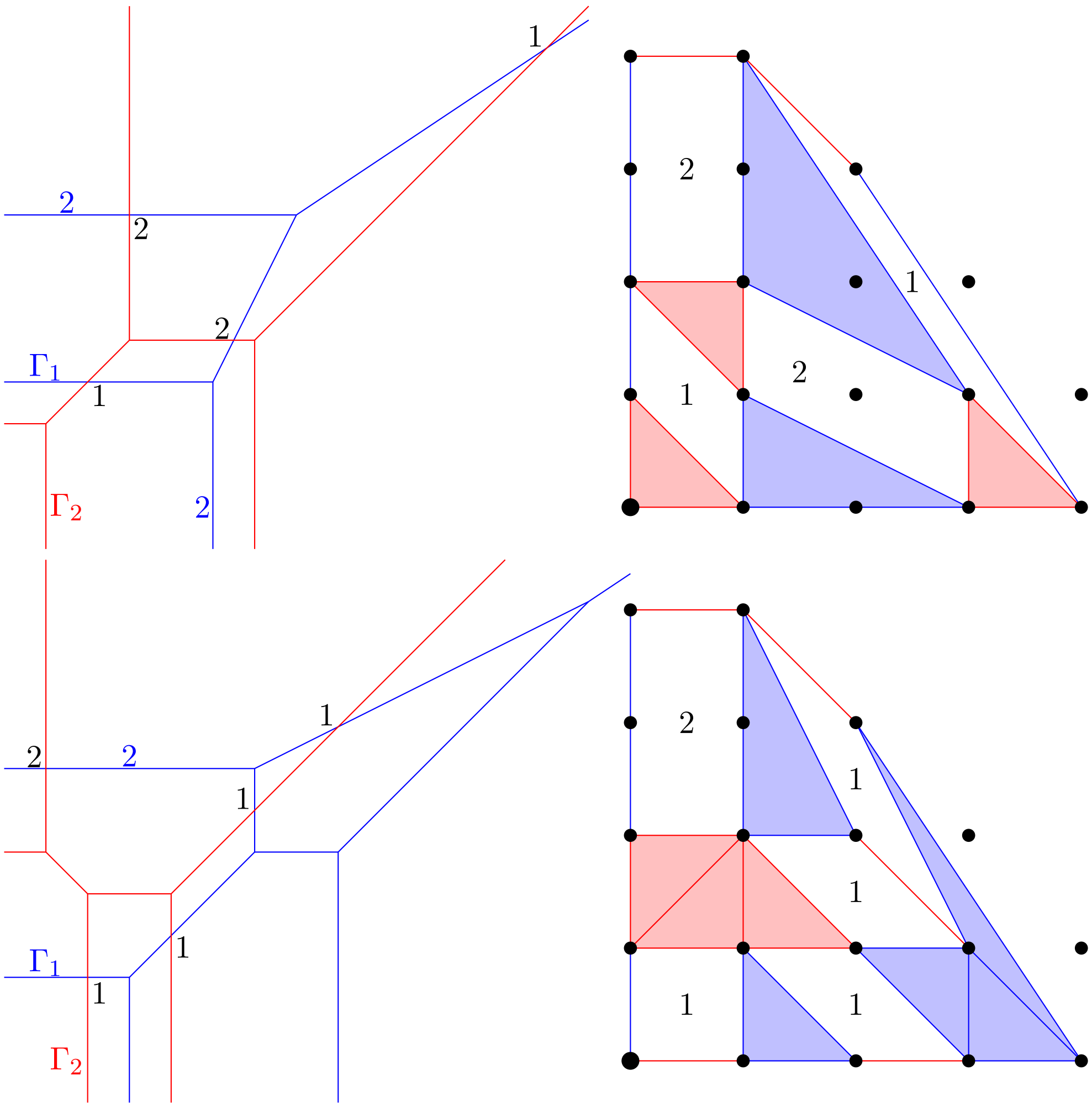
	\end{example}
\subsection{Quadratically enriched tropical geometry}
	From now on, let $\kk$ be a, not necessarily algebraically closed, field of characteristic zero or big, that is, larger than the sum of the degrees of the polynomials involved, and $\KK$ still be its field of Puiseux series.
	As mentioned in the introduction, tropical geometry can be used to solve enumerative problems, and more recently, it has also been applied to enumerative problems over not necessarily algebraically closed fields. We give a quick introduction; for more details we again refer to \cite{JPPbezout} and \cite{pauli2026pcmilecturenotesmotivic}.\\
	Interestingly, the approach is not to restrict tropicalization to the field over which one wants to solve the problem, but rather decorate the tropical side with more information.
	\begin{definition}
		A \emph{(quadratically) enriched tropical (Laurent) polynomial} is a tuple $\widetilde{F} = (\bigoplus_{i=0}^k a_i \odot\chi^{m_i},\alpha)$ where $\alpha\colon\{m_0,\dots,m_k\}\to \kkq$ is a map. We denote by $F$ the underlying tropical polynomial $\bigoplus_{i=0}^k a_i \odot\chi^{m_i}$.\\
		As a shorthand we will write $\widetilde{F} = \bigoplus_{i=0}^k \alpha_ia_i \odot\chi^{m_i}$.
	\end{definition}
	\begin{example}
		\label{ex:enrichedPolynomial}
		Let $\kk = \Q_p$ with $p = 3$. Then $\widetilde{F} = [1]3\odot x^{2}\oplus [p]4\odot xy\oplus [-1]3\odot x\oplus [-p]1\odot y\oplus [1](-2)$ is an enriched tropical polynomial.
	\end{example}
	\begin{definition}
		Given a Laurent polynomial $f = \sum_{i=0}^{k}a_i\chi^{m_i}\in\KK[M]$, we define its \emph{(quadratically) enriched tropicalization} to be $\widetilde\trop f = \bigoplus_{i=0}^k [a_i]\nu (a_i) \odot\chi^m$.
	\end{definition}
	\begin{remark}
		Here, the ``extra coefficient'' in the enriched polynomial is to be understood as remembering additional information and not as ``things one can add''.\\
		In hyperfield theory, one can make sense of enriched tropicalization as a map into tropical extensions of hyperfields, see \cite{smith2024matroidstropicalextensionstracts}. In our case, the hyperfield in question would be $\kkqo[\kk][\R]$.\\
		In \cite{JPPbezout}, instead of enriched tropical polynomials, the authors define enriched Viro polynomials, that is, a preimage under tropicalization together with a map $\alpha$ as above. The stored information is the same as in our definition, however, our definition allows us later to further ``stay on the tropical side''.
	\end{remark}
	\begin{remark}
		\label{rem:whykphase}
		Notice that multiplying a polynomial in $\KK[x_1,\dots,x_n]$ by a constant $c\in \KK$ does not change its zero set, but does change the enrichment of its tropicalization by the factor $[c]\in\kkq$. Consequently, the enriched tropicalization of the set $V(f)$ is only well-defined up to multiplication with elements of $\kkq$. In particular, the polynomial-way to enrich hypersurfaces carries too much information. In Section \ref{section:kphase} we will see a different way to enrich hypersurfaces that eliminates this ambiguity.
	\end{remark}
	\begin{definition}
		A \emph{(quadratically) enriched tropical hypersurface} $\widetilde{V} = (V,(\alpha_i)_{i=0}^k)$ is a tropical hypersurface $V$ in $N_\R$ together with a coefficient $\alpha_i\in\kkq$ assigned to each connected component of $N_\R\backslash V$.\\
		Given an enriched tropical polynomial $\widetilde{F}= \bigoplus_{i=0}^k \alpha_i a_i \odot\chi^{m_i}$, we define $\widetilde{V}^{\trop}(\widetilde{F})$ to be the enriched tropical hypersurface with $V = V^{\trop}(F)$ and to each connected component we assign the $\alpha_i$ for which $a_i\odot\chi^{m_i}$ is maximal on this component.
	\end{definition}
	Equivalently, we can assign the quadratic coefficients to the lattice points in the dual subdivision of the Newton polytope.
	\begin{example}
		\label{ex:pictureOfEnrCurve}
		For the enriched tropical polynomial from Example \ref{ex:enrichedPolynomial}, we get the enriched tropical curve seen in Figure \ref{fig:oneEnrCurve}.
		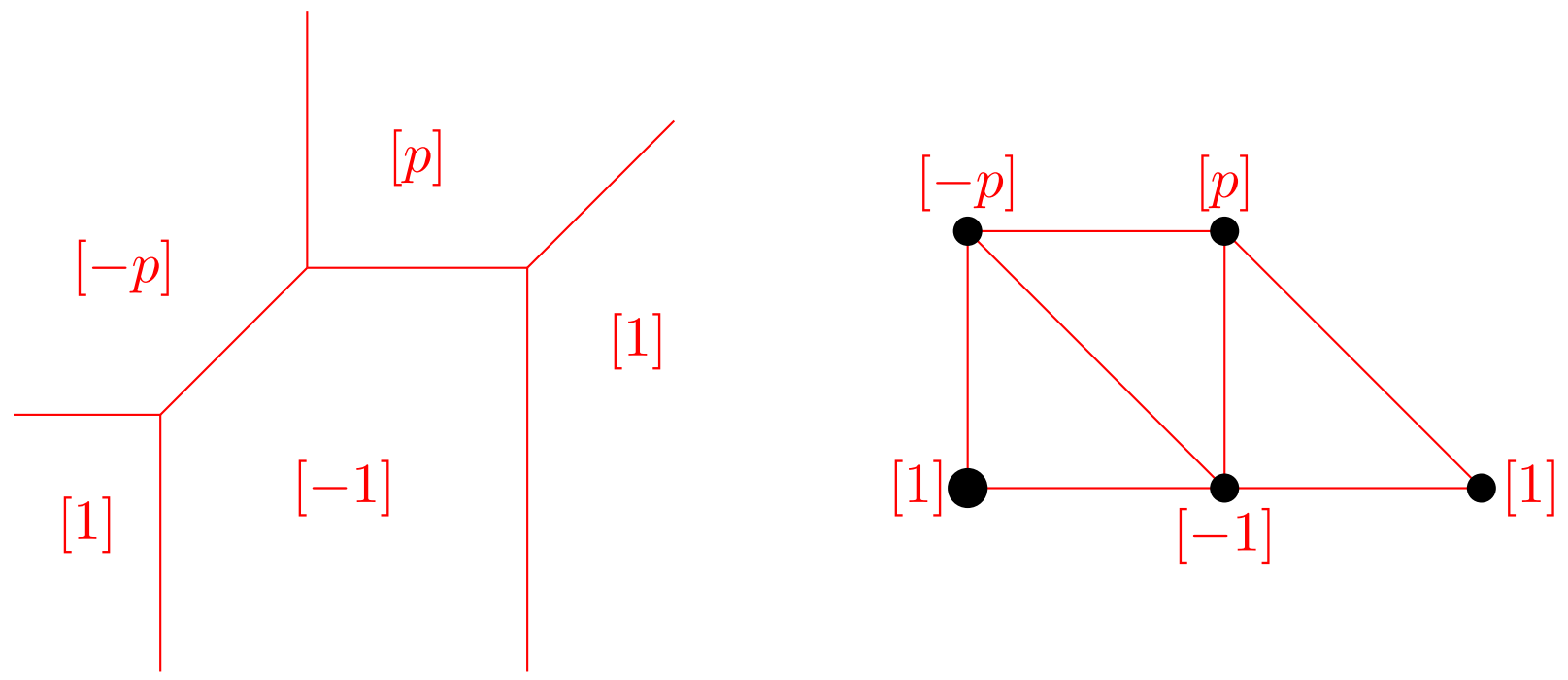
\begin{figure}
		\begin{center}
		\pgfmathsetmacro{\zoom}{0.5}
		\begin{tikzpicture}
			\begin{scope}[shift={(0,0)},scale={\zoom}, every node/.style={scale=1}]
				\draw[red] (0.0,-1.0) -- (-3.0,-1.0);
				\draw[red] (0.0,-1.0) -- (2.0,1.0);
				\draw[red] (0.0,-1.0) -- (0.0,-6.5);
				\draw[red] (-3.0,-1.0) -- (-5.0,-3.0);
				\draw[red] (-3.0,-1.0) -- (-3.0,2.5);
				\draw[red] (-5.0,-3.0) -- (-5.0,-6.5);
				\draw[red] (-5.0,-3.0) -- (-7.0,-3.0);
				
				\node[red] at (-6,-4.5) {$[1]$};
				\node[red] at (-2.5,-4) {$[-1]$};
				\node[red] at (1.5,-2) {$[1]$};
				\node[red] at (-5.5,-1) {$[-p]$};
				\node[red] at (-1.5,0.5) {$[p]$};				
			\end{scope}
			\begin{scope}[shift={(3,-2)},scale={3.5 * \zoom}, every node/.style={scale=1}]
				\draw[red] (1,0) -- (1,1);
				\draw[red] (1,1) -- (2,0);
				\draw[red] (1,0) -- (2,0);
				\draw[red] (1,0) -- (0,1);
				\draw[red] (0,1) -- (1,1);
				\draw[red] (1,0) -- (0,0);
				\draw[red] (0,1) -- (0,0);
				\filldraw[color = black, fill=black] (0,0) circle (2.1pt);
				\filldraw[color = black, fill=black] (0,1) circle (1.5pt);
				\filldraw[color = black, fill=black] (1,0) circle (1.5pt);
				\filldraw[color = black, fill=black] (1,1) circle (1.5pt);
				\filldraw[color = black, fill=black] (2,0) circle (1.5pt);
				
				\node[red,anchor=east] at (0,0) {$[1]$};
				\node[red,anchor=north] at (1,0) {$[-1]$};
				\node[red,anchor=west] at (2,0) {$[1]$};
				\node[red,anchor=south] at (0,1) {$[-p]$};
				\node[red,anchor=south] at (1,1) {$[p]$};				
			\end{scope}
		\end{tikzpicture}
		\caption{The enriched curve from Example \ref{ex:pictureOfEnrCurve}.}
		\label{fig:oneEnrCurve}
		\end{center}
		\end{figure}
	\end{example}
	After choosing coordinates, we can, just as in the classical case, define the intersection multiplicity to be the number of preimages under the valuation map -- except that we now count with an $\Al$-multiplicity. We will only consider intersections in the torus, so by choosing the chart $U_{\sigma_0}$ we can avoid all sign problems caused by the twisted coordinates of Section \ref{section:ExplicitRelOrientOnTorVar}. By Lemma \ref{lem:degAsTrace}, we get the following definition.
	\begin{definition}
		\label{def:enrTropIntMult}
		Let $\widetilde{V}_1,\dots,\widetilde{V}_n$ be $n$ enriched tropical hypersurfaces in $\R^n$ with defining polynomials $\widetilde{F}_i = \widetilde\trop f_i$ where $f_1,\dots,f_n\in \KK[x_1,\dots,x_n]$. Let $p$ be an intersection point $p\in V_1\cap \dots \cap V_n$ at which the $V_i$ intersect tropically transversely. Assume that there are exactly $s$ closed points $z^{(1)},\dots,z^{(s)}$ in $\G_m^n$ such that $\nu (z^{(i)}) = p$. We define the \emph{enriched intersection multiplicity} of $\widetilde{V}_1,\dots,\widetilde{V}_n$ at $p$ to be 
		$$\widetilde{\mult}(\widetilde{V}_1,\dots,\widetilde{V}_n)\coloneqq \sum_{i=1}^{s}\Tr_{\kappa(z^{(i)}) / \KK}\left(\gw{\det \Jac(f_1,\dots,f_n)(z^{(i)})}\right)\in\GW(\KK) \cong \GWk.$$
	\end{definition}
	Taking the rank of this enriched multiplicity recovers the classical one. \cite{JPPbezout} prove a combinatorial formula for it, which implies in particular, that the definition above only depends on the enriched tropical hypersurfaces and not on the chosen preimages of the polynomials. To give this formula, we need some preparations.\\
	Let $p$ be an intersection point as in the definition above and let $P$ be the parallelepiped in the dual subdivision of $V_1\cup\dots\cup V_n$ dual to $p$. For $i=1,\dots,n$, let $I^i,J^i$ be the exponents of the monomials of $F_i$ that are maximal at $p$ and $\alpha_i,\beta_i$ be their enrichments. If $v$ is a corner vertex of $P$, it can uniquely be expressed as
	$$v = v_\circ-\sum_{i\in A} (I^i-J^i) \eqqcolon  v_A$$
	for some subset $A \subseteq\{1,\dots,n\}$, where $v_\circ = \sum_{i=1}^{n} I^i$. It can be shown (\cite[Lemma 3.8]{JPPbezout}), that the coefficient  of $v_A$ equals
	$$\alpha_{v_A} = \prod_{i\not\in  A}\alpha_i\cdot \prod_{i\in A}\beta_i.$$
	We define the sign of the vertex $v_A$ with respect to the parallelepiped $P$ to be 
	$$\epsilon_P(v_A)\coloneqq  (-1)^{\#A}\cdot \operatorname{sign}(\det (I^1-J^1,\dots,I^n-J^n)).$$
	\begin{remark}
		Geometrically, the  sign of the  vertex $v$ with respect to the parallelepiped $P$ is the sign of the determinant of the edges of $P$ adjacent to $v$ oriented outwards to $v$.
	\end{remark}
	\begin{theorem}[{\cite[Theorem 4.7]{JPPbezout}}]
		\label{thm:combFormulaEnrTropIntersect}
		Assume that $\kk$ is a field of characteristic $0$ or characteristic bigger than the sum of the diameters of the Newton polytopes of the enriched polynomials $\widetilde{F}_i$ of $\widetilde{V}_i$ for $i=1,\dots,n$. Further, assume $\operatorname{char}\kk \not=2.$ Let $p$ be a tropically transverse intersection point of $\widetilde{V}_1,\dots,\widetilde{V}_n$. Let $P$ be the parallelepiped in the dual subdivision of $\widetilde{V}_1\cup\dots\cup\widetilde{V}_n$ corresponding to $p$ and let $v_1,\dots ,v_q$ be the odd corner vertices of $P$. Assume the non-enriched tropical intersection multiplicity equals $m$.\\
		Then
		$$\widetilde{\mult}_p(\widetilde{V}_1,\dots,\widetilde{V}_n) = \sum_{l=1}^{q}\langle \epsilon_P(v_l)\alpha_{v_l}\rangle + \frac{m-q}{2}h\in\GW(\kk),$$
		where $\alpha_{v_l}$ is the coefficient of the odd vertex $v_l$ in the dual subdivision of $\widetilde{V}_1\cup\dots\cup\widetilde{V}_n$ for $l = 1,\dots, q.$
	\end{theorem}
	We will need this proposition later. 
	\begin{proposition}[{\cite[Proposition 5.1]{JPPbezout}}]
		\label{prop:Prop5.1JPP}
		Let $V_1,\dots,V_n$ be $n$ tropical hypersurfaces in $\R^n$ intersecting tropically transversal with associated Newton polytopes $\Delta_1,\dots,\Delta_n$. Let $v$ be a lattice point in the interior of the Minkowski sum $\Delta_1+\dots+\Delta_n$ and let $P_v$ denote the set of parallelepipeds in the dual subdivision of $V_1\cup\dots\cup V_n$ dual to some intersection point of $V_1,\dots,V_n$such that $v$ is a corner vertex of $P$.\\ Then the cardinality of $P_v$ is even and there are equally many parallelepipeds $P$ in $P_v$ such that $\epsilon_P(v)$ is  positive as there are with negative sign.
	\end{proposition}
	This statement can easily be verified for the intersecting curves in Figure \ref{fig:tropBKK}. All ingredients for the following theorem were already in \cite{JPPbezout} and its proof is similar to theirs of their Corollary 5.2.
	\begin{theorem}[Enriched tropical BKK theorem]
		\label{thm:enrtropBKK}
		Let $\widetilde{F}_1,\dots,\widetilde{F}_n$ be enriched tropical polynomials in $n$ variables with Newton polytopes $\Delta_1,\dots,\Delta_n$ and $\widetilde{V}_i = \widetilde{V}^{\trop}(\widetilde{F}_i)$ for $i=1,\dots,n$ such that the intersections are tropically transversely. Further, let there be no odd lattice points on the boundary of $\Delta_1+\dots+\Delta_n$. Then
		$$\sum_{p\in V_1\cap\dots\cap V_n}\widetilde{\mult}(\widetilde{V}_1,\dots,\widetilde{V}_n) = \frac{\MV(\Delta_1,\dots,\Delta_n)}{2}h\in\GWk.$$
	\end{theorem}
	\begin{proof}
		From the classical tropical BKK theorem \ref{thm:tropBKK}, we know that the sum of the multiplicities of all intersection points is equal to $\MV(\Delta_1,\dots,\Delta_n)$. Hence, this is the rank of the enriched count. Applying the combinatorial formula of Theorem \ref{thm:combFormulaEnrTropIntersect} and Proposition \ref{prop:Prop5.1JPP}, we get that 
		$$\sum_{p\in V_1\cap\dots\cap V_n} \multA (\widetilde{V}_1,\dots,\widetilde{V}_n) = \sum_{v} \left(\frac{N(v)}{2}\langle\alpha_v\rangle+\frac{N(v)}{2}\langle-\alpha_v\rangle\right) \in \W(\kk),$$
		where the sum on the right is taken over the odd vertices in the interior of $\Delta_1+\dots+\Delta_n$ and $N(v)$ is the number of parallelepipeds in the dual subdivision of $\widetilde{V}_1\cup\dots\cup \widetilde{V}_n$ corresponding to an intersection point and $\alpha_v$ the coefficient of $v$ in $\widetilde{V}_1\cup\dots\cup \widetilde{V}_n$.\\
		Since $\langle a \rangle +  \langle -a \rangle = 0$ in $\W(\kk)$, we get that the right hand side is zero in $\W(\kk)$, and consequentially a multiple of $h$ in $\GWk$.\\
		As an element in $\GWk$ is fully determined by its rank and image in $\W(\kk)$, we are done.
	\end{proof}
	\begin{example}
		\label{ex:tropEnrBKK}
		In Figure \ref{fig:tropEnrBKK} we see two enriched tropical curves intersecting transversely. In the two-dimensional case, the sign at an odd vertex $v$ can be computed as follows. While walking in $P$ anticlockwise around the vertex, if the edge you start at is dual to an edge of $\Gamma_1$ then the sign is $+1$. Else, it is $-1$.\\
		In the example, we have multiplicities 
		\begin{itemize}
			\item $\widetilde{\mult}_{u_1}(\Gamma_1,\Gamma_2)=h$,
			\item $\widetilde{\mult}_{u_2}(\Gamma_1,\Gamma_2)=\gw{p}$,
			\item $\widetilde{\mult}_{u_3}(\Gamma_1,\Gamma_2)=\gw{-p}$,
			\item $\widetilde{\mult}_{u_4}(\Gamma_1,\Gamma_2)=\gw{1}$, and
			\item $\widetilde{\mult}_{u_5}(\Gamma_1,\Gamma_2)=\gw{-1}$.
		\end{itemize}
		 Since $\gw{a} + \gw{-a} =  h$ for all $a\in \kk^\times$, see Remark \ref{rem:a-aish}, the multiplicities sum up to $3h$, as expected.
	\begin{figure}
	\begin{center}
		\pgfmathsetmacro{\zoom}{0.5}
		\pgfmathsetmacro{\ips}{0.2}
		\pgfmathsetmacro{\eps}{0.28}
		\begin{tikzpicture}
			\begin{scope}[shift={(0,0)},scale={\zoom}, every node/.style={scale=1}]
				\node[blue] at (2.7,-5) {$\Gamma_1$};
				\node[red] at (-5.5,6) {$\Gamma_2$};
				
				\draw[blue] (2.0,0.0) -- (8.0,6.0);
				\draw[blue] (2.0,0.0) -- (0.0,0.0);
				\draw[blue] (2.0,0.0) -- (2.0,-6.0);
				\draw[blue] (8.0,6.0) -- (0.0,2.0);
				\draw[blue] (8.0,6.0) -- (9.0,6.666666666666666);
				\draw[blue] (0.0,0.0) -- (0.0,2.0);
				\draw[blue] (0.0,0.0) -- (-3.0,-3.0);
				\node[blue] at (-3.0,2.3) {2};
				\draw[blue] (0.0,2.0) -- (-6.0,2.0);
				\draw[blue] (-3.0,-3.0) -- (-3.0,-6.0);
				\draw[blue] (-3.0,-3.0) -- (-6.0,-3.0);
				
				\draw[red] (-2.0,-1.0) -- (-4.0,-1.0);
				\draw[red] (-2.0,-1.0) -- (6.0,7.0);
				\draw[red] (-2.0,-1.0) -- (-2.0,-6.0);
				\draw[red] (-4.0,-1.0) -- (-5.0,0.0);
				\draw[red] (-4.0,-1.0) -- (-4.0,-6.0);
				\draw[red] (-5.0,0.0) -- (-5.0,7.0);
				\draw[red] (-5.0,0.0) -- (-6.0,0.0);
				
				\node[red] at (-6,-1.5) {$[1]$};
				\node[red] at (-3,-2) {$[-1]$};
				\node[red] at (1.5,-2) {$[1]$};
				\node[red] at (-6,3) {$[-p]$};
				\node[red] at (-2,5) {$[p]$};	
				
				\node[blue] at (-5,-5) {$[p]$};
				\node[blue] at (-0.5,-4) {$[1]$};
				\node[blue] at (5,-1) {$[-1]$};
				\node[blue] at (-3,0.5) {$[1]$};				
				\node[blue] at (1.5,1.5) {$[p]$};				
				\node[blue] at (0,4.5) {$[1]$};
				
				\filldraw[color = black, fill=black] (-5,2) circle (2pt);
				\filldraw[color = black, fill=black] (2,3) circle (2pt);
				\filldraw[color = black, fill=black] (0,1) circle (2pt);
				\filldraw[color = black, fill=black] (-4,-3) circle (2pt);
				\filldraw[color = black, fill=black] (-2,-2) circle (2pt);
				\node[black] at ({-5-\eps},{2+\eps}) {$u_1$};
				\node[black] at ({2-\eps},{3+\eps}) {$u_2$};
				\node[black] at ({0-\eps},{1+\eps}) {$u_3$};
				\node[black] at ({-2+\eps},{-2-\eps}) {$u_5$};
				\node[black] at ({-4+\eps},{-3.1-\eps}) {$u_4$};
				
			\end{scope}
			\begin{scope}[shift={(4.75,-3)},scale={3.5 * \zoom}, every node/.style={scale=1}]
				\fill[blue, opacity=0.25] (1,2) -- (2,2) -- (1,4) -- cycle;
				\fill[blue, opacity=0.25] (1,0) -- (2,0) -- (1,1) -- cycle;
				\fill[blue, opacity=0.25] (2,1) -- (3,1) -- (2,3) -- (4,0) --(3,0) -- cycle;
				\fill[red, opacity=0.25] (0,1) -- (0,2) -- (1,2) --(2,1) -- cycle;
				
				\draw[red] (0,0) -- (1,0);
				\draw[blue] (2,0) -- (1,0);
				\draw[red] (2,0) -- (3,0);
				\draw[blue] (3,0) -- (4,0);
				\draw[blue] (0,0) -- (0,1);
				\draw[red] (0,1) -- (0,2);
				\draw[blue] (0,4) -- (0,2);
				\draw[blue] (1,1) -- (1,0);
				\draw[red] (0,1) -- (1,2);
				\draw[red] (1,1) -- (1,2);
				\draw[red] (0,1) -- (1,1);
				\draw[red] (0,2) -- (1,2);
				\draw[red] (0,4) -- (1,4);
				\draw[blue] (1,2) -- (1,4);
				\draw[red] (1,4) -- (2,3);
				\draw[blue] (2,3) -- (4,0);
				\draw[blue] (4,0) -- (3,1);
				\draw[blue] (2,2) -- (1,4);
				\draw[blue] (2,2) -- (1,2);
				\draw[blue] (3,1) -- (2,3);
				\draw[red] (3,1) -- (2,2);
				\draw[blue] (1,1) -- (2,0);
				\draw[red] (1,1) -- (2,1);
				\draw[red] (1,2) -- (2,1);
				\draw[blue] (2,1) -- (3,0);
				\draw[blue] (2,1) -- (3,1);
				\draw[blue] (3,0) -- (3,1);
				
				\filldraw[color = black, fill=black] (0,0) circle (2.1pt);
				\filldraw[color = black, fill=black] (0,1) circle (1.5pt);
				\filldraw[color = black, fill=black] (0,2) circle (1.5pt);
				\filldraw[color = black, fill=black] (0,3) circle (1.5pt);
				\filldraw[color = black, fill=black] (0,4) circle (1.5pt);
				\filldraw[color = black, fill=black] (1,0) circle (1.5pt);
				\filldraw[color=black, fill=black] (1,1) ++(-1.5pt,-1.5pt) rectangle +(3pt,3pt);
				\filldraw[color = black, fill=black] (1,2) circle (1.5pt);
				\filldraw[color=black, fill=black] (1,3) ++(-1.5pt,-1.5pt) rectangle +(3pt,3pt);
				\filldraw[color = black, fill=black] (1,4) circle (1.5pt);
				\filldraw[color = black, fill=black] (2,0) circle (1.5pt);
				\filldraw[color = black, fill=black] (2,1) circle (1.5pt);
				\filldraw[color = black, fill=black] (2,2) circle (1.5pt);
				\filldraw[color = black, fill=black] (2,3) circle (1.5pt);
				\filldraw[color = black, fill=black] (3,0) circle (1.5pt);
				\filldraw[color=black, fill=black] (3,1) ++(-1.5pt,-1.5pt) rectangle +(3pt,3pt);
				\filldraw[color = black, fill=black] (3,2) circle (1.5pt);
				\filldraw[color = black, fill=black] (4,0) circle (1.5pt);
				\filldraw[color = black, fill=black] (4,1) circle (1.5pt);
				
				\node[black] at ({-\ips},{-\ips}) {$[p]$};
				\node[black] at (1,{-\ips}) {$[-p]$};
				\node[black] at (2,{-\ips}) {$[-1]$};
				\node[black] at (3,{-\ips}) {$[1]$};
				\node[black] at (4,{-\ips}) {$[-1]$};
				
				\node[black] at ({-\ips},1) {$[1]$};
				\node[black] at ({1-\ips},{1-\ips}) {$[-1]$};
				\node[black] at ({2-\ips},{1-\ips}) {$[1]$};
				\node[black] at ({3-\ips},{1-\ips}) {$[p]$};
				
				\node[black] at ({-0.075-\ips},2) {$[-p]$};
				\node[black] at ({1-\ips},{2+\ips}) {$[p]$};
				\node[black] at ({2-\ips},{2-\ips}) {$[1]$};
				
				\node[black] at ({2-\ips},{3.075-\ips}) {$[1]$};
				
				\node[black] at ({-0.075-\ips},4) {$[-p]$};
				\node[black] at ({1-\ips},{4-\ips}) {$[p]$};
				
				\node[black] at (0.5,3) {$u_1$};
				\node[black] at (0.5,0.5) {$u_4$};
				\node[black] at (2,0.5) {$u_5$};
				\node[black] at (2,1.5) {$u_3$};
				\node[black] at (2,2.5) {$u_2$};
			\end{scope}
		\end{tikzpicture}
	\end{center}
	\caption{Two intersecting enriched curves mentioned in Example \ref{ex:tropEnrBKK}. The enrichments of the curve $\Gamma_1\cup\Gamma_2$ are noted in black in the dual subdivision. Odd vertices are squares instead of circles.}
	\label{fig:tropEnrBKK}
	\end{figure}
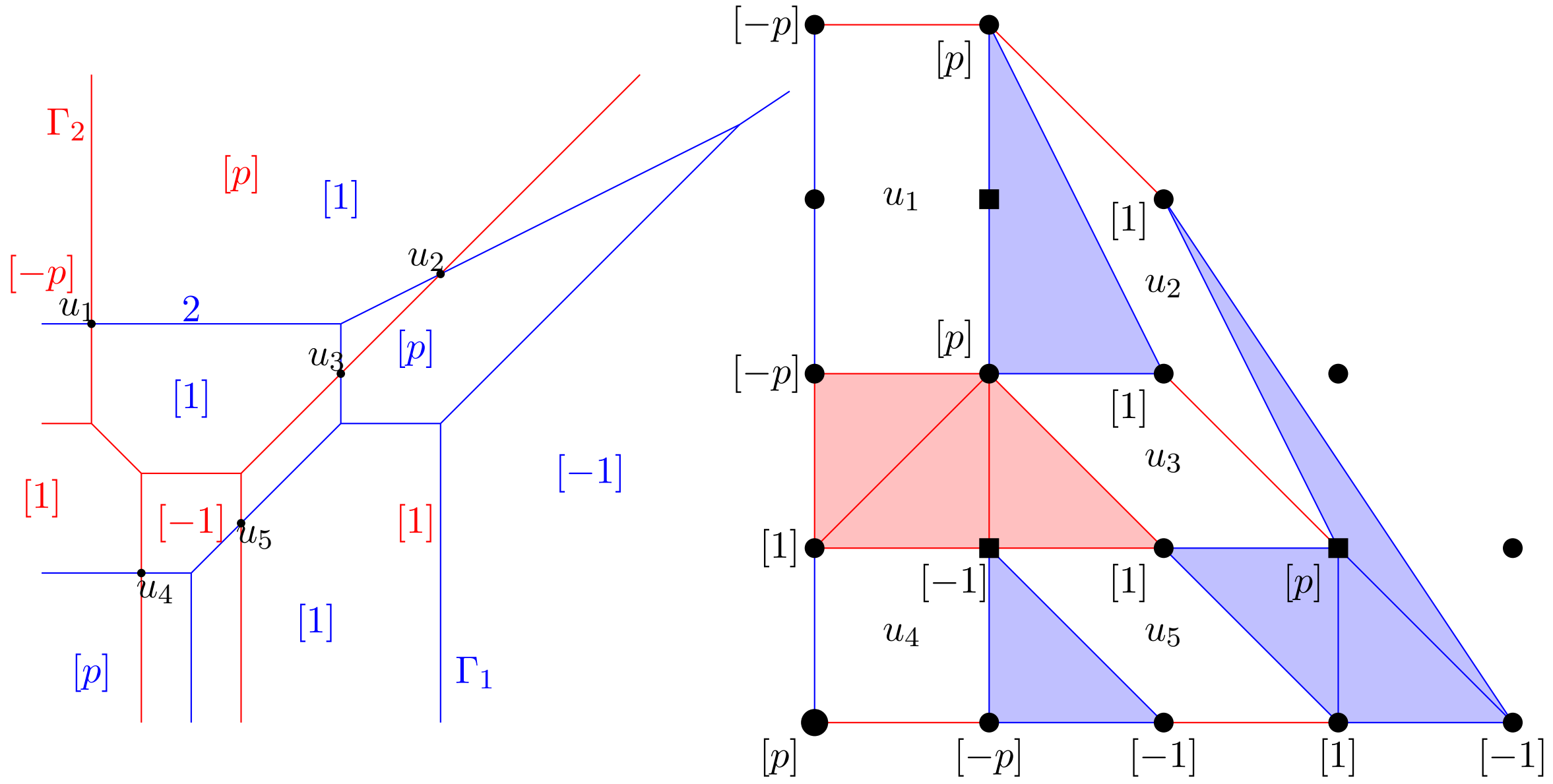
	\end{example}
\subsection{Combinatorial orientability in tropical geometry}
\label{subsection:CombOrientInTropGeo}
	Consider a toric vector bundle $\mcOX(D)$ on a toric variety $\XSig$, and a general global section $f = \sum_{m\in  P_D\cap M}a_m\chi^m$ of it. If $P_D$ is a lattice polytope, it agrees with the Newton polytope of $\trop f$.\\
	Recall that, given a balanced polyhedral complex $Y$ of codimension one with positive weights, there exists a tropical polynomial $F$ such that $Y = V^{\trop}(F)$, see e.g.\@ \cite[Theorem 2.4.10]{MikhalkinRau_TropicalGeometry}. This polynomial is unique up to multiplication with a tropical monomial, consequently its Newton polytope is unique up to translation. Fortunately, our definition of combinatorial orientability is invariant under translation. This makes the following well-defined.
	\begin{definition}
		Let $Y\subseteq N_\R$ be a balanced polyhedral complex of codimension one with positive weights. We say that $Y$ is \emph{combinatorially orientable} if $\NP(F)$ is combinatorially orientable for some tropical polynomial $F$ with $V^{\trop}(F) = Y$.
	\end{definition}
	\begin{theorem}
		\label{thm:tropicalizeOrient}
		Let $\XSig$ be a smooth proper toric variety and $D\in \Div_T(\XSig)$ be relatively orientable with $P_D$ a full-dimensional lattice polytope.	
		Then, given a general section $f$ of $\mcOX(D)$, the Newton polytope of its tropicalization is combinatorially orientable. Further, in each chart $U_\sigma$, its Newton polytope has no odd boundary points.\\
		In particular, \cite{JPPbezout}'s Conjecture 5.14 is true if $D$ is the sum of $n$ ample divisors. 
	\end{theorem}
	\begin{proof}
		This is a direct consequence of Theorems \ref{thm:RelOrientImpliesCombOrient} and \ref{thm:orientAndOddBoundary}.
	\end{proof}
	The converse is not true, see Example \ref{ex:FirstHirzNotRelOrientRedundantLine} and Example \ref{ex:NoOddNotRelOrient}.
\subsection{The $\Al$-BKK theorem}
\label{subsection:A1BKKthm}
	We re-prove the $\Al$-BKK theorem by using tropical geometry and the results obtained above. As stated in the introduction, we do not claim originality of this result, it follows e.g.\@ from \cite[Proposition 19]{SrinivasanWickelgren}, which is based on work in \cite{LevineAspects}, together with the classical result. For our proof, we need the following observation.
	\begin{lemma}
		\label{lem:orientationOverZ}
		Let $\Sigma$ be a smooth complete fan and let $X_{\Sigma,\Z}$ be the
		associated toric scheme over $\Spec\Z$. If $E_1,\dots,E_n\in\Div_T(X_\Sigma)$
		are such that $V_\Z = \mathcal{O}(E_1)\oplus\dots\oplus\mathcal{O}(E_n)$ is
		relatively orientable, then the relative orientation and the Nisnevich coordinates of
		Section \ref{section:ExplicitRelOrientOnTorVar} are defined over $\Z$, and
		their base change to any field $\kk$ is the relative orientation and the
		Nisnevich coordinates constructed over $\kk$.
	\end{lemma}
	\begin{proof}
		All ingredients of the construction -- the Cartier data $b_\sigma^{D}$,
		$b_\sigma^\omega$, $b_\sigma^E$, the monomials $\bfx^{m}$ and the signs
		$d^\sigma_{\sigma_0}\in\{\pm1\}$ -- are defined over $\Z$, and the computation
		of the transition functions is the same.
	\end{proof}
	\begin{theorem}[$\Al$-BKK theorem]
		\label{thm:A1BKK}
		Let $\Sigma$ be a smooth complete fan of dimension $n$ and let
		$E_1,\dots,E_n\in\Div_T(X_\Sigma)$ be nef such that
		$V = \mathcal{O}(E_1)\oplus\dots\oplus\mathcal{O}(E_n)$ is relatively
		orientable. Then for every field $\kk$
		$$n^{\Al}(V_\kk) = \frac{\MV(P_{E_1},\dots,P_{E_n})}{2}\,h\in\GW(\kk).$$
		In particular, if $\kk$ is infinite, then for a general section $s$ of $V_\kk$
		all zeros are simple and lie in the torus, and
		$\sum_{z} \ind_z s = \frac{\MV(P_{E_1},\dots,P_{E_n})}{2}\,h$.
	\end{theorem}
	\begin{proof}
		We use the same approach as in \cite[Remark 5.4.]{JPPbezout}. By Lemma \ref{lem:orientationOverZ}, $V_\Z\to X_{\Sigma,\Z}$ is a relatively
		oriented vector bundle on a smooth proper scheme over $\Z$, so its Euler number
		$n^{\Al}(V_\Z)\in\GW(\Z)$ is defined, and for every field $\kk$ its image under
		$\GW(\Z)\to\GW(\kk)$ is $n^{\Al}(V_\kk)$ \cite{BachmannWickelgren}.
		The map $\GW(\Z)\to\GW(\R)$ is an isomorphism, as
		$\GW(\Z)\cong\Z\gw{1}\oplus\Z\gw{-1}$, further $\GW(\R)\to\GW(\KKR)$ for $\KKR\coloneqq\R\{\{t\}\}$ is an isomorphism, too, by \cite{MarkwigPayneShaw}. It therefore suffices
		to show $n^{\Al}(V_{\KKR}) = \frac{\MV(P_{E_1},\dots,P_{E_n})}{2}h$.\\
		Since each $E_i$ is nef, each $P_{E_i}$ is a lattice polytope
		(Theorem~\ref{thm:BasepointFreeCharacterizations}). For $i=1,\dots,n$ choose $h_i\colon P_{E_i}\cap M\to\Q$ and coefficients $c_{i,m}\in\R^\times$ and set
		$$f_i = \sum_{m\in P_{E_i}\cap M} c_{i,m}\,t^{-h_i(m)}\chi^m
		\in\Gamma(X_{\Sigma,\KKR},\mathcal{O}(E_i)).$$
		Then $\trop f_i = \bigoplus_m h_i(m)\odot\chi^m$ and $\NP(\trop f_i) = P_{E_i}$.
		For generic $h_i$, the tropical hypersurfaces $V_i = V^{\trop}(\trop f_i)$ intersect
		tropically transversely and for generic $c_{i,m}$ the
		section $f=(f_1,\dots,f_n)$ has only simple zeros, all of them in the torus
		$\Spec\KKR[M]$.
		As $n^{\Al}(V_\KK)$ does not depend on the section \cite{BachmannWickelgren},
		we have $n^{\Al}(V_{\KKR})=\sum_z\ind_z f$.\\
		We use the chart $U_{\sigma_0}$ to avoid sign issues. Let $\widetilde V_i = \widetilde V^{\trop}(\widetilde\trop f_i^{\sigma_0})$. Notice that $f^{\sigma_0}_i$ has Newton polytope $\NP_{\sigma_0}(E_i)$. By
		Lemma \ref{lem:degAsTrace}, $\ind_z f = \Tr_{\kappa(z)/\KK}
		\gw{\det\Jac(f^{\sigma_0}_1,\dots,f^{\sigma_0}_n)(z)}$. Every zero $z$
		tropicalizes to a point of
		$V_1\cap \dots \cap V_n$,
		so by Definition \ref{def:enrTropIntMult}
		$$n^{\Al}(V_{\KKR}) = \sum_{p\in V_1\cap \dots\cap V_n} \widetilde{\mult}_p(\widetilde V_1,\dots,\widetilde V_n).$$
		Let $D=\sum_i E_i$. As the $E_i$ are nef, $P_D = \sum_i P_{E_i}$ and its a lattice polytope, further 		$$\sum_i\NP_{\sigma_0}(E_i) = \psi^D_{\sigma_0}\Big(\sum_i P_{E_i}\Big)
		= \NP_{\sigma_0}(D),$$
		using $b^D_{\sigma_0}=\sum_i b^{E_i}_{\sigma_0}$. Since $D$ is relatively
		orientable, $\NP_{\sigma_0}(D)$ has no odd boundary points by
		Theorem \ref{thm:orientAndOddBoundary}. As $\operatorname{char}\KKR=0$, the
		enriched tropical BKK theorem (Theorem~\ref{thm:enrtropBKK}) gives
		$$n^{\Al}(V_\KK) = \frac{\MV(\NP_{\sigma_0}(E_1),\dots,\NP_{\sigma_0}(E_n))}{2}\,h
		= \frac{\MV(P_{E_1},\dots,P_{E_n})}{2}\,h,$$
		where the last equality holds as $\psi^{E_i}_{\sigma_0}$ are affine maps with
		the same unimodular linear part.\\
		Last, if $\kk$ is infinite, the argument above over $\kk$ shows that a general
		section of $V_\kk$ has only simple zeros, all in the torus, and
		$n^{\Al}(V_\kk)=\sum_z\ind_z s$ by \cite{BachmannWickelgren}.
	\end{proof}
\section{Phase structures}
\label{section:kphase}
	The definition of enriched tropical hypersurfaces is motivated by Viro's patchworking over $\R$, which uses distributions of signs - this is the same thing as a quadratic enrichment over the real numbers in our language. Viro's patchworking,see e.g. \cite{viro2000dequantization}, is a tool to construct curves with certain topological properties. In higher codimension, \cite{RRS_patchwork} use real phase structures for patchworking. A real phase stores the extra information on the surface itself and not in their complement, so it can be defined in higher codimension, too. The idea behind the definition  of real  phase structures is to remember the orthant from which the segment of the tropical variety came under tropicalization, instead of the coefficient of some defining polynomial.\\
	As we have seen earlier, enriching the hypersurfaces with the coefficients of the monomials involves a choice; instead of a polynomial $f$ one could also tropicalize the polynomial $cf$, $c$ a non-zero constant, and would end up with an enriched hypersurface with different enrichment. As mentioned in  Remark \ref{rem:whykphase}, for $\kk$-phases, this ambiguity does not exist. Additionally, the definition of $\kk$-phases on hypersurfaces is a first step towards a possible generalization of tropical enriched intersection multiplicity to  not-$0$-dimensional intersections. Further, it could help, analogously to the real case, to generalize the quadratically enriched refined patchworking theorem from \cite{JPPcorresp}, which was inspired by \cite{shustin2005patchworkingsingularalgebraiccurves}. Most of the content of this section were part of the author's Master's thesis.
\subsection{Real phase structures}
	We give the definition following \cite{RRS_patchwork}. For that, let $\operatorname{Aff}_d((\Z/2\Z)^n)$ denote the set of affine subspaces of $(\Z/2\Z)^n$. Given an affine subspace $A$, denote its tangent space by $T(A)$. We say that $A\in\operatorname{Aff}_d((\Z/2\Z)^n)$ is parallel to a rational cone $\sigma$ generated by $\rho_1,\dots,\rho_s$, if $T(A)$ equals the reduction modulo $2$ of the integer points of $T(\sigma)\coloneqq\spann_{\R}\{\rho_1,\dots,\rho_s\}$. 
	\begin{definition}
		Let $E$ be a set. A collection $C$ of subsets of $E$ is called an \emph{even covering}, if every element in $E$ the union is contained in an even number of elements of $C$.
	\end{definition}
	Recall that for a polyhedral complex $Y$ of pure dimension $d$, a face of $Y$ is called a facet if it is of dimension $d$.
	\begin{definition}
		\label{def:realPhaseRRS}
		Let $Y$ be a polyhedral complex of pure dimension $d$ in $\R^n$. A \emph{real phase structure} $\mce$ on $Y$ is a map 
		$$\mce\colon \text{Facets of }Y\to \operatorname{Aff}_d((\Z/2\Z)^n)$$
		such that
		\begin{enumerate}[(i)]
			\item for every facet $\sigma$ of $Y$, the affine subspace $\mce(\sigma)$ is parallel to $\sigma$.
			\item for every codimension one face $\tau$ of $\mcx$ with adjacent facets $\sigma_1,\dots,\sigma_k$, the sets $\mce(\sigma_1),\dots,\mce(\sigma_k)$ are an even covering.
		\end{enumerate}
	\end{definition}
	This definition of real phase structure uses the fact that $(\kkq[\R],\cdot)\to(\Z/2\Z,+);\ [-1]\mapsto 1,[1]\mapsto  0$ is an isomorphism. So, intuitively, a point in a facet $\sigma$ with $(0,\dots,0)\in\mce(\sigma)$ has a preimage under tropicalization that lies in the positive orthant.
\subsection{General phase structures}
	In the following, we will look at equations in the group $\kkq$. Obviously, $\kkq$ is an abelian group of exponent $2$. In particular, we can consider all exponents modulo $2$. For reading convenience, we introduce the notation $-\alpha \coloneqq [-1]\cdot \alpha$. Note that $[-1] = [1]$ is possible. We use multi-index notation.\\
	As we want to follow the same path of thought as for the real phase, we have to think about the following question: If $\alpha_Ix^I$ and $\alpha_Jx^J$, where $I,J\in\N^n,\ \alpha_I,\alpha_J\in\kkq$, are the two enriched monomials for which the minimum is attained, what do we need $[x_1],\dots,[x_n]\in\kkq$ to be to allow the existence of some $a_I\in\alpha_{I},\ a_J\in\alpha_{J},\ y_i\in[x_i],\ i\in\{1,\dots,n\}$ such that $a_Iy^I + a_Jy^J = 0$?\\
	Equivalently, we can ask for the solutions to the equation 
	\begin{equation*}
		[a_I][x^I] = [-a_J][x^J] \Leftrightarrow [x^{I-J\mod{2}}] = -\alpha_I\alpha_J
	\end{equation*}
	in $(\kkq)^n$. Consequently, we want to assign the set $\{x\in(\kkq)^n\mid x^{I-J} = -\alpha_I\alpha_J\}$ to the facet defined by the enriched monomials $\alpha_Ix^I$ and $\alpha_Jx^J$.\\
	Let $\tau$ be a codimension one face of the nonsingular tropical hypersurface defined by the polynomial $F = \bigoplus_{I\subseteq\N_0^n} c_I\odot x^I$ with adjacent facets 
	$\sigma_1 = \{x\in\R^n\mid F(x) = c_{I_1}\odot x^{I_1} = c_{J_1}\odot x^{J_1}\},\dots,\sigma_3 = \{x\in\R^n\mid F(x) = c_{I_3}\odot x^{I_3} = c_{J_3}\odot x^{J_3}\}$. Because of the balancing condition, we know that $\sum_{i=1}^{3}(I_i-J_i) = 0$.\\
	Hence, 
	$$[1] = x^{\sum_{i=1}^{3}(I_i-J_i)} = \prod_{i=1}^{3} x^{I_i-J_i} = [-1]^3\prod_{i=1}^{3}\alpha_{I_i}\alpha_{J_i}.$$
	In the nonsingular case, the even covering condition requires the sets associated to the facets to have empty intersection. This inspires the following definition.
	\begin{definition}
		A \emph{$\kk$-phase} on a tropical hypersurface $\mcx$ in $\R^n$ is a map $\varphi\colon \text{Facets of }Y\to \kkq$ such that for each codimension $1$ face with adjacent facets $\sigma_1,\dots,\sigma_k$ it is $\varphi(\sigma_1)\cdots\varphi(\sigma_k) = [1]$.
	\end{definition}
	To get a set of ``signs'' out of the associated element, we need to solve the equation $x^{I-J} = -\varphi(\sigma)$ in $(\kkq)^n$. Note that this set might be empty or the whole group. If it is not, it is a coset in the group $(\kkq)^n$ containing $(\kkq)^{n-1}$ elements. In the real case this gives an affine space parallel to the facet. This observation allows to define the concept of such cosets being parallel to a vector subspace.\\
	Note that $\R$-phase, i.e. $\kk$-phase for the real numbers, and real phase do not agree if the tropical hypersurface contains a facet of weight $2$, as the solution set might be empty. However, we get the following theorem that ensures the even covering property, (ii) in Definition \ref{def:realPhaseRRS}.
	\begin{theorem}
		\label{thm:evenCovering}
		For a tropical hypersurface, the even covering property, (ii) in Definition \ref{def:realPhaseRRS}, is equivalent to the property in the definition of $\R$-phase that $\prod_{i=1}^{k}\varphi(\sigma_i) = [1]$.
	\end{theorem}
	\begin{proof}
		Let $y\in(\kkq[\R])^n$ be a solution to the equations $x^{I_i-J_i} = -\varphi(\sigma_i)$ for w.l.o.g.\@ $i=1,\dots,l$. For $i>l$ follows $x^{I_i-J_i} = \varphi(\sigma_i)$, as $|\kkq[\R]| = 2$. Then 
		$$1 = y^0 = y^{I_1-J_1}\cdots y^{I_l-J_l}y^{I_{l+1}-J_{l+1}}\cdots y^{I_k-J_k} = [-1]^l \prod_{i=1}^{k}\varphi(\sigma_i) = [-1]^l,$$
		which is true if, and only if, $l$ is even.
	\end{proof}
	\begin{corollary}
		For a tropical hypersurface with no even weights, real phase structures and $\R$-phase structures are in a one-to-one correspondence.
	\end{corollary}
	\begin{example}
		\label{Exm:EmptyRPhase}
		Consider the tropical curve $V^{\trop}(0\oplus x \oplus y^2)$, see Figure \ref{Fig:RPhaseEmpty}. For the $\R$-phase we see in the figure, the set assigned to $\sigma_1$ is $\{(x,y) \in \kkq[\R] \mid [1] = x^0y^2 = [-1]\} = \emptyset$.\\
		For $\sigma_3$ and $\sigma_2$ the set is $\{(x,y) \in \kkq[\R] \mid x^1y^2 = x = [-1]\} = \{(x,y) \in \kkq[\R] \mid x^1y^0 = [-1]\} = \{(-1,1),(-1,-1)\}$. They had to agree, by the even covering property.
		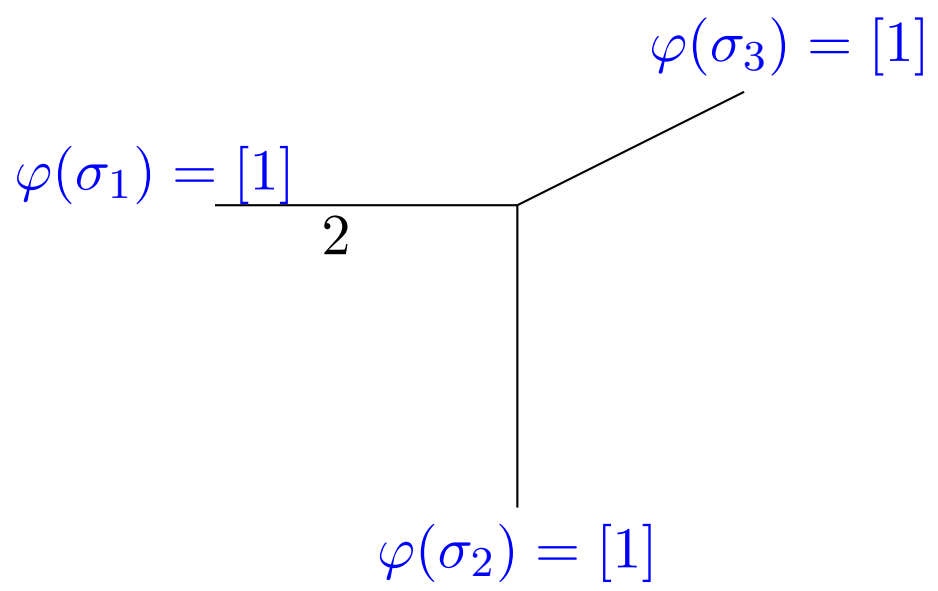
\begin{figure}[htb]
			\centering
			\begin{tikzpicture}
				\draw [black] (1.5,0.75) -- (0,0);
				\draw [black] (-2,0) -- (0,0);
				\draw [black] (0,-2) -- (0,0);
				\node[black] at (-1.2,-0.2) {$2$};
				\node[blue] at (-2.4,0.2) {$\varphi(\sigma_1) = [1]$};
				\node[blue] at (0,-2.3) {$\varphi(\sigma_2) = [1]$};
				\node[blue] at (1.8,1.05) {$\varphi(\sigma_3) = [1]$};
			\end{tikzpicture}
			\caption{Example \ref{Exm:EmptyRPhase}'s tropical curve in $\R^2$ with $\R$-phase.}
			\label{Fig:RPhaseEmpty}
		\end{figure}
	\end{example}
	\begin{example}
	\label{Exm:NotSupSet}
		Let $\kk = \Q_p$ with $p=3$. Recall that in this case $\kkq = \{1,-1,p,-p\}$. Consider the tropical line seen in Figure \ref{Fig:TropLine} with $\kk$-phase defined by $\varphi(\sigma_1) = [p],\varphi(\sigma_2)=[-p]$ and $\varphi(\sigma_3)=[-1]$. In this case we get
		\begin{itemize}
			\item $\{(x,y)\in (\kkq)^2\mid y = [-p]\} = \{([1],[-p]),([-1],[-p]),([p],[-p]),([-p],[-p])\}$ for $\sigma_1$
			\item $\{(x,y)\in (\kkq)^2\mid x = [p]\} = \{([p],[1]),([p],[-1]),([p],[p]),([p],[-p])\}$ for $\sigma_2$, and
			\item $\{(x,y)\in (\kkq)^2\mid xy = [1]\} = \{([1],[1]),([-1],[-1]),([p],[p]),([-p],[-p])\}$ for $\sigma_3$.
		\end{itemize}
		These sets obviously do not form an even covering. The reason the proof of Theorem \ref{thm:evenCovering} fails in this case is that $\kkq$ has more than two elements.
	\end{example}
	\begin{example}
		\label{Exm:NonEmptyTrippleIntersection}
		Let $\kk=\Z/5\Z$. Then $\kkq$ has two elements, namely $[1] = \{1,4\}$ and $[2] = \{2,3\}$, in particular $[1] = [-1]$. Consider, again, the tropical line seen in Figure \ref{Fig:TropLine} with $\kk$-phase assigning to each ray the element $[1]$. Then 
		\begin{itemize}
			\item for $\sigma_1$ we get $\{(x,y)\in (\kkq)^2\mid y = [1]\} = \{([1],[1]),([2],[1])\}$,
			\item for $\sigma_2$ we get $\{(x,y)\in (\kkq)^2\mid x = [1]\} = \{([1],[1]),([1],[2])\}$, and
			\item for $\sigma_3$ we get $\{(x,y)\in (\kkq)^2\mid xy = [1]\} = \{([1],[1]),([2],[2])\}$.
		\end{itemize}
		In particular, the even covering condition is not satisfied for any occurring element. The proof of Theorem \ref{thm:evenCovering} does not apply in this example, as here $[1] = [-1]$ in $\kkq$.
	\end{example}
	\begin{remark}
		Interestingly, the set of $\kk$-phase structures on a given tropical hypersurface with point-wise multiplication is a group (which is abelian and of exponent $2$).
	\end{remark}
	\begin{remark}
		In \cite{RRS_orientedMatroid}, the authors establish a one-to-one correspondence between orientations of a matroid and real phase structures on its matroid fan (also known as its Bergman fan). Their constructions described in sections 3.1 and 3.4 can easily be adapted to our case of $\kk$-phase and matroids over the hyperfield $\kkqo$. Notice however, that, as we only consider codimension one fans, the underlying matroid for our case always has to be $U_{n-1,n}$. In this case, one can show that the one-to-one correspondence is true.\\
		In \cite{smith2024matroidstropicalextensionstracts} the author constructs enriched tropical linear spaces out of a given matroid over certain hyperfields. On matroid fans, the spaces obtained there agree with our construction, with hyperfield $\kkqo$.\\
		It would be desirable to get a matroid-independent $\kk$-phase definition on fans of higher codimension and generalize \cite{RRS_orientedMatroid}'s result. This is work in progress.
	\end{remark}
\subsection{Phase structures and enrichments}
	In this subsection, we establish a correspondence between $\kk$-phases and enrichments analogously to the real case treated in \cite{renaudineauRealPhaseAndDistribution}.
	Let $\Delta$ be a dual subdivision to a tropical hypersurface $\mcx$ with Newton polytope $P$. We denote by $\Delta^0$ the set of vertices, by $\Delta^1$ the set of edges, and by $\Delta^2$ the set of two-dimensional cells of $\Delta$. For an edge connecting $v,w\in\Delta^0$ we write $vw\in\Delta^1$. Denote by $\Delta^{(2)}$ the union of the cells of $\Delta$ with dimension lesser or equal to $2$.\\
	With this notation, an enrichment of $\mcx$ is the same as a map $\alpha\colon\Delta^0\to\kkq$ and a $\kk$-phase a map $\varphi\colon\Delta^1\to \kkq$ such that $\varphi(e_1)\cdots\varphi(e_k)$ for $e_1,\dots,e_k\in\Delta^1$ the boundary edges of some $F\in \Delta^2$.
	\begin{definition}
		We call an enrichment $\alpha$ \emph{compatible with} a $\kk$-phase $\varphi$ if for any $v\in\Delta^0$ the following condition is satisfied:\\
		For any vertex $w\in\Delta^0$ adjacent to $v$, one has $\alpha(w) = \varphi(vw)\alpha(v)$.
	\end{definition}
	\begin{example}
		\label{Exm:TwoDimNotLineKPhase}
		See Figure \ref{Fig:TwoDimNotLineKPhase} for an example of a curve with compatible enrichment and $\kk$-phase for $\kk = \Q_p$ with $p =3$.
		\begin{figure}[htb]
			\centering
			\begin{tikzpicture}[scale=1, every node/.style={scale=1}]
				\draw[black] (-1,0) -- (0,-1);
				\draw[black] (-1,1) -- (-1,0);
				\draw[black] (-1,1) -- (-2,1);
				\draw[black] (-1,0) -- (-2,0);
				\draw[black] (-1,1) -- (0,2);
				\draw[black] (0,-1) -- (1,-1);
				\draw[black] (0,-1) -- (0,-2);
				\draw[black] (1,-1) -- (1,-2);
				\draw[black] (1,-1) -- (2,0);
				
				\node[blue, anchor=north east] at (-0.5,-0.5) {$[p]$};
				\node[blue, anchor=east] at (-2,0) {$[-p]$};
				\node[blue, anchor=east] at (-2,1) {$[-1]$};
				\node[blue, anchor=north] at (0,-2) {$[-p]$};
				\node[blue, anchor=north] at (1,-2) {$[p]$};
				\node[blue, anchor=south west] at (2,0) {$[-p]$};
				\node[blue, anchor=south west] at (0,2) {$[1]$};
				\node[blue, anchor=south] at (0.5,-1) {$[-1]$};
				\node[blue, anchor=west] at (-1,0.5) {$[-1]$};
				
				\node[black] at (-1.5,-1.5) {$[1]$};
				\node[black] at (0.5,-1.5) {$[-p]$};
				\node[black] at (1.5,-1.5) {$[-1]$};
				\node[black] at (-1.5,0.5) {$[-p]$};
				\node[black] at (-1.5,1.5) {$[p]$};
				\node[black] at (1.5,1.5) {$[p]$};
			\end{tikzpicture}\hspace{0.5cm}
			\begin{tikzpicture}[scale=2, every node/.style={scale=1}]
				\filldraw[color = black, fill=black] (0,0) circle (1.5pt);
				\filldraw[color = black, fill=black] (0,1) circle (1.5pt);
				\filldraw[color = black, fill=black] (0,2) circle (1.5pt);
				\filldraw[color = black, fill=black] (1,0) circle (1.5pt);
				\filldraw[color = black, fill=black] (1,1) circle (1.5pt);
				\filldraw[color = black, fill=black] (2,0) circle (1.5pt);
				
				\draw[black] (0,0) -- (0,2);
				\draw[black] (0,0) -- (2,0);
				\draw[black] (2,0) -- (0,2);
				\draw[black] (0,0) -- (1,1);
				\draw[black] (1,1) -- (0,2);
				\draw[black] (2,0) -- (1,1);
				\draw[black] (1,0) -- (1,1);
				\draw[black] (0,1) -- (1,1);
				
				\node[anchor=north east] at (0,0) {$[1]$};
				\node[anchor=east] at (-0,1) {$[-t]$};
				\node[anchor=south east] at (-0,2) {$[p]$};
				\node[anchor=north] at (1,-0) {$[-p]$};
				\node[anchor=north west] at (2,-0) {$[-1]$};
				\node[anchor=south west] at (1,1) {$[p]$};
				
				\node[blue, anchor=east] at (-0,0.5) {$[-p]$};
				\node[blue, anchor=east] at (-0,1.5) {$[-1]$};
				\node[blue, anchor=north] at (0.5,-0) {$[-p]$};
				\node[blue, anchor=north] at (1.5,-0) {$[p]$};
				\node[blue, anchor=south west] at (1.5,0.5) {$[-p]$};
				\node[blue, anchor=west] at (1,0.5) {$[-1]$};
				\node[blue, anchor=south east] at (0.5,0.5) {$[p]$};
				\node[blue, anchor=south] at (0.5,1) {$[-1]$};
				\node[blue, anchor=south west] at (0.5,1.5) {$[1]$};
			\end{tikzpicture}
			\caption{A tropical hypersurfaces in $\R^2$ enriched over $\kk = \Q_p$ with $p=3$ and their duals. The black coefficients are the enrichment, the blue coefficients define a $\kk$-phase structure compatible with it.}
			\label{Fig:TwoDimNotLineKPhase}
		\end{figure}
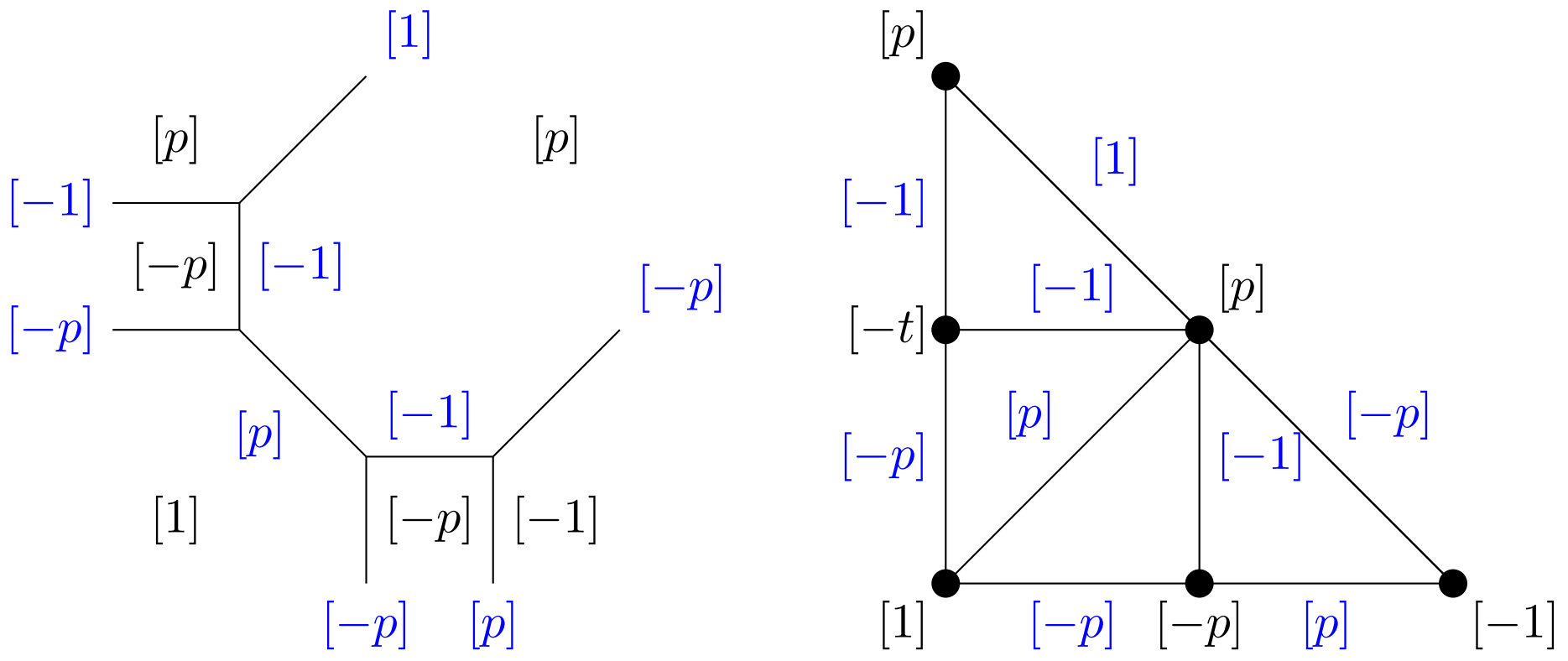
	\end{example}
	\begin{theorem}
		Let $\mcx$ be a tropical hypersurface with dual subdivision $\Delta$ of its Newton polytope $P$.\\
		For any $\kk$-phase $\varphi$ on $\mcx$, there exist exactly $|\kkq|$ enrichments compatible with $\varphi$. Reciprocally, given an enrichment $\alpha$ of $\mcx$, there exists a unique $\kk$-phase compatible with $\alpha$.
	\end{theorem}
	\begin{proof}
		Let $\varphi\colon \Delta^1\to \kkq$ be a $\kk$-phase on $\mcx$. Choose an arbitrary vertex $v_0$ of $\Delta$ and assign to it an arbitrary element $\beta\in\kkq$. Starting at $v_0$, given a vertex $v$ of $\Delta$ already equipped with an element $\alpha(v)\in\kkq$, for every adjacent vertex $w$ assign an element of $\kkq$ by using the compatibility condition, i.e.\@ $\alpha(w) = \varphi(vw)\alpha(v)$.\\
		We have to verify that, no matter which path we choose from $v_0$ to a vertex $w$, we assign to it the same element of $\kkq$. Equivalently, we have to show that for every loop $L = (e_1,\dots,e_l)$ it is $\prod L\coloneqq\prod\varphi(e_l) = [1]$. The property of $\varphi$ being a $\kk$-phase immediately implies that homotoping a loop $L$ across a two-dimensional cell does not change $\prod L$. Now, notice that the fundamental group $\pi_1(P)$ of the Newton polytope is trivial, as $P$ is a convex set. Further, $\pi_1(P) = \pi_1(\Delta^{(2)})$. Hence, given a loop $L$, we know that $\prod L = [1]$.\\	
		Given an enrichment $\alpha$, we get a $\kk$-phase by defining $\varphi\colon\Delta^1\to\kkq;\ vw\mapsto \alpha(v)\cdot\alpha(w)$: If $e_1=v_1v_2,\dots,e_k=v_kv_1\in\Delta^1$ are the boundary edges of some $F\in \Delta^2$, it is
		$$\varphi(e_1)\cdots\varphi(e_k) = \alpha(v_1)\alpha(v_2)\cdot \alpha(v_2)\alpha(v_3)\cdots\alpha(v_k)\alpha(v_1) = \prod_{i=1}^k\alpha(v_i)^2 = [1].$$ It is compatible with $\alpha$, as $\alpha(w) = \alpha(v)^2\alpha(w) = \varphi(vw)\alpha(v)$ for adjacent $v,w\in\Delta^0$.
	\end{proof}
	This result is not surprising: Given a hypersurface $V(f)$, we can change the defining polynomial to $af$ for any $a\in\kk^\times$ and not change the hypersurface. So we get $|\kkq|$ different possibilities for the residue classes of the coefficients of $f$ per hypersurface.
	\bibliographystyle{amsalpha}
	\bibliography{library}
\end{document}